\input harvmac
\input psfig
\def\rhob{{\rho\kern-0.465em \rho}}

\def\Z{{\bf Z}}

\def\ontopss#1#2#3#4{\raise#4ex \hbox{#1}\mkern-#3mu {#2}}

\def\nom#1#2{{#1 \atopwithdelims[] #2}}

\def\mod{~{\rm mod}~}
\def\e#1{{\bf e}_{#1}}
\def\be#1{{\bf \bar e}_{#1}}
\def\E#1#2{{\bf E}^{(t)}_{#1,#2}}
\def\bE#1#2{{\bf \bar E}^{(t)}_{#1,#2}}

\def\bEl#1#2{{\bf \bar E}_{#1,#2}}
\def\u#1{{\bf u}'_{#1}}
\def\aw#1{\mathop{\longrightarrow}^{#1}}
\def\ev#1#2{\mathop{\longrightarrow}^{#1}_{#2}}
\def\evs#1{\mathop{\Longrightarrow}^{#1}}
\def\vn{{\bf n}}
\def\vm{{\bf m}}
\def\vA{{\bf A}}
\def\vB{{\bf B}}
\def\vu{{\bf u}}
\def\vw{{\bf w}}
\def\O{{\bf O}}
\def\ve{{\bf e}}
\def\bra#1#2{{#1 \atopwithdelims\{\} #2}}
\def\kio#1{\theta(#1~{\rm odd})}
\def\kie#1{\theta(#1~{\rm even})}
\def\m{\mu}
\def\bu{{\bf\bar{u}}}

\setbox\strutbox=\hbox{\vrule height12pt depth5pt width0pt}

\def\strut{\relax\ifmmode\copy\strutbox\else\unhcopy\strutbox\fi}

\nref\rroga{L.J. Rogers, Proc. Lond. Math. Soc. 25 (1894) 318.}

\nref\rschur{I. Schur, S.--B Preuss. Akad. Wiss. Phys.--Math. Kl
(1917) 302.}

\nref\rram{L.J. Rogers and 
S. Ramanujan, Proc. Camb. Phil. Soc. 19 (1919) 214.}

\nref\rmac{P.A. MacMahon, {\it Combinatory Analysis,} vol. 2
(Cambridge University Press, Cambridge 1916).}

\nref\randd{G.E. Andrews,{\it The Theory of Partitions}, Encyclopedia
of Mathematics and its Applications, vol 2. G.--C. Rota ed. (Addison
Wesley 1976).}

\nref\rslatb{L.J. Slater, Proc. Lond. Math. Soc. (2) {\bf 54} (1951--52)
147.}

\nref\rbax{R.J. Baxter, J. Stat. Phys. {\bf 26} (1981) 427.}

\nref\rabf{G.E. Andrews, R.J. Baxter and P. J. Forrester,
J. Stat. Phys. 35 (1984) 193.}

\nref\rdate{E. Date, M. Jimbo, T. Miwa and M. Okado, 
	Phys. Rev. B{\bf 35} (1987) 2105;
E. Date, M. Jimbo, A. Kuniba, T. Miwa and M. Okado, 
	Nucl. Phys B{\bf 290} (1987) 231 and
	Adv. Stud. in Pure Math. {\bf 16} (1988) 17.}

\nref\rbpz{A.A. Belavin, A.M. Polyakov and A.B. Zamolodchikov,
J. Stat. Phys. {\bf 34} (1984) 763 and Nucl. Phys. B{\bf 241} (1984) 333.}

\nref\rff{B. Feigin and D.B. Fuchs, Funct. Anal. Appl. {\bf 17} (1983) 241.}

\nref\rlf{A. Feingold and J. Lepowsky, Adv. in Math. {\bf 29} (1978) 271.}

\nref\rlp{J. Lepowsky and M. Primc, {\it Structure of the standard
modules for the affine Lie algebra $A^{(1)}_1$}, Contemporary
Mathematics, Vol. 46 (AMS, Providence, 1985).}

\nref\rkm{R. Kedem and B.M. McCoy, J. Stat. Phys. {\bf 71} (1993) 865.}

\nref\rdkmm{S. Dasmahapatra, R. Kedem, B.M. McCoy and E. Melzer,
J. Stat. Phys. {\bf 74} (1994) 239.}

\nref\rkkmma{R. Kedem, T.R. Klassen, B.M. McCoy and E. Melzer,
Phys. Letts. B {\bf 304} (1993) 263.}

\nref\rkkmmb{R. Kedem, T.R. Klassen, B.M. McCoy and E. Melzer,
Phys. Letts. B {\bf307} (1993) 68.}

\nref\rdkkmm{S. Dasmahapatra, R. Kedem, T.R. Klassen, B.M. McCoy and
E. Melzer, Int. J. Mod. Phys. B{\bf 7} (1993) 3617.} 

\nref\rrc{A. Rocha--Caridi, in {\it Vertex Operators in Mathematics
and Physics}, ed. J. Lepowsky, S. Mandelstam and I.M. Singer (Springer
Berlin 1985).}

\nref\rand{G.E. Andrews, Proc. Nat. Sci. USA {\bf 71} (1974) 4082.}

\nref\randpac{G.E. Andrews, Pac. J. Math. 114 (1984) 267.}

\nref\rbmo{A. Berkovich, B.M. McCoy and W. Orrick,
J. Stat. Phys. {\bf 83} (1996) 795.}

\nref\rb{A. Berkovich, Nucl. Phys. B{\bf 431} (1994) 315.}

\nref\rFW{O. Foda and S.O. Warnaar, Lett. Math. Phys. {\bf 36} (1996) 145.}

\nref\rw{S.O. Warnaar, J. Stat. Phys. {\bf 82} (1996) 657.}

\nref\rfoda{O. Foda and Y-H. Quano, 
	Int. J. Mod. Phys. A 12 (1997) 1651.}

\nref\rbm{A. Berkovich and B.M. McCoy, Lett. Math. Phys. {\bf 37} (1996) 49.}

\nref\rmela{E. Melzer, Int. J. Mod. Phys. A{\bf 9} (1994) 1115.}

\nref\ranne{A. Schilling, Nucl. Phys. B{\bf 459} (1996) 393 and 
Nucl. Phys. B{\bf 467} (1996) 247.}

\nref\randb{G.E. Andrews, Scripta Math. {\bf 28} (1970) 297.}

\nref\rFQ{O. Foda and Y.-H. Quano, Int. J. Mod. Phys. A {\bf 10} (1995) 2291.}

\nref\rKirillov{A.N. Kirillov, Prog. Theor. Phys. Suppl. {\bf 118} (1995) 61.}

\nref\role{S.O. Warnaar, Commun. Math. Phys. {\bf 184} (1997) 203.}

\nref\rfb{P.J. Forrester and R.J. Baxter, J. Stat. Phys. {\bf 38} (1985) 435.}

\nref\rabbbfv{G.E. Andrews, R.J. Baxter, D.M. Bressoud, W.H. Burge,
P.J. Forrester and G. Viennot, Europ. J. Combinatorics {\bf 8} (1987) 341.}

\nref\rwp{S.O. Warnaar, P.A. Pearce, K.A. Seaton and B. Nienhuis,
J. Stat. Phys. {\bf 74} (1994) 469.}

\nref\rts{M. Takahashi and M. Suzuki, Prog. of Theo. Phys. {\bf 48}
 (1972) 2187.}

\nref\rgr{G. Gaspar and M. Rahman, {\it Basic Hypergeometric Series},
(Cambridge Univ. Press 1990), Appendix I.}

\nref\rbmb{A. Berkovich and B.M. McCoy, Int. J. of Math. and Comp.
Modeling (in press), hepth/9508110.}

\nref\rbg{ E. Baver and D. Gepner, Phys. Lett. B{\bf 372} (1996) 231.}

\nref\rbgs{A. Berkovich, C. Gomez and G. Sierra, Nucl. Phys. B{\bf 415}
(1994) 681.}

\nref\rbail{W.N. Bailey, Proc. Lond. Math. Soc. (2) 50 (1949) 1.} 
 
\nref\raab{A.K. Agarwal, G.E. Andrews and D.M. Bressoud, J. Indian
Math. Soc. 51 (1987) 57.}

\nref\rbmsb{A. Berkovich, B.M. McCoy and A. Schilling,
Physica A {\bf 228} (1996) 33.}

\nref\rsw{A. Schilling and S.O. Warnaar, Int. J. Mod. Phys. B {\bf 11}
(1997) 189 and {\it A Higher--Level Bailey Lemma: Proof and Application},
to appear in The Ramanujan Journal (q-alg/9607014).}

\nref\rbmsw{A. Berkovich, B.M. McCoy, A. Schilling and S.O Warnaar,
Bailey flows and Bose-Fermi identities for the conformal coset models
$(A_1^{(1)})_N\times(A_1^{(1)})_{N'}/(A_1^{(1)})_{N+N'},$ submitted
to Nucl. Phys. B (hep-th/9702026).}

\Title{\vbox{\baselineskip12pt\hbox{BONN-TH-96-07}
  \hbox{ITPSB 96-35}
  \hbox{q-alg 9607020}}}
  {\vbox{\centerline{Rogers-Schur-Ramanujan Type Identities for the}
\centerline{$M(p,p')$ minimal models of Conformal Field Theory}}}

  \centerline{Alexander Berkovich~\foot{berkov\_a@math.psu.edu}}

  \medskip\centerline{\sl Physikalisches Institut der}
  \centerline{\sl Rheinischen Friedrich-Wilhelms Universit{\"a}t Bonn}
  \centerline{\sl Nussallee 12}
  \centerline{\sl D-53115 Bonn, Germany}  
  \bigskip
  \centerline{and}
  \bigskip
  \centerline{ Barry~M.~McCoy~\foot{mccoy@max.physics.sunysb.edu} and
               Anne Schilling\foot{anne@insti.physics.sunysb.edu}}

  \medskip\centerline{\sl Institute for Theoretical Physics}
  \centerline{\sl State University of New York}
  \centerline{\sl Stony Brook,  NY 11794-3840}
  \bigskip
  \bigskip
  \medskip
  \centerline{\it Dedicated to the memory of Poline Gorkova} 
  \bigskip
  \Date{\hfill 07/96}
  
  \eject

\centerline{\bf Abstract}

We present and prove  Rogers--Schur--Ramanujan (Bose/Fermi) type
identities for the  Virasoro characters of the minimal model $M(p,p').$
The proof uses the continued fraction decomposition of $p'/p$
introduced by Takahashi and Suzuki for the study of the Bethe's
Ansatz equations of the XXZ model and gives a general method to
construct polynomial generalizations of the fermionic form of
the characters which satisfy the same recursion relations as the
bosonic polynomials of Forrester and Baxter. We use this method to get
fermionic representations of the characters $\chi_{r,s}^{(p,p')}$ for
many classes of $r$ and $s.$

\newsec{Introduction}
Rogers--Schur--Ramanujan type identities is the generic mathematical name
given to the identities which have been developed in the last 100 years
from the work of Rogers~\rroga, Schur~\rschur~ and Ramanujan \rram~who proved, 
among other things,
that for $a=0,1$
\eqn\rog{\eqalign{\sum_{j=0}^{\infty}{q^{j^2+aj}\over (q)_j}&=
\prod_{j=1}^{\infty}{1\over (1-q^{5j-1-a})(1-q^{5j-4+a})}\cr
&={1\over
(q)_{\infty}}\sum_{j=-\infty}^{\infty}(q^{j(10j+1+2a)}-q^{(5j+2-a)(2j+1)})\cr}}
where
\eqn\qdfn{(q)_k=\prod_{j=1}^{k}(1-q^j),~~k>0;~~~(q)_0=1.}
These identities are of the greatest importance in the theory of
partitions~\rmac--\randd~and number theory and by 
the early 50's at least 130 of them were known~\rslatb. 

The emergence of these identities in physics is much more recent,
starting with the work of Baxter~\rbax, Andrews, Baxter and
Forrester~\rabf~ and the Kyoto group \rdate~
in the 80's on the order parameters of solvable
statistical mechanical models.

An even more recent relation to physics is the application to
conformal field theories invented by Belavin, Polyakov and
Zamolodchikov~\rbpz~in 1984. Here the left hand side is obtained by
using a fermionic basis and the right hand side is obtained by using a
bosonic basis for the Fock space. The bosonic constructions are done
in a universal fashion using the methods of Feigin and
Fuchs~\rff. The construction of the fermionic basis is more
involved.  The earliest examples of such
fermionic representations are for modules of the affine 
Lie algebra~$A^{(1)}_1$~\rlf-\rlp, but the general theory of this 
application of the identities has only
been explicitly developed in the last several years~\rkm-\rdkkmm.
However some of the mathematics of these constructions is already present 
in the original identities~\rog~which are now recognized as being the
fermi/bose identities for the conformal field theory $M(2,5).$ 
Thus in some sense one might say
that the 1894 paper of Rogers~\rroga~ is one
of the first mathematical contributions to conformal field theory even
though conformal field theory as a physical theory was invented 
only in 1984 \rbpz.

The theory of bosonic representations of conformal field theory
characters is well developed. In  particular the characters of all
$M(p,p')$ minimal models are given by the Rocha-Caridi formula~\rrc
\eqn\rca{{\hat{\chi}}_{r,s}^{(p,p')}(q)=q^{\Delta_{r,s}^{(p,p')}-c/24}
B_{r,s}(q)}
where
\eqn\rocb{B_{r,s}(q)={1\over
(q)_{\infty}}\sum_{j=-\infty}^{\infty}(q^{j(jpp'+rp'-sp)}-q^{(jp'+s)(jp+r)}),}
with conformal dimensions
\eqn\dim{\Delta_{r,s}^{(p,p')}={(rp'-sp)^2-(p-p')^2\over
4pp'}~~(1\leq r\leq p-1,~1\leq s\leq p'-1),}
and central charge
\eqn\ccharge{c=1-{6(p-p')^2\over pp'}.}
$p$ and $p'$ are relatively prime and we note the symmetry property 
$B_{r,s}(q)=B_{p-r,p'-s}(q).$ It is obvious  that ~\rocb~generalizes the sum 
on the right hand side of ~\rog. 

The generalization of the $q$-series on the left hand side of ~\rog~has
a longer history. The first major advance was made in the 70's when
Andrews realized~\rand~that there were generalizations of ~\rog~ in terms
of multiple sums of the form
\eqn\multsum{\sum_{m_1,\cdots, m_{k-1}}q^{{1\over 2}{\bf m}^T {\bf Bm}
+{\bf A}^T{\bf m}}\prod_{i=1}^{k-1}{1\over (q)_{m_i}}~~~~(k\geq 2)}
where $\bf B$ is a $(k-1)$ by $(k-1)$ matrix, $\bf A$ is a 
$(k-1)$--dimensional vector and the summation variables $m_i$ run over
positive integers. These
results are now recognized as the characters of the $M(2,2k+1)$ models.
In the interpretation of ~\rkm-\rdkkmm~we say that each $m_i$
represents the number of fermionic quasi-particles of type $i$. 
A second generalization of the
form~\multsum~ is that the summation variables may obey restrictions such as
$m_i$ being even or odd or having linear combinations being congruent 
to some value $Q$ (mod $N$) (where $N$ is some integer). 
The odd/even restriction is present in the original work
of Rogers~\rroga~and the (mod $N$) restrictions were first found by
Lepowsky and Primc~\rlp.
 
Further generalizations of ~\multsum~are needed to
represent the most general character. In particular we need 
the general form of
what we call the ``fundamental fermionic form'' which was
first found in ~\rdkmm~(and generalizes the special case (5.5) of
~\randpac which in retrospect is $M(5,5k+2)$)
\eqn\femiform{f_{r,s}(\vu,q)
=\sum_{{\bf m},{\rm restr.}}q^{{1\over 2}\vm^T \vB\vm+\vA^T\vm}
\prod_{i=1}^{k-1}{(({\bf I}_{k-1}-\vB)\vm+\vu)_i\atopwithdelims[]m_i}_q}
where ${\bf I}_{k-1}$ is the $(k-1)$ by $(k-1)$ dimensional unit matrix and 
$\vu$ is a $(k-1)$--dimensional vector with components $(\vu)_i$.
(In general we impose the notation that the components of a vector $\vu$
are either denoted by $(\vu)_i$ or $u_i$. $\vu_i$ would denote a vector
labeled by $i$ and not its $i$th component).
We define the $q$-binomial coefficients for nonnegative $m$ and $n$ as
\eqn\qbindef{{m+n \atopwithdelims[] m}_q=
{m+n \atopwithdelims[] n}_q=
\cases{{(q)_{m+n}\over (q)_m(q)_n}& if $m,n\geq 0$,\cr
0&otherwise.}}
There exist generalizations of \qbindef~to negative $n$, and their use 
in the context of fermionic characters was first found in~\rbmo.
We also note that using the property
\eqn\binprop{\lim_{n\rightarrow \infty} {m+n\atopwithdelims[] m}_q=
{1\over (q)_m}}
the general form \femiform~reduces to \multsum~when $u_i\rightarrow
\infty$ for all $1\leq i\leq k-1$.
Then in terms of these fundamental fermionic forms the generic form of
the generalization of ~\rog~ is now given as the linear combination 
\eqn\eidbf{F_{r,s}(q)=\sum_{i}q^{c_i}f_{r,s}({\bf u}_i;q)=
q^{{\cal N}_{r,s}}B_{r,s}(q)}
where ${\cal N}_{r,s}$ is a normalization constant.

Character identities of this form which generalize the results of 
~\rand~were conjectured for some special cases of $M(p,p')$ in~\rkkmmb~
including $p'=p+1$. Proofs of the identities for $M(p,p+1)$ are given in
\rb--\rw. Several other special cases of $p$ and $p'$ and particular 
values of $r$ 
and $s$ are proven in~\rfoda. In a previous letter~\rbm~two of
the present authors gave results for the case of arbitrary $p$
and $p'$ for certain selected values of $r$ and $s$. Here we
generalize and prove the results of that letter. 
 
Our method of proof is to generalize the infinite series for the 
bosonic and fermionic forms of the characters in~\eidbf~to polynomials 
$B_{r,s}(L,q)$ and $F_{r,s}(L,q)$ whose order depends on an integer $L$ and 
then to prove that both $B_{r,s}(L,q)$ and $F_{r,s}(L,q)$ satisfy the same 
difference equations in $L$ with the same boundary conditions. The 
generalization from infinite series to a set of polynomials is referred to as 
``finitization''~\rmela. 

For the proof of the $L$--difference equations we utilize the technique
of telescopic expansion first introduced in ~\rb~to prove the conjecture 
of~\rkkmmb, \rmela~for the
$M(p,p+1)$ model and subsequently used to prove identities for the
$N=1$ supersymmetric model $SM(2,4\nu)$~\rbmo~and for  general
series of the A$_1^{(1)}$ coset models with integer levels~\ranne. 
 This method is the  extension to many quasi
particles of the recursive proof of ~\rog~ given by~\rschur,~\rmac~ and
~\randb. (Somewhat different methods have been used to prove 
polynomial analogues of the Andrews--Gordon
identities in \rFQ,~\rKirillov,~\role.) 
  
There are many ways to finitize the fermionic and bosonic forms of
the characters. For example the bosonic character~\rocb~has a
polynomial generalization
in terms of $q$-binomial coefficients
\eqn\polyfb{B_{r(b),s}(L,b;q)=\sum_{j=-\infty}^{\infty}
\biggl(q^{j(jpp'+r(b)p'-sp)}{L\atopwithdelims[]
{L+s-b\over 2}-jp'}_q-q^{(jp+r(b))(jp'+s)}
{L\atopwithdelims[] {L-s-b\over 2}-jp'}_q\biggr)}
where $L+s-b$ is even and  $r(b)$ is a prescribed function of $b$ with  
$1 \leq b \leq p'-1$ (see (3.4) below). This generalization first
appeared in the work of Andrews, Baxter and Forester~\rabf~for
$p'=p+1$ and for general $p~p'$ in the work of Forrester and Baxter~\rfb. 
The bosonic polynomials in 
\polyfb~have the symmetry
\eqn\blsym{B_{r(b),s}(L,b;q)=B_{p-r(b),p'-s}(L,p'-b;q).}
Thanks to~\binprop~equation~\polyfb~tends to~\rocb~as 
$L\to\infty$. Notice that in this limit the dependence on $b$ drops out
and hence for each character identity there are several different polynomial
identities with the same limit.
The polynomials~\polyfb~generalize the polynomials used by
Schur~\rschur~ in connection with difference two partitions and is the
finitization   we use in this paper. They
satisfy a simple recursion relation in $L$ and 
can be interpreted as the generating functions for partitions with
prescribed hook differences~\rabbbfv. However, there are several other
known polynomial finitizations~\rbmo, \rw, \role, \rwp~
which satisfy other $L$ difference
equations and prove to be useful in other contexts.

In this paper we will present a method which allows the construction
of fermionic polynomials which  satisfy the identities
\eqn\rfbpolyid{F_{r(b),s}(L,b;q)=q^{{\cal N}_{r(b),s}}B_{r(b),s}(L,b;q)}
for the general minimal model $M(p,p')$ in principle for all $b$ and $s$.
It is however difficult to find a notation which allows for a compact
treatment of all values of $b$ and $s$ at the same time. Consequently
even though our methods in this paper are general, we  
present results only for certain classes of $b$ and $s.$
However, we emphasize that all cases can be treated by the same methods. 
Additional results will be presented elsewhere.

The polynomials appearing on the left hand side of~\rfbpolyid~
generalize the polynomials
originally used by MacMahon~\rmac~ in his analysis 
of~\rog.
In contrast to the bosonic polynomials~\polyfb~the form of
the fermionic polynomials depends on the values of $b$ and $s.$ This
is because the fermionic polynomials depend on the continued fraction
decomposition of $p'/p.$ We will present the formalism of this
decomposition in sec. 2 and defer the presentation of our results to
sec. 3. 

The proof of our results for $p'>2p$ is given in secs. 4--11. The
case $p<p'<2p$ is obtained from the case $p'>2p$ in sec. 12 by the 
method of the dual transformation discussed in ~\rbm. We close in sec.
13 with a discussion of several ways our results can be extended and
with an interpretation of our polynomial identities in terms of new
Bailey pairs. 

\newsec{Summary of the formalism of Takahashi and Suzuki}
 
We begin with the observation that the models
$M(p,p')$ are obtained as a reduction of the XXZ spin chain
\eqn\hamxxz{H_{XXZ}=-\sum_{k}(\sigma^x_k\sigma^x_{k+1}+
\sigma_k^{y}\sigma^{y}_{k+1}+\Delta\sigma^z_{k}\sigma^z_{k+1})} 
where $\sigma_k^i~(i=x,y,z)$ are the Pauli spin matrices and
\eqn\deldef{\Delta=-\cos\pi{p\over p'}.} 
Consequently we may use the results
of the classic study of the thermodynamics of the XXZ chain made by
Takahashi and Suzuki~\rts~ in 1972.
This treatment begins by introducing, for $p'>2p$, the $n+1$ integers 
$\nu_0,\nu_1,\cdots,\nu_n$ 
from the following
continued fraction decomposition of $p'/p$
\eqn\contfrac{{p'\over p}=\nu_0+1+{1\over \nu_1+{1\over
\nu_2+\cdots+{1\over \nu_n+2}}}}
where $\nu_n\geq 0$ and all other $\nu_j\geq 1.$ For the
case $p'<2p$ we replace $p$ by $p'-p.$
We say that this is a $n+1$ zone decomposition and that there
are $\nu_j$ types of quasi particles in  zone $j$. From these integers we
define (where $\m$ is an integer, $0\leq \m\leq n+1$)
\eqn\defnt{t_\m=\cases{\sum_{j=0}^{\m-1}\nu_j&for $1\leq \m \leq n+1$\cr
-1&for $\m=0$\cr}.}
We refer to $t_{n+1}$ as the number of types of quasi particles in the system.
When an index $j$ satisfies 
\eqn\zone{t_\m+1\leq j \leq t_{\m+1}+\delta_{n,\m}}
we say that the index $j$ is in the $\m^{th}$ zone and that $1+t_\m$ 
and $t_{1+\m}$ are
the boundaries of this zone. Note that by definition zone $0 ~(n)$ has
$\nu_0+1~(\nu_n+1)$ allowed values of $j$ while all other zones have
$\nu_{\m}$ allowed values of $j.$ We will sometimes refer to $j=t_{n+1}+1$
and $j=0$ as  ``virtual'' positions. We will explicitly consider below the
case $p'>2p.$ The case $p'<2p$ will be treated  separately in sec 12.

According to Takahashi and Suzuki~\rts, there are also $\nu_i$ types
of quasi particles in zones $i=0,\ldots,n-1$ for the XXZ chain. In zone $n$
there are, however,  $\nu_{n}+2$ types of quasiparticles, and in addition
there is an extra zone $n+1$ with one quasi particle in the XXZ chain.
It is the omission of the three quasi particles of zone $n$ and $n+1$
which truncates the XXZ chain to the model $M(p,p').$

{}From the $\nu_j$ we define the set of integers $y_{\m}$ recursively as
\eqn\ydef{y_{-1}=0,~~y_0=1,~~y_1=\nu_0+1,
  ~~y_{\m+1}=y_{\m-1}+(\nu_{\m}+2\delta_{\m,n})y_{\m},~~~(1\leq \m\leq n).}
and further set
\eqn\ldef{\l_j=y_{\m-1}+(j-1-t_{\m}) y_{\m}
  ~~{\rm for}~~1+t_{\m}\leq j \leq t_{\m+1}+\delta_{n,\m}.}
We then define what we call the Takahashi length 
$\{l^{(\m)}_{1+j}| 0\leq \m\leq n, j$ is in $\m$th zone $\}$ 
\eqn\stringdef{l^{(\m)}_{1+j}=\cases{j+1&for $\m=0$ 
and $0\leq j \leq t_1$\cr
y_{\m-1}+(j-t_{\m})y_{\m}& for $1\leq \m\leq n$ and $1+t_{\m}\leq
j\leq t_{1+\m}+\delta_{n,\m}$\cr}}
which differs from $l_{1+j}$ only at the boundaries $j=t_{\m}.$ Note
that $l_{1+j}^{(\m)}$ is monotonic in $j$ while $l_{1+j}$ is not.
To indicate that $j$ lies in the $\m$th zone, i.e. $1+t_{\m}\leq j\leq
t_{1+\m}$, we will in the following always write $j_{\m}$ instead of $j$. 

We also define a second set of integers $z_{\m}$ for as
\eqn\zdef{z_{-1}=0,~~z_0=1,~~z_1=\nu_1+2\delta_{1,n},~~
 z_{\m+1}=z_{\m-1}+(\nu_{\m+1}+2\delta_{\m+1,n})z_{\m},~~~(1\leq \m\leq n-1)}
and 
\eqn\ltildem{{\tilde l}^{(\m)}_{1+j_{\m}}=\cases{z_{\m-2}+(j_{\m}-t_{\m})
 z_{\m-1}&
for $1\leq \m\leq n$ and $1+t_{\m}\leq j_{\m}\leq t_{1+\m}+\delta_{n,\m}$,\cr
0&for $\m=0.$\cr}}  
We refer to $\tilde {l}_{1+j_{\m}}^{(\m)}$ as a truncated Takahashi length.
It is clear that $z_{\m}$ is obtained  from
the same set of recursion relations as the $y_{\m}$ except that zone zero
is removed in the partial fraction decomposition of $p'/p.$ The
removal of this zone zero is equivalent to considering a new XXZ chain
with an anisotropy $\Delta'=-\cos\pi\{{p'\over p}\}$ where $\{x\}$ denotes
the fractional part of $x.$
We note that the $l^{(\m)}_{1+j_{\m}}~~({\tilde l}^{(\m)}_{1+j_{\m}})$ 
are the dimensions of the unitary
representations of the quantum group $su(2)_{q_\pm}$ with $q_+=
e^{i\pi {p\over p'}}~~(q_-=e^{i\pi{ \{{p'\over p}\}}})$. 


The final result we need from~\rts~is the specialization of their
equation (1.10) to the case of the $0^{th}$ Fourier component. Then,
using the notation where the integers $n_k~~(m_k)$ with $1\leq k\leq t_{n+1}$
are the number of
particle (hole) excitations of type $k$ we find what we call
$(\vm,\vn)$ system (eqn (2.18) of ~\rbm)
\eqn\mnsys{\eqalign{n_k+m_k&={1\over 2}(m_{k-1}+m_{k+1})+
{1\over 2}\bar{u}_k~~{\rm for}~ 1\leq k
\leq t_{n+1}-1~{\rm and}~k\neq t_i,i=1,\cdots ,n\cr
n_{t_i}+m_{t_i}&={1\over 2}(m_{t_i-1}+m_{t_i}-m_{t_i+1})+{1\over
2}\bar{u}_{t_i}
,~~{\rm
for}~~i=1,\cdots ,n\cr
n_{t_{n+1}}+m_{t_{n+1}}&={1\over 2}
(m_{t_{n+1}-1}+m_{t_{n+1}}\delta_{\nu_n,0})+{1\over 2}\bar{u}_{t_{n+1}}\cr}}
where by definition $m_0=L$, all $\bar{u}_k$ are integers and all $m_k$ are
nonnegative integers. Let us emphasize here that whereas $m_k$ is always
nonnegative, $n_k$ may at times take on negative values. 
We denote~\mnsys~symbolically as
\eqn\syn{{\bf n}={\bf M} {\bf m}+{L\over 2}{\bf{\bar e}}_1+{1\over 2}{\bu}}
where we define the  $t_{n+1}$-dimensional vectors ${\bf{\bar e}}_k$ by
\eqn\iunit{({\bf{\bar e}}_k)_j=\cases{\delta_{j,k}& $1\leq k\leq t_{n+1}$\cr
0& $k=0, 1+t_{n+1}.$}}
Additionally, we shall require the $1+t_{n+1}$-dimensional vectors
${\bf e}_k$ defined as 
\eqn\newiunit{({\bf e}_k)_j=\cases{\delta_{j,k}& $1\leq k\leq 1+ t_{n+1}$\cr
0& $k=0.$}}
Solving \syn~ for $\vm$ yields
\eqn\synm{\vm={\bf M}^{-1}(\vn-{L\over 2}{\bf{\bar e}}_1-{1\over 2}{\bu}).}

The asymmetric matrix $-2\bf M$ is almost of block diagonal form with each
block except the last being the tadpole Cartan matrix $A_{\nu_i},$ and the
last block being a regular $A_{\nu_n}$ Cartan matrix (unless $\nu_n=0$). 
Note that in the simplest case $p=1$ the matrix $-2\bf M$ becomes the Cartan 
matrix $A_{p'-3}.$ For these reasons it is natural to think of $-2\bf M$ as a
new mathematical construct: a Cartan matrix of fractional
size. Algebraic structures associated with $-2\bf M$ will be
investigated elsewhere.

{}From the system of equations~\mnsys~one may deduce the following partition
problem for $L$ by multiplying equation $i$ in~\mnsys~by $l_i$ and summing
all equations up
\eqn\partprob{\sum_{i=1}^{t_{n+1}} (n_i-{\bar{u}_i\over 2}) 
l_i+{m_{(t_{n+1})}\over
2}l_{(1+t_{n+1})}={L\over 2}}
where $l_{j}$ is given by ~\ldef~and $l_{(1+t_{n+1})}=p'-2y_n.$
This partition problem can be employed to carry out an analysis along the 
lines of \rw. We will also need it later in the proof of the initial 
conditions and to determine allowed variable changes in the proof of
the recurrences.

Finally we note that it is useful to consider the cases $p'>2p$ 
and $p<p'<2p$ separately because
the case $p'>2p$ may be obtained from the case $p'<2p$ by a
``duality map''~\rbm. We concentrate first on the case
$p'>2 p$ where $\nu_0 \geq 1.$ In this case the $\nu_0$ types of particles in
the $0^{th}$ zone are treated differently than the remaining
particles. At times we will find it necessary to
 use $n_k,~1\leq k \leq \nu_0$ as the independent variables in the
$0^{th}$ zone and $m_k$ as the independent variables in all other
zones.  Thus we define the vector of independent variables as
\eqn\mtil{{\tilde {\bf m}}=\{n_1,n_2,\cdots,n_{\nu_0},m_{\nu_0+1},
m_{\nu_0+2},\cdots,m_{t_{n+1}}\}} 
and the vector of dependent variables as
\eqn\ntil{{\tilde {\bf n}}=\{m_1,m_2,\cdots,m_{\nu_0},n_{\nu_0+1},
\cdots, n_{t_{n+1}}\}.}
{}From~\mnsys~we find that
\eqn\dep{{\tilde{\bf n}}=-{\bf B}{{\tilde{\bf m}}}+L{\bf
{\bar E}}_{1,\nu_0}+{L\over 2}{\bf {\bar e}}_{\nu_0+1}
+{{\bf B}\over 2}{\bu}_{+}
+{1\over 2}{\bf {\bar e}}_{\nu_0+1}({\bu}^T_{+}\cdot{\bf V})+{1\over
2}{\bu}_{-}}
where the nonzero elements of the matrix  ${\bf B}$ are
given in terms of the $\nu_0\times \nu_0$ matrix
\eqn\ctm{C_T^{-1}=
\pmatrix{1&1&1&\ldots&1\cr
         1&2&2&\ldots&2\cr
	 1&2&3&\ldots&3\cr
	 \vdots&\vdots&\vdots&\ddots&\vdots\cr
	 1&2&3&\ldots&\nu_0\cr}}
(where $-{1\over 2}C_T$ is the matrix of the first $\nu_0$ rows and columns of
the matrix $M$ of \syn)~as
\eqn\bmat{\eqalign{{}&B_{i,j}=2(C_T^{-1})_{i,j}~~{\rm for}~~1\leq
  i,j\leq \nu_0\cr
  {}&B_{\nu_0+1,j}=B_{j,\nu_0+1}=j~~{\rm for}~~1\leq j \leq
  \nu_0\cr
  {}&B_{j,j}={{\nu_0}\over 2}\delta_{j,\nu_0+1}
  +(1-{1\over 2}\sum_{i=2}^{t_n}\delta_{j,t_i})~~{\rm for}~~\nu_0+1\leq j 
  \leq t_{n+1}-1\cr
  {}&B_{t_{n+1},t_{n+1}}=1\cr
  {}&B_{j,j+1}=-{1\over 2}+\sum_{i=2}^n \delta_{j,t_i}~~{\rm
  for}~j>\nu_0\cr
  {}&B_{j+1,j}=-{1\over 2}~{\rm for}~j>\nu_0.}}
The vector ${\bu}$ is decomposed as
\eqn\udecop{{\bu}={\bu}_+ +{\bu}_-}
with
\eqn\uplus{{\bu}_{+}=\sum_{i=1}^{\nu_0}\bar{u}_i {\bf {\bar e}}_i}
\eqn\uminus{{\bu}_-=\sum_{i={\nu_0+1}}^{t_{n+1}}\bar{u}_i{\bf {\bar e}}_i}
\eqn\rhovec{{\bf V}=\sum_{i=1}^{\nu_0}i{\bf {\bar e }}_{i}}
and
\eqn\eevecdef{{\bf {\bar E}}_{a,b}=\sum_{i=a}^{b}{\bf {\bar e}}_i.}
The splitting of~ \udecop~can be done for any vector.
We will, for example, also use this splitting for the vector 
$\tilde{\vm}$ in the following.

\newsec{Summary of results}
 
Now that we have summarized the needed results of ~\rts~we may
complete the specification of the bosonic and fermionic sums which
appear in our identities. We first complete the bosonic specification
in subsection A and give the definitions needed for the fermionic
sums in subsection B. The final results will be outlined in subsection C,
but the detailed identities will be presented in sec. 10.
We conclude in subsection D with a discussion of the
special case $q=1.$

{\bf A. Bosonic polynomials}

In order to complete the specification of the bosonic polynomials
$B_{r,s}(L,b)$ of ~\polyfb~ we represent both $b$ and $s$ as series
in Takahashi lengths. Thus we write
\eqn\becopm{\eqalign{b&=\sum_{i=1}^{\beta}l_{1+j_{\m_i}}^{(\m_i)} \cr
&{\rm
with}~~0 \leq \m_1< \m_2 < \cdots < \m_{\beta}\leq n
~~{\rm and}~~1+t_{\m_i}\leq j_{\m_i}
\leq t_{1+\m_i}+\delta_{n,\m_i}}}
with the further restriction that 
\eqn\restriction{{\rm if}~j_{\m_i}=t_{1+\m_i}~~{\rm then}~\m_{i+1}\geq \m_i+2.}
This decomposition is unique. We often say that a $b$ of this form lies 
in zone $\m_{\beta}$. Hence if $b$ is in zone $\mu$ then
\eqn\bzone{l^{(\m)}_{2+t_\m}\leq b\leq l^{(\m+1)}_{2+t_{\m+1}}-1.}

A similar decomposition is made for $s.$

{}From the decomposition of $b$ we specify $r(b)$ as follows
\eqn\rbgen{r(b)=\cases{\sum_{i=1}^{\beta} {\tilde
l}_{1+j_{\m_i}}^{(\m_i)}+\delta_{\m_1,0}&
if $1\leq b \leq p'-\nu_0-1,$\cr
\sum_{i=1}^{\beta }{\tilde l}^{(\m_i)}_{1+j_{\m_i}}& if
$p'-\nu_0 \leq b \leq p'-1.$\cr}}
where ${\tilde l}_{1+j_{\m}}^{(\m)}$  as defined by~\ltildem.~
This generalizes the case considered in~\rbm~of
 $b=l^{(\m)}_{1+j_{\m}}$ being a single
Takahashi length with $\m\geq 1$ where we had
\eqn\rpure{r(l_{1+j_{\m}}^{(\m)})={\tilde l}_{1+j_{\m}}^{(\m)}.}

One may prove that
\eqn\alter{\lfloor b {p\over p'}\rfloor=\sum_{i=1}^{\beta}{\tilde
l}^{(\m_i)}_{1+j_{\m_i}}-a_{\m_1}}
where $\lfloor x \rfloor$ is the greatest integer contained in $x$ and
\eqn\amdefn{a_{\m_1}=\cases{1& if $\m_1\geq 2$ and $\m_1$ even,\cr
0& otherwise.}}
Using~\alter~the map~\rbgen~ may be alternatively expressed as
\eqn\alterrnew{r(b)=\cases{\lfloor b{p\over p'}\rfloor+\kie{\m_1}&for $1\leq
b\leq p'-\nu_0-1$\cr
\lfloor b{p\over p'}\rfloor&for $p'-\nu_0\leq b \leq p'-1$\cr}}
where
\eqn\fdefine{\theta(A)=\cases{1& if $A$ is true,\cr 0& if $A$ is false.}}
When $r$ is expressed using relation~\alterrnew~we refer to the expression as
the Forrester-Baxter form for $r$ after the closely related formula of \rfb.
The proof of ~\alter~is straightforward and will be omitted.
The following elementary properties of the $b\rightarrow
r$ map are clear from ~\rbgen~and ~\ltildem: 
\eqn\rpropi{{\rm 1.}~~~r(b+1)=r(b)~~{\rm  or}~~ r(b)+1,} 
\eqn\rpropii{{\rm 2.~~~if}~ r(b+1)=r(b)+1~{\rm then}~ r(b+2)=r(b+1)} 
3.~~~if $b$ is as in ~\becopm,~\restriction~with $\mu_1\geq 1$ then 
\eqn\rpropiii{r(b-1)=r(b)=r(b+1)-1=r(b+2)-1}
and
\eqn\rpropiv{r(b+x)-r(x)=\cases{\sum_{i=1}^{\beta}{\tilde
l}^{(\mu_i)}_{1+j_{\mu_i}}&for~$b<p'-y_n=l^{(n)}_{2+t_{1+n}},~1\leq
x\leq y_{\mu_1}$\cr
{\tilde l}^{(n)}_{2+t_{1+n}}-\theta(y_n-y_1 < x)&for
$b=l^{(n)}_{2+t_{1+n}},~1\leq x\leq y_n-1$\cr}}
\eqn\rpropv{{\rm 4.}~~~r(1)=r(2)=1,~~~r(p'-1)=p-1}

Whereas the bosonic polynomials $B_{r,s}(L,b;q)$ depend on the Takahashi
decomposition of $b$ only through the map $r(b)$, the corresponding fermionic
polynomials depend sensitively on the details of the Takahashi decomposition
of $b$ and $s.$ In ~\rbm~ we treated the simplest case where both
$b$ and $s$ consist of one single Takahashi length. In this paper we
will still consider $s$ to be of the form
\eqn\sdecom{s=l^{(\m_s)}_{1+j_{s}}~({\rm or}~p'-l^{(\m_s)}_{1+j_{s}})
~~~1+t_{\m_s} \leq j_s \leq t_{\m_s+1}}
but $b$ is often left arbitrary ($\m_s$ and $j_s$ are from now on reserved
to specify $s$ as in \sdecom~and are not to be confused with $\m_i$ or
$j_{\mu_i}$ of the Takahashi decomposition \becopm).
However, the complexity of the final 
fermionic polynomials will depend on the details of the
decomposition of $b.$

{\bf B. Fermionic polynomials}
 
Our fermionic sums $F_{r,s}(L,b;q)$ will be constructed out of
two elementary fermionic objects $f_s(L,{\bf u};q)$ and
${\tilde f}_s(L,{\bf {\tilde u}};q)$ where the vectors ${\bf u}$ 
and ${\bf{ \tilde u}}$ depend on $b.$
The objects $f_s(L,{\bf u};q)$ and ${\tilde f}_s(L,{\bf {\tilde u}};q)$
are polynomials (in $q^{1/4}$)  and
differ from the polynomials of ~\rbm~in that 
the value at $q=0$ is not normalized to 1. This choice of
normalization is made for later convenience and we trust that no
confusion will result from referring to $f_s(L,{\bf u};q)$ and
${\tilde f}_s(L,{\bf {\tilde u}};q)$ as polynomials in the sequel. 

In order to define $f_s(L,\vu;q)$ and
$\tilde{f}_s(L,{\bf {\tilde u}};q)$ 
we need to discuss two generalizations of the 
$q$-binomials in~\qbindef~which allow $n$ to be negative for $n+m<0$, but will
automatically vanish if $m<0$.

{\bf Definition of $q$-binomials}

We use here the definitions of ~\rgr~
\eqn\newqbind{{m+n\atopwithdelims[] n}_q^{(0)}=
{m+n\atopwithdelims[] m}_q^{(1)}=
\cases{ {(q^{n+1})_m\over (q)_m} & for $m\geq 0$, $n$ integer,\cr
0 & otherwise,}}
where 
\eqn\qsym{(a;q)_n=(a)_n=\prod_{j=0}^{n-1}(1-aq^j),~1\leq n; (a)_0=1.}
The reason for distinguishing between the two $q$-binomials with superscript
zero or one is that the first one is more convenient if the variables
$n_j$ are taken as independent in the $(\vm,\vn)$-system whereas the latter
one is more convenient if the $m_j$ are taken as independent. 
We remark that when $m$ and $n$ are  nonnegative we have the symmetry
\eqn\qbinprop{{m+n\atopwithdelims[] n}_q^{(0,1)}=
{m+n\atopwithdelims[] m}_q^{(0,1)}.}
We also note that we have the special values for $L$ positive
\eqn\defnvala{\eqalign{{-L\atopwithdelims[] -L}_q^{(0)}
&=1,~~~{-L \atopwithdelims[] 0}_q^{(0)}=0,\cr
{-L \atopwithdelims[] -L}_q^{(1)}
&=0,~~~{-L \atopwithdelims[] 0}_q^{(1)}=1.\cr}}
Both functions satisfy the two recursion relations
without restrictions on $L$ or $m$
\eqn\qbinrra{{L \atopwithdelims[] m}_q^{(0,1)}
={L-1 \atopwithdelims[] m}_q^{(0,1)}+q^{L-m}{L-1
\atopwithdelims[] m-1}_q^{(0,1)},}
\eqn\qbinrrb{{L \atopwithdelims[] m}_q^{(0,1)}=
q^m{L-1 \atopwithdelims[] m}_q^{(0,1)}+{L-1
\atopwithdelims[] m-1}_q^{(0,1)}.}

This extension of the $q$-binomial coefficients is needed for 
the fermionic function $F_{r(b),s}(L,b;q)$ for general
values of $b$ and $s$. In ~\rbm~we used the definition~\qbindef~and thus
obtained more limited results. 

{\bf {Definition of} $f_s(L,{\bf u};q)$}

Let us now define $f_s(L,\vu;q)$ with ${\bf u}\in Z^{1+t_{1+n}}$
where $L$ is a nonnegative integer. 
There are actually several equivalent
forms, differing only by which set of $m_k$ and $n_k$ in the 
$(\vm,\vn)$--system~\mnsys~is taken to be independent. In the following
three equivalent forms for $f_s(L,\vu;q)$ are given with
$\tilde{\vm},\vn$ and $\vm$ as independent variables, respectively.
The corresponding variables $\tilde{\vn},\vm$ and $\vn$ can be determined
by equations \dep, \synm~and \syn, respectively.

\eqn\flq{\eqalign{f_s&(L,{\bf u};q)=\sum_{\tilde{\vm}_- \in 2\Z^{t_{n+1}-\nu_0}
+\vw_-(u_{1+t_{n+1}},{\bu}) \atop \tilde{m}_+ \in \Z^{\nu_0}}
q^{\tilde{Q}({\tilde{\bf m}})+{\bf A}^T{\tilde{\vm}}}\cr
\times&\prod_{j=1}^{t_{n+1}}{(({\bf I}_{t_{n+1}}-{\bf B}){\tilde {\bf
m}}+L{\bf{\bar E}}_{1,\nu_0}+{L\over 2}{\bf{\bar e}}_{\nu_0+1}+
{1\over 2}{\bf B}{\bf\bar{u}}_{+}
+{1\over 2}{\bf{\bar e}}_{\nu_0+1}({\bf{\bar u}}^T_{+}\cdot{\bf V})+{1\over
2}{\bf {\bar u}}_{-})_j \atopwithdelims[] 
{\tilde m}_j}^{(\theta(j>\nu_0))}_q\cr
&=\sum_{\vn \in \Z^{t_{n+1}}}\theta (m_{t_{n+1}}\equiv u_{1+t_{n+1}}({\rm
mod}~ 2))
 q^{Q(\vn,\vm)+\vA^T\tilde{\vm}}
\prod_{j=1}^{t_{n+1}}{m_j+n_j\atopwithdelims[] n_j}^{(0)}_q\cr
&=\sum_{\vm \in 2 \Z^{t_{n+1}}+\vw(u_{1+t_{n+1}},{\bu})}
q^{Q(\vn,\vm)+\vA^T\tilde{\vm}}
\prod_{j=1}^{t_{n+1}}{(({\bf I}_{t_{n+1}}
+{\bf M})\vm+{{\bu}\over 2}+{L\over
2}{\bf{\bar e}}_1)_j\atopwithdelims[] m_j}_q^{(1)}.}}
Here $\vB,{\bf M},{\bf I}_{t_{n+1}},
{\bf {\bar E}}_{a,b}$ and ${\bf V}$ have been defined
in secs. 1,2 and
$\Z$ is the set of integers $\{\ldots,-2,-1,0,1,2,\ldots\}.$
The vectors $\vm,\vn$ satisfy the system~\mnsys~with ${\bu}$
\eqn\ubardf{{\bf \bar{u}}={\bf \bar{u}}(s)+
\sum_{i=1}^{t_{n+1}}{\bf{\bar e}}_i u_i}
where for $s$ as in \sdecom~and $1\leq k \leq t_{1+n}$
\eqn\sstring{(\bu(s))_k=\cases{\delta_{k,j_{s}}-\sum_{i={\m_s}+1}^n
\delta_{k,t_i}& for $1+t_{\m_s}\leq j_{s} \leq t_{\m_s+1}$ and $\m_s\leq
n-1$,\cr
\delta_{k,j_{s}}& for $t_n<j_s$, $\m_s=n$.\cr}}
Note that the vectors ${\bf \bar u}$ and $ {\bf \bar{u}}(s)$  in \ubardf~have 
dimension $t_{n+1}$ while $\bf u$ has dimension $1+t_{n+1}.$ 

Let us now define the quadratic forms $\tilde{Q}(\tilde{\vm})$ and 
$Q(\vn,\vm)$, the linear term $\vA^T \tilde{\vm}$ and the parity restriction
vector ${\bf w}(u_{1+t_{n+1}},{\bu})$:

{\bf Definition of the quadratic forms} $Q$ {\bf and} $\tilde{Q}$

We define the quadratic form $\tilde{Q}(\tilde{\vm})$ as
\eqn\qform{\tilde{Q}({{\tilde {\bf m}}})={1\over 2}\tilde{\vm}^T\vB\tilde{\vm}}
where $\vB$ is defined by ~\bmat.
An equivalent form makes use of ~\mtil--\dep~ to write
\eqn\phib{Q(\vn,\vm)=-{1\over 2}
\tilde{\vm}^T({\bf {\tilde n}}-L{\bf{\bar E}}_{1,\nu_0}-{L\over 2}{\bf
{\bar e}}_{\nu_0+1}-{{\bf B}\over 2}{\bu}_+-
{1\over 2}{\bf {\bar  e}_{\nu_0+1}}({\bu}_+^T\cdot{\bf V})
-{{\bf{\bar u}_-}\over 2}).}

{\bf The linear term} $\vA^T\tilde{\vm}$

We write
\eqn\lincomplete{{\bf A}= {\bf A}^{(b)}+ {\bf A}^{(s)}}
where the $t_{n+1}$-dimensional vector ${\bf A}^{(b)}$ 
is obtained from $\bf u$ as
\eqn\lintermr{A_k^{(b)}=\cases{-{1\over 2} u_k& for $k$ in an
even, nonzero zone,\cr
0& otherwise,\cr}}
and the $t_{n+1}$-dimensional vector ${\bf A}^{(s)}$ 
is obtained from ${\bf \bar{u}}(s)$ as
\eqn\linterms{A_{k}^{(s)}=\cases{-{1\over 2}\bar{u}(s)_k
-{1\over 2}{\bf V}^T\cdot  
 \bu(s)\delta_{k,\nu_0+1}&for $k$ in an odd zone,\cr
-{1\over 2}{\bf {\bar e}}_k^T {\bf B} \bu(s)_+
& for $k$ in an even zone.\cr}}

{\bf The parity restrictions and the vector ${\bf w}(u_{1+t_{n+1}},{\bu})$}

The summation variables $\tilde{m}_j$ $(j=\nu_0+1,\ldots,t_{1+n})$ in the first
line of~\flq~and $m_j$ $(j=1,\ldots,t_{1+n})$ in the third sum of 
~\flq~are subject to even/odd
restrictions which are determined from $\bf u$ and $s$ by the requirement that
the entries of all $q$-binomials in ~\flq~ are integers as long as ${\bf
u}\in {\bf Z}^{1+t_{n+1}}$.
To formulate this analytically we define 
$w_k^{(j)}$ for 
\eqn\wren{1\leq k\leq t_{n+1},~~t_{\m_0}+\delta_{\m_0,0}+1\leq
j\leq t_{\m_0+1}+\delta_{\m_0,n},~~{\rm for~some}~~0\leq \m_0\leq n}
as
\eqn\wvecdf{w_k^{(j)}=\cases{0&$k \geq j$\cr
j-k&$t_{\m_0}\leq k \leq j-1$\cr
w_{k+1}^{(j)}+w^{(j)}_{1+t_{\m+1}}&$t_{\m}\leq k<t_{\m+1},
~{\rm for~some}~\m<\m_0$\cr}}
and then define
\eqn\wpub{{\bf w}(u_{1+t_{n+1}},{\bu})=\sum_{k=1}^{t_{n+1}}{\bf
{\bar e}}_k\left(u_{1+t_{1+n}}w_k^{(1+t_{n+1})}+
\sum_{j=1}^{t_{n+1}}w^{(j)}_k{\bar u}_j\right).}
Then the variables $m_j$ in ~\flq~satisfy
\eqn\eores{{\bf m}\equiv{\bf w}(u_{1+t_{1+n}},
{\bf \bar u})~({\rm mod}~2\Z^{t_{1+n}}).}

Note that $m_{t_{n+1}}\equiv u_{1+t_{n+1}}\equiv P({\rm mod}~2)$ 
where parity $P
\in \{0,1\}.$ Clearly, the $1+t_{1+n}$ component of the vector $\bf u$
determines the parity of the $m_{t_{1+n}}$ variable in the fundamental
fermionic polynomials~\flq.

{\bf The limit} $L\rightarrow \infty$

In order to obtain character identities from the polynomial identities
we need to consider the limit $L\rightarrow \infty.$ 
Only the first expression for $f_s(L,b;q)$ in~\flq~is suitable for this
limit since $Q(\vn,\vm)$ which appears in the other two expressions depends
on $L$. Hence in the limit $L\to\infty$ the
expression for the polynomial  ~\flq~ reduces to
\eqn\limfpoly{\eqalign{&f_s({\bf u};q)=\lim_{L\rightarrow
\infty}f_s(L,{\bf u};q)
=\sum_{\tilde{\vm}_+ \in \Z_{\geq 0}^{\nu_0} \atop \tilde{\vm}_- \in
2\Z^{t_{n+1}-\nu_0}+\vw(u_{1+t_{n+1}},{\bu})}
q^{\tilde{Q}({\tilde{\bf  m}})+{\bf A}^T{\tilde {\bf
m}}}\cr
&~~~~~~~~~~~~~\times
\left( \prod_{j=1}^{\nu_0+1}{1\over (q)_{\tilde{m}_j}} \right)
\prod_{j=2+\nu_0}^{t_{n+1}}{(({\bf I}_{t_{n+1}}-{\bf B}){\bf {\tilde
m}}+{1\over 2}{\bu}_{-})_j\atopwithdelims[] \tilde{m}_j}_q^{(1)}\cr}}
which depends on the vector ${\bu}_{+}$ only through the linear term
${\bf A}.$

{\bf Definition of} ${\tilde f}_s(L,{\bf \tilde u};q)$

The fermionic polynomial ${\tilde f}_s(L,{\bf \tilde u};q)$ is defined for
vectors of the form ${\bf \tilde u}=\e{\nu_0-j_0-1}-\e{\nu_0}+\vu'_1$ where
$0\leq j_0\leq \nu_0$ is in the zone $0$ and $\vu'_1\in Z^{1+t_{1+n}}$ is any
vector with no components in the zero zone, i.e. $({\bf u}'_1)_i=0$ for
$1\leq i\leq \nu_0$. We define
\eqn\ftildedf{ {\tilde f}_s(L,\e{\nu_0-j_0-1}-\e{\nu_0}+\vu'_1;q)=\cases{
q^{-{L+j_0\over 2}}
[f_s(L+1,{\bf e}_{\nu_0-j_0}-{\bf e}_{\nu_0}+\vu'_1;q)&~\cr
~~~~~~~~~~-f_s(L,{\bf e}_{1+\nu_0-j_0}-{\bf e}_{\nu_0}+\vu'_1;q)]&for
$1\leq j_0 \leq \nu_0,$\cr
(q^{L\over 2}-q^{-{L\over 2}})f_s(L-1,\vu'_1;q)&~\cr
~~~~~~~~~+q^{-{L\over2}}f_s(L,{\bf e}_{\nu_0-1}-{\bf e}_{\nu_0}+\vu'_1,q)&
for $j_0=0.$}}
When $j_0=\nu_0$, the left hand side of this definition loses
its meaning because we have no ${\bf e}_{-1}$. However, for conformity,
we introduce the notation 
$\tilde f_s(L,{\bf e}_{-1}-{\bf e}_{\nu_0}+{\bf u}'_1;q)$
to mean the first line of the right hand side of ~\ftildedf~ with $j_0=\nu_0$.

We will show in sec. 6 that for $j_0=\nu_0$ this definition reduces to
\eqn\ftildeff{{\tilde f}_s(L,{\bf e}_{-1}-{\bf e}_{\nu_0}+\vu'_1;q)=0,}
while for $j_0=\nu_0-1$ 
\eqn\fftilde{\tilde{f}_s(L,\ve_0-\e{\nu_0}+\vu'_1;q)=q^{{L-(\nu_0-1)\over 2}}
f_s(L-1,\e{1}-\e{\nu_0}+\vu'_1;q)=q^{L-(\nu_0-1)\over 2}
f_s(L,-\e{\nu_0}+\vu'_1;q).}
In the limit $L\rightarrow \infty$ we have
\eqn\limft{\lim_{L\rightarrow \infty}{\tilde f}_s(L,{\bf \tilde u};q)=0.}

{\bf C. The identities.}

The polynomial identities proven in this paper are all of the form
\eqn\polid{F_{r(b),s}(L,b;q)=q^{{\cal N}_{r(b),s}} B_{r(b),s}(L,b;q)}
where $r(b)$ as given in~\rbgen~and $F_{r(b),s}(L,b;q)$ is of the form
\eqn\polyform{F_{r(b),s}(L,b;q)=\sum_{\vu\in U(b)} q^{c_{\bf u}+
{L\over 2}g_{\bf u}}f_s(L,{\bf u},q)
+\sum_{\tilde{\vu}\in \tilde{U}(b)}
q^{{\tilde c}_{{\bf {\tilde u}}}}{\tilde f}_s(L,{\bf{\tilde u}},q).}
Here $U(b)$ and $\tilde{U}(b)$ are sets of vectors determined by $b$.
The sets $U(b)$, $\tilde{U}(b)$ and the exponents $c_{\vu},~
\tilde{c}_{\tilde{\vu}},~g_{\vu}\geq 0$ and ${\cal N}_{r(b),s}$ 
will be explicitly computed.
We note however that there are very special values for $b$ for 
which the representation of $F_{r(b),s}(L,b;q)$  in \polyform~is not correct.
This phenomenon will be discussed in sec. 8.

The results depend sensitively on the details of the expansion of $b$
~\becopm. For the cases of
1,2 and 3 zones the identities for all values of $b$ are given in
sec. 7. However for the general case of $n+1$ zones with $n\geq 3$  
there are many special cases to consider. The process of obtaining a
complete set of results is tedious and in this paper we present
explicit results only for the following special cases of $b$ (and $r$)
\eqn\results{\eqalign{1:&~~b=l^{(\m)}_{1+j_{\m}},~~
r={\tilde l}^{(\m)}_{1+j_{\m}}+\delta_{\m ,0}~~
\cr
2a:&~~~l^{(\m)}_{1+j_{\m}}-\nu_0+1\leq b \leq l^{(\m)}_{1+j_{\m}}-1,
~~r={\tilde l}^{(\m)}_{1+j_{\m}}~~\cr
&~~~~~1\leq \m,~~t_{\m}+1\leq j_{\m}\leq t_{1+\m}-1+\delta_{\m,n}\cr
2b:&~~~l_{1+j_{\m}}^{(\m)}+1\leq b\leq l^{(\m)}_{1+{j_{\m}}}+\nu_0+1,
~~r=1+{\tilde l}^{(\m)}_{1+j_{\m}}~~\cr
&~~~~~1\leq \m,~~t_{\m}+1\leq j_{\m}\leq t_{1+\m}-1+\delta_{\m,n}\cr
3:&~~b=\sum_{\m=0}^{\beta}l^{(\m)}_{1+j_{\m}},~~~
r=1+\sum_{\m=1}^{\beta}{\tilde l}^{(\m)}_{1+j_{\m}},~~~
1\leq \beta \leq n\cr
&~~~~0\leq j_0\leq \nu_0-1\cr
&~~~~1+t_{\m}\leq j_{\m}\leq t_{\m+1}-3~~1\leq \m \leq \beta-2\cr
&~~~~1+t_{\beta-1}\leq j_{\beta-1}\leq t_{\beta}-2\cr
&~~~~1+t_{\beta}\leq j_{\beta}\leq t_{\beta+1}-1+\delta_{\beta,n}\cr
4:&~~b=\sum_{\m=\alpha}^{\beta}l^{(\m)}_{1+j_{\m}},~~
r=\sum_{\m=\alpha}^{\beta}{\tilde l}^{(\m)}_{1+j_{\m}},
~~~1 \leq \alpha,~\alpha+1\leq \beta\leq n\cr
&~~~~1+t_{\m}\leq j_{\m}\leq t_{\m+1}-3~~\alpha\leq \m \leq \beta-2\cr
&~~~~1+t_{\beta-1}\leq j_{\beta-1}\leq t_{\beta}-2\cr
&~~~~1+t_{\beta}\leq j_{\beta}\leq t_{\beta+1}-1+\delta_{\beta,n}\cr
&~~~{\rm where}~\alpha=1~{\rm and}~\alpha \geq 2~
{\rm are~ treated~ separately.}}}

The results are given in sec. 10 as follows:

{\bf Polynomial Identities}
\eqn\anspoly{\eqalign{1:&~~~~~10.3\cr
2a:&~~~~~10.4\cr
2b:&~~~~~10.5\cr
3:&~~~~~10.10\cr
4:&~~~~~10.15.\cr}}

{\bf Character identities}
\eqn\anschar{\eqalign{1~{\rm and}~ 2a:&~~~~~10.20\cr
2b:&~~~~~10.21\cr
3:&~~~~~10.22\cr
4,~\alpha=1:&~~~~~10.23\cr
4,~\alpha \geq 2:&~~~~~10.24.\cr}}
These results hold for $\nu_i \neq 1,~i=1,\cdots,n-1$ and in addition
(10.3) and (10.20) are valid even if some or all $\nu_i=1$ or $\nu_n=0$
(or both) with minor modifications of overall factors. Further results are 
to be found in sec. 11, sec. 12 and (C.19) of appendix C. 
The identities of sec.10 (11) are valid for $p'>2p$ and for
$s=l^{(\mu_s)}_{1+j_s}~(p'-l^{(\mu_s)}_{1+j_s}).$ The identities of
sec. 12 are valid for $p'<2p.$

{\bf D. The special case} $q=1$
 
We close this section with a brief discussion of the special case
$q=1.$ The details of this case will be presented elsewhere, but
it is useful to sketch the results here in order to give a
characterization of the vectors $\bf u$ and $\bf {\tilde
u}$ that appear in ~\polyform~which is complementary to the
constructive procedure developed in sections 7 and 8.

When $q=1$ the bosonic function ~\polyfb~ simplifies because the
dependence on $r$ vanishes and the fermionic form~\polyform~simplifies
because
\eqn\fqone{{\tilde f}_s(L,{\bf{\tilde u}};q=1)=f_s(L,{\bf {\tilde  u}};q=1).}
For these reasons many of the distinctions between the special cases
noted in the previous section disappear and we find the single
identity valid for all $b$ 
\eqn\qoneid{B_{r(b),s}(L,b;q=1)=\sum_{\vu\in W(b)}f_s(L,\vu;q=1)}
where the set of vectors $W(b)$ is obtained from the Takahashi
decomposition of $b$ ~\becopm~as follows:

If in ~\becopm~$\beta=1$ so that 
\eqn\veca{b=l^{(\m)}_{1+j_{\m}}}
then we define as in ~\rbm~$W(b)$ having only the single element
\eqn\usinbsingle{\vu={\bf e}_{j_{\m}}-\sum_{i=\m+1}^{n}{\bf e}_{t_i},
~~1+t_{\m}\leq j_{\m}\leq t_{\m+1}+\delta_{\mu,n},~~\mu=0,1,\cdots ,n.}
This vector was called an $r$ string in ~\rbm. 

For the general case  we write the
Takahashi decomposition of  an arbitrary $b$ as
\eqn\anx{b=l_{1+j_{\m_{\beta}}}^{(\m_{\beta})}+x~~{\rm where}~~x=
\sum_{i=1}^{\beta-1}l_{1+j_{\m_i}}^{(\m_i)}}
where $\beta $ is defined from ~\becopm.
If $x=0$ we define a vector ${\bf v}^{(0)}$ as:
\eqn\vitera{{\bf v}^{(0)}={\bf e}_{j_{\m_{\beta}}}.}
If $x \neq 0$ and it is not true that $\m_\beta=n$ and $j_{n}
=1+t_{1+n}$ then we define two vectors as
\eqn\viterb{\eqalign{{\bf v}^{(0)}&=\cases{{\bf
e}_{1+j_{\m_{\beta}}}& for $t_{\m_{\beta}}+1\leq
j_{\m_{\beta}}\leq -1+ t_{1+\m_\beta},$\cr
{\bf e}_{1+t_{1+\m_{\beta}}}+(1-\delta_{\mu_{\beta},n}){\bf
e}_{t_{1+\m_{\beta}}}
& for $j_{\m_{\beta}}=
t_{1+{\m_\beta}},$\cr}\cr
{\bf v}^{(1)}&={\bf e}_{j_{\m_{\beta}}}.\cr}}
Furthermore if $x\neq 0$  and we do not
have $\m_{\beta}=n$ and $j_{n}=1+t_{1+n}$ we define the two numbers
\eqn\biter{b^{(0)}=x~~~{\rm and }~~~b^{(1)}=y_{\m_\beta}-x.}
In the exceptional case where $x\neq 0, ~\m_\beta=n$ and
$j_{n}=1+t_{1+n}$ we define 
only one vector and one number
\eqn\exceptb{\eqalign{{\bf v}^{(1)}&={\bf e}_{1+t_{1+n}},\cr
b^{(1)}&=y_{n}-x.\cr}} 

We then expand both $b^{(0)}$ and $b^{(1)}$ (and in the exceptional case only
$b^{(1)}$) again in a Takahashi series
and repeat the process as many times as needed until $x=0$ in which case
the process terminates. 
This recursion leads to a branched chain of vectors ${\bf
v}^{(i_1)},~{\bf v}^{(i_1,i_2)},\cdots ,{\bf v}^{(i_1,i_2,\cdots
,i_f)}$ where $1\leq f\leq \m_{\beta}+1;~i_1,\cdots , i_f=0,1$.
Notice however that $f$ might vary from branch to branch.
Let us define $\m_f$ to be the lowest $\m$ such that there exist an $i$
such that $t_{\m}<i\leq t_{1+\m}+\delta_{\m,n}$ and
$v_i^{(i_1,\ldots,i_f)}\neq 0$.
Then the set $W(b)$ consists of all vectors
\eqn\uvec{{\bf v}^{(i_1)}+{\bf v}^{(i_1,i_2)}+\cdots
+{\bf v}^{(i_1,i_2,\cdots ,i_f)}-\sum_{k=\m_f+1}^n{\bf e}_{t_k},}
where all vectors ${\bf v}$ are determined by the above described recursive
procedure.

The explicit solution of this recursion relation involves the
recognition that there are many separate cases of the Takahashi
decomposition of $b$ which lead to sets $W(b)$ which
may differ even in the number of vectors in the set. Certain of these
special cases correspond to the cases distinguished in the previous
section. The complete solution of this recursive definition will be
given elsewhere where we will use it to give explicit
Rogers-Schur-Ramanujan  identities for general values of $b$. We note 
in the cases considered in the previous subsection that when
$l^{(0)}_{1+j_0}$ is present in the Takahashi decomposition then there
are vectors which have components in zone zero of the form 
\eqn\mirror{\tilde{\vu}=\e{\nu_0-j_0-a}-\e{\nu_0}+{\bf u}'_1,~~~a=1,2}
where the vector ${\bf u}'_1$ is defined above ~\ftildedf. 
These are the vectors $\tilde{\vu} \in { \tilde U}(b)$ of
~\polyform. Note that the set ${ \tilde U}(b)$ is empty if in the
Takahashi decomposition of $b$ ~\becopm~we have $\mu_1\geq 1.$  
The remaining vectors are the vectors ${\bf u} 
\in {U}(b)$ of ~\polyform. 

\newsec{The bosonic recursion relations}

In this section we derive recursion relations for the bosonic
polynomials $B_{r,s}(L,b;q)$ defined by ~\polyfb~and the $b\rightarrow r$
map ~\rbgen. Here and in the remainder of the paper we will often suppress the
argument $q$ of all functions unless explicitly needed. Moreover
 we find it convenient to remove an $L$
independent power of $q$ by defining 
\eqn\mobose{{\tilde B}_{r(b),s}(L,b)
=q^{{1\over 2}({\phi_{r(b),s}}-{\phi_{r(s),s}})} B_{r(b),s}(L,b)}
where ${\phi}_{x,y}$ is defined as
\eqn\frdefine{\phi_{x+1,y}-\phi_{x,y}=
1-y+\sum_{i=1}^{\xi}l_{1+j_{\eta_i}}^{(\eta_i)}}
when $x$ has the decomposition in terms of truncated Takahashi lengths
\eqn\xdec{x=\sum_{i=1}^{\xi}\tilde{l}^{(\eta_i)}_{1+j_{\eta_i}},}
with $1\leq \eta_1<\eta_2<\cdots<\eta_{\xi}\leq n$ and 
$1+t_{\eta_i}\leq j_{\eta_i}\leq t_{1+\eta_i}+\delta_{n,\eta_i}$.
We note that because ${\phi}_{x,y}$ appears only
as a difference in~\mobose~that boundary conditions on ${\phi}_{x,y}$
are not needed. This change in normalization allows us to prove the
following recursion relations which contain no explicit dependence on $s$.
For $2\leq b\leq p'-\nu_0-1$
\eqn\brecrelt{{\tilde B}_{r(b),s}(L,b)=
\cases{{\tilde B}_{r(b),s}(L-1,b+1)+{\tilde B}_{r(b),s}(L-1,b-1)&~\cr
~~~~+(q^{L-1}-1){\tilde B}_{r(b),s}(L-2,b)&if $\m_1=0,~j_0\geq 1$\cr
q^{L-1\over 2}{\tilde B}_{r(b)+1,s}(L-1,b+1)+{\tilde B}_{r(b),s}
(L-1,b-1)&if $\m_1 \geq 1$\cr
{\tilde B}_{r(b),s}(L-1,b+1)+q^{L\over 2}{\tilde B}_{r(b)-1,s}
(L-1,b-1)&if $\m_1=0$ 
$j_0=0$\cr}}
where $\m_1$ and $j_0$ are obtained from the Takahashi decomposition
of $b$ \becopm.
For the remaining cases we have 
\eqn\brecreli{{\tilde B}_{1,s}(L,1)={\tilde B}_{1,s}(L-1,2)}
and
\eqn\brecreltd{{\tilde B}_{p-1,s}(L,b)
=\cases{{\tilde B}_{p-1,s}(L-1,b+1)+{\tilde B}_{p-1,s}(L-1,b-1)&~\cr
~~~~~+(q^{L-1}-1){\tilde B}_{p-1,s}(L-2,b)&if $p'-\nu_0 \leq b \leq p'-2$\cr
{\tilde B}_{p-1,s}(L-1,p'-2)&if $b=p'-1$.}
}

The recursion relations~\brecrelt--\brecreltd~have a unique 
solution once appropriate boundary conditions are given. One set of
boundary conditions which will specify
${\tilde B}_{r(b),s}(L,b)$ as the unique solution are the values obtained from
\eqn\binitial{B_{r(b),s}(0,b)=\delta_{s,b},~~1\leq b \leq p'-1.}
However, it will prove useful to recognize that this is not the only
way in which boundary conditions may be imposed on ~\brecrelt--\brecreltd~
to give ~\polyfb~as the unique solution. One alternative
is the condition which is readily obtained from the term $j=0$ in ~\polyfb
\eqn\newinitial{B_{r(b),s}(L,b)=\cases{0& for $L=0,1,2,\cdots,|s-b|-1$\cr
1&for $L=|s-b|$\cr}}

To prove ~\brecrelt-\brecreltd~for the $b\rightarrow r$ map of sec. 3 we
first consider the case where $b$ and $r$ are unrelated and recall the
elementary recursion relations for $q$--binomial coefficients~\qbinrra~
and \qbinrrb.

If we use~\qbinrra~ in the definition~\polyfb~ 
 for $B_{r,s}(L,b)$ we find
\eqn\genrbrra{B_{r,s}(L,b)=B_{r,s}(L-1,b-1)+q^{L+b-s\over
2}B_{r+1,s}(L-1,b+1)}
and if we use~\qbinrrb~we find
\eqn\genrbrrb{B_{r,s}(L,b)=q^{L-b+s\over 2}B_{r-1,s}(L-1,b-1)+
B_{r,s}(L-1,b+1).}
Furthermore if we write~\genrbrra~as
\eqn\brect{q^{-{L+b-s\over 2}}
[B_{r,s}(L,b)-B_{r,s}(L-1,b-1)]=B_{r+1,s}(L-1,b+1)}
and
\eqn\brecu{q^{-{L+b-s\over 2}}
[B_{r,s}(L-1,b+1)-B_{r,s}(L-2,b)]
=B_{r+1,s}(L-2,b+2)}
and subtract \brecu~from \brect~
we may use~\genrbrrb~on the right hand side to find
\eqn\brecv{\eqalign{q^{-{L+b-s\over 2}}
[B_{r,s}&(L,b)-B_{r,s}(L-1,b-1)-
B_{r,s}(L-1,b+1)+B_{r,s}(L-2,b)]\cr
&=q^{{L+s-b\over 2}-1} B_{r,s}(L-2,b)}}
and thus we have
\eqn\brecw{B_{r,s}(L,b)=B_{r,s}(L-1,b-1)+B_{r,s}(L-1,b+1)+
(q^{L-1}-1)B_{r,s}(L-2,b).}
If we now relate $r$ to $b$ using the map of sec. 3 and use the definition~\mobose~ then
equations ~\genrbrra,~\genrbrrb~ and ~\brecw~are the three equations of
~\brecrelt~and the first equation of ~\brecreltd.

We also note from ~\polyfb~that
\eqn\brecx{\eqalign{B_{0,s}(L,0)
&= \sum_{j=-\infty}^{\infty}\left( q^{j^2pp'-jsp}{L
\atopwithdelims[] {L+s\over 2}-jp'}_q-q^{j^2pp'+jsp}{L \atopwithdelims[]
{L-s\over 2}-jp'}_q \right)\cr
&=\sum_{j=-\infty}^{\infty}\left( q^{j^2pp'-jsp}{L
\atopwithdelims[] {L+s\over 2}-jp'}_q-q^{j^2pp'+jsp}{L \atopwithdelims[]
{L+s\over 2}+jp'}_q \right)=0\cr}}
where to get the second line we first use ~\qbinprop~and then  
let $j\rightarrow -j.$
Combining this with ~\genrbrrb~ at $r=1,~b=1$ we obtain \brecreli.
 Finally we note that in an analogous fashion we may prove 
\eqn\brecy{B_{p,s}(L,p')=0}
and thus the last equation of ~\brecreltd~follows.

\newsec{Recursive properties for fundamental fermionic polynomials 
$f_s(L,{\bf u})$}

Our goal is to construct fermionic objects from the fundamental
fermionic polynomials ~\flq~which will satisfy the bosonic recursion 
relations~~\brecrelt-\brecreltd~ obeyed by 
${\tilde B}_{r(b),s}(L,b).$ To do this we will use the following recursive
properties for the fundamental fermionic polynomials $f_s(L,{\bf u})$
where here and in the remainder of the paper 
we restrict our attention to $\nu_i\geq 2$ for $i=1,\cdots , n-1.$ 
The analysis for $\nu_i=1$ is analogous to that for $\nu_i\geq 2$, but
since the recursion relations given below should be slightly modified, this case will not be treated here for reasons of space.

Let us first introduce the following notation for $t_{\m}+1\leq j_{\m}
\leq t_{\m+1}+\delta_{\m,n}$
\eqn\uvec{\eqalign{
\vu_0(j_{\m})&=\e{j_{\m}}+\vu'-\e{t_{1+\m}}(1-\delta_{\m,n})\cr
\vu_{\pm 1}(j_{\m})&=\e{j_{\m}\pm 1}+\vu'-\e{t_{1+\m}}(1-\delta_{\m,n})
}}
where $\vu'$ is an arbitrary vector only restricted by
$\sum_{i=1}^{\nu_0}u'_i=0$ for $\m=0$ and for
$\m\geq 1$ the components must satisfy $u'_j=0$ 
for $j\leq j_{\m}+1-\delta_{j_{\mu},t_{1+n}}$. 
These conditions are used in the proof of the recursive properties.
Further define
\eqn\etdefn{{\bf E}^{(t)}_{a,b}=\sum_{i=a}^{b}{\bf e}_{t_i}.}

Then if $j_0$ is in the zone $\m=0$ where $0\leq j_0\leq \nu_0=t_1$ we find
(for $L\geq 1$)
\eqn\freczero{f_s(L,\vu_0(j_0))
=\cases
{f_s(L-1,\vu_1(0)),~~~~j_0=0&~\cr
f_s(L-1,\vu_1(j_0))
+f_s(L-1,\vu_{-1}(j_0))&~\cr
~~+(q^{L-1}-1)f_s(L-2,\vu_0(j_0)),
~~~1\leq j_0 \leq \nu_0-1&~\cr
f_s(L-1,\vu_{1}(\nu_0)+\e{\nu_0})
+f_s(L-1,\vu_{-1}(\nu_0))&~\cr
~~+(q^{L-1}-1)f_s(L-2,\vu_0(\nu_0)),~~~~~j_0=\nu_0=t_1.&~\cr}}

For $j_{\m}$ in the zone $1\leq \m \leq n$ where 
$1+t_{\m}\leq j_{\m} \leq  t_{\m+1}+\delta_{\m,n}$ 
we have four separate recursive properties (for $L\geq 1$):

1. For $j_{\m}=1+t_{\m}$
\eqn\frecrelb{\eqalign{
f_s(L&,\vu_0(1+t_{\m}))=q^{{L-1\over 2}-{\nu_0+1\over
4}+{3\over 4}\kio{\m}}f_s(L-1,\vu_1(1+t_{\m})-\E{1}{\m})\cr
&+q^{\theta (\m \geq 2)({L-1\over 2}-{\nu_0-\kie{\m}\over 4})}
f_s(L-1,\vu_{-1}(1+t_{\m})-\E{1}{\m})\cr
&+\theta(\m \geq 2)q^{{L-1\over 2}-{\nu_0-\kio{\m}\over 4}}
f_s(L-1,\e{-1+t_{\m}}+\vu_0(1+t_{\m})-\E{1}{\m})\cr
&+2\theta(\m\geq 3)
\sum_{i=2}^{\m-1}q^{{L-1\over 2}-{\nu_0-\kio{i}\over 4}}
f_s(L-1,-\E{1}{i}+\e{-1+t_i}+\vu_0(1+t_{\m}))\cr
&+[1+\theta(\m\geq 2)]f_s(L-1,\e{-1+\nu_0}-\e{\nu_0}+\vu_0(1+t_{\m}))\cr
&+(q^{L-1}-1)f_s(L-2,\vu_0(1+t_{\m}));\cr}}

2. For $2+t_{\m}\leq j_{\m} \leq -1+t_{\m+1}+\delta_{\m,n}$
\eqn\frecrelc{\eqalign{f_s(L&,\vu_0(j_{\m}))
=q^{{L-1\over 2}-{\nu_0+1\over 4}+{3\over
4}\kio{\m}}f_s(L-1,\vu_1(j_{\m})-\E{1}{\m})\cr
&+q^{{L-1\over 2}-{\nu_0-\kie{\m}\over 4}}f_s(L-1,\vu_{-1}(j_{\m})
-\E{1}{\m})\cr
&+2\theta(\m\geq 2)\sum_{i=2}^{\m}q^{{L-1\over
2}-{\nu_0-\kio{i}\over 4}}
f_s(L-1,\vu_0(j_{\m})-\E{1}{i}+\e{-1+t_i})\cr
&+2f_s(L-1,\e{\nu_0-1}-\e{\nu_0}+\vu_0(j_{\m}))+(q^{L-1}-1)
f_s(L-2,\vu_0(j_{\m}));
\cr}}

3. For $1\leq \m \leq n-1,$ and $j_{\m}=t_{\m+1}$
\eqn\frecreld{\eqalign{
f_s(L&,\vu_0(t_{1+\m}))=q^{{L-1\over
2}-{\nu_0-\kio{\m}\over 4}}f_s(L-1,-\E{1}{\m}+\e{t_{1+\m}}
+\vu_{1}(t_{\m+1}))\cr
&+q^{{L-1\over 2}-{\nu_0-\kie{\m}\over 4}}
f_s(L-1,-\E{1}{\m}+\vu_{-1}(t_{\m+1}))\cr
&+2\theta (\m\geq 2)
\sum_{i=2}^{\m}q^{{L-1\over 2}-{\nu_0-\kio{i}\over
4}}f_s(L-1,-\E{1}{i}+\e{-1+t_i}+\vu_0(t_{1+\m}))\cr
&+2f_s(L-1,\e{-1+\nu_0}-\e{\nu_0}+\vu_0(t_{1+\m}))
+(q^{L-1}-1)f_s(L-2,\vu_0(t_{1+\m}));\cr}}  

4. For $\m=n$ and $j_n=1+t_{n+1}$
\eqn\frecrele{\eqalign{
f_s(L&,{\bf e}_{1+t_{1+n}})=q^{{L-1\over
2}-{\nu_0-\kie{n}\over 4}}f_s(L-1,\e{t_{1+n}}-\E{1}{n})\cr
&+2\theta(n\geq 2)\sum_{i=2}^{n}q^{{L-1\over
2}-{\nu_0-\kio{i}\over 4}}
f_s(L-1,-\E{1}{i}+\e{-1+t_i}+\e{1+t_{1+n}})\cr
&+2f_s(L-1,\e{\nu_0-1}-\e{\nu_0}+\e{1+t_{1+n}})
+(q^{L-1}-1)f_s(L-2,\e{1+t_{1+n}})\cr}}
where we remind the reader that the  $1+t_{1+n}$ component of the vector ${\bf
u}$ determines
the parity of the variable $m_{t_{1+n}}$ in the fundamental
fermionic polynomial~\flq.

We prove these recursive properties in Appendix A.

\newsec{Properties of ${\tilde f}_s(L,{\bf \tilde {u}})$}

In this section we will  prove the 

{\bf Recursive properties of ${\tilde f}_s(L,{\bf \tilde{u}})$}

{\it The function ${\tilde f}_s(L,\vu_0(\nu_0-j_0-1))$ 
defined by ~\ftildedf~and~\uvec~with ${\bf u}'\rightarrow {\bf u}'_1$
satisfies the recursive properties for $1\leq j_0\leq \nu_0-1$
\eqn\reccft{\eqalign{{\tilde f}_s(L,\vu_0(\nu_0-j_0-1))&=
{\tilde f}_s(L-1,\vu_{-1}(\nu_0-j_0-1))+
{\tilde f}_s(L-1,\vu_1(\nu_0-j_0-1))\cr
&+(q^{L-1}-1){\tilde f}_s(L-2,\vu_0(\nu_0-j_0-1)).\cr}}}

To prove ~\reccft~ we adopt the simplified notation $f_{j_0}(L)=
f_s(L,\vu_0(j_0)).$ Then for
$1\leq j_0\leq \nu_0-1$ the relation ~\freczero~ reads
\eqn\againfrec{f_{j_0}(L)=f_{1+j_0}(L-1)+f_{-1+j_0}(L-1)+
(q^{L-1}-1)f_{j_0}(L-2)}
which we rewrite as
\eqn\rewrite{f_{j_0}(L)-f_{1+j_0}(L-1)=
f_{-1+j_0}(L-1)-f_{j_0}(L-2)+q^{L-1}f_{j_0}(L-2).}
Then if we replace $j_0$ by $\nu_0-j_0,$ $L$ by $L+1$ and set
\eqn\idf{I(L,j_0)=f_{\nu_0-j_0}(L+1)-f_{1+\nu_0-j_0}(L)}
we have 
\eqn\irec{I(L,j_0)=I(L-1,j_0+1)+q^Lf_{\nu_0-j_0}(L-1).}
We now eliminate $f_{\nu_0-j_0}(L-1)$ by first writing ~\irec~as
\eqn\ireca{[I(L,j_0)-I(L-1,1+j_0)]q^{-L}=f_{\nu_0-j_0}(L-1)}
and then by setting $L\rightarrow L-1$ and $j_0\rightarrow j_0-1$ to
get the companion equation
\eqn\irecb{[I(L-1,j_0-1)-I(L-2,j_0)]q^{-L+1}=f_{1+\nu_0-j_0}(L-2).}
Subtracting ~\irecb~from ~\ireca~ and multiplying by
$q^{L-j_0\over 2}$we obtain
\eqn\irecc{\eqalign{&[I(L,j_0)-I(L-1,1+j_0)]q^{-{L+j_0\over 2}}\cr
&-[I(L-1,-1+j_0)-I(L-2,j_0)]q^{-{L-2+j_0\over 2}}=q^{L-j_0\over
2}I(L-2,j_0).}}
Recalling the definition~\ftildedf~ we see  that 
\eqn\moredefin{{\tilde f}_s(L,\vu_0(\nu_0-j_0-1))=q^{-{L+j_0\over 2}}I(L,j_0)}
and using this along with ~\irecc~we obtain~\reccft.
Note from~\freczero~that $f_{0}(L)=f_1(L-1)$ which implies  that
\eqn\irecd{{\tilde f}_s(L,\e{-1}-\e{\nu_0}+\vu'_1)=0=I(L,\nu_0)}
and using ~\ireca~we have
\eqn\irece{{\tilde f}_s(L,\e{0}-\e{\nu_0}+\vu'_1)=
q^{L-(\nu_0-1)\over 2}f_s(L-1,\e{1}-\e{\nu_0}+\vu'_1)=
q^{L-(\nu_0-1)\over 2}f_s(L,-\e{\nu_0}+\vu'_1).}

We also need the 

{\bf Limiting Relation}
\eqn\limitft{\lim_{L\rightarrow \infty}{\tilde f}_s(L,{\bf\tilde{u}})=0.}

To prove ~\limitft~we use ~\ireca~and ~\irecd~to obtain the system of
$\nu_0-j_0$ equations
\eqn\isysa{\eqalign{I(L,j_0)-I(L-1,j_0+1)&=q^Lf_{\nu_0-j_0}(L-1)\cr
I(L-1,j_0+1)-I(L-2,j_0+2)&=q^{L-1}f_{\nu_0-j_0-1}(L-2)\cr
I(L-2,j_0+2)-I(L-3,j_0+3)&=q^{L-2}f_{\nu_0-j_0-2}(L-3)\cr
&\cdots\cr
I(L-(\nu_0-j_0-1),\nu_0-1)-0&=q^{L-(\nu_0-j_0-1)}f_{1}(L-(\nu_0-j_0))\cr}}
which, if added together give
\eqn\ilima{I(L,j_0)=\sum_{i=0}^{\nu_0-j_0-1}q^{L-i}f_{\nu_0-j_0-1}(L-1-i).}
Thus multiplying by $q^{-{L+j_0\over 2}}$ we find 
\eqn\ilimb{{\tilde f}_s(L,\vu_0(\nu_0-j_0-1))
=q^{L\over 2}\sum_{i=0}^{\nu_0-j_0-1}q^{-{j_0\over
2}-i}f_{\nu_0-j_0-i}(L-1-i)}
and hence since $\lim_{L\rightarrow \infty}f_{j_0}(L)$ is finite we
see from~\ilimb~that the limiting relation~\limitft~holds for ${\bf
u}={\bf u}_0(\nu_0-j_0-1),~1\leq j_0 \leq \nu_0-1$ 

The case $j_0=0$ requires a separate treatment. First we note that
\eqn\neqlim{{\tilde f}_s(L,{\bf u}_0(\nu_0-1))=
{\tilde f}_s(L,{\bf u}_0(\nu_0-2))+q^{L\over 2}f_s(L-1,{\bf u}_{0}(\nu_0))}
which follows from ~\ftildedf. In this equation we let $L\rightarrow
\infty$ and using ~\ilimb~with $j_0=1$ we obtain ~\limitft~for ${\bf
\tilde u}={\bf u}_0(\nu_0-1).$

\newsec{Construction of $F_{r(b),s}(L,b)$}

We now turn to the details of the construction of the fermionic
representations $F_{r(b),s}(L,b),~L\geq 0$ of the bosonic 
polynomials ${\tilde B}_{r(b),s}(L,b).$ Our
method is to construct  fermionic functions $F_{r(b),s}(L,b)$ in terms
of $f_{s}(L,{\bf u})$ and $\tilde{f}_s(L,{\tilde \vu})$ 
which by construction satisfy the bosonic
recursion relations ~\brecrelt~and ~\brecreltd.

We begin by choosing  as a starting value for $b=1$
\eqn\prfb{F_{1,s}(L,1)=f_s(L,-\E{1}{n}).}
The boundary bosonic recursion relation~\brecreli~ 
requires that 
\eqn\prfc{F_{1,s}(L-1,2)=f_s(L,-\E{1}{n})}
from which if we let $L\rightarrow L+1$ and use the first 
fermionic recursive properties in~\freczero~we find with $\vu'=-\E{2}{n}$
\eqn\prfd{F_{1,s}(L,2)=f_s(L,{\bf e}_1-\E{1}{n}).}

We now continue this procedure in a recursive fashion. 
We  construct
$F_{r(b),s}(L,b)$ for all $b$ by defining $F_{r(b),s}(L,b)$  in terms
of $F_{r(b-1),s}(L,b-1)$ and $F_{r(b-2),s}(L,b-2)$ 
through the bosonic recursion relations
\brecrelt~for the $b\rightarrow r$ map defined in \rbgen~and then 
simplifying the expressions by using the fermionic recursive properties
of sec. 5 and 6. The $b\rightarrow r$ map
is important since it prescribes which bosonic recursion relation
is being used to construct $F_{r(b),s}(L,b)$ (i.e. 
whether one uses the depth-2
recurrence or one of the depth-1 recurrences).

This process is continued until we reach $b=p'-1$ for which the last
identity of ~\brecreltd~must hold. This recursive construction can be
carried out for any starting function $F_{1,s}(L,1)$ but the last
equation of ~\brecreltd~will only hold if $F_{1,s}(L,1)$ has been
properly chosen. The recursive process used to generate fermionic
polynomials will be referred to as an evolution. The map from two
initial polynomials $F_{1,s}(L,1),~F_{1,s}(L,2)$ to polynomials
$F_{r(b),s}(L,b),~F_{r(b+1),s}(L,b+1)$ will be called a flow of length
$b$. 


After $F_{r(b),s}(L,b)$ has been constructed for all $p'-1$ values of
$b$ to satisfy the bosonic recursion relations ~\brecrelt~--\brecreltd~we
will complete the analysis by studying the behavior for small values of $L$ to
show that with a suitable normalization $q^{{\cal N}_{r,s}}$ the
initial conditions ~\binitial~or ~\newinitial~are satisfied.

For $1\leq b\leq \nu_0+1$ (i.e. $b$ in zone 0) 
this general recursive construction can be
explicitly carried out for an arbitrary $n+1$ zone
problem. In this case the Takahashi decomposition ~\becopm~ of $b$
consists of the  single term
\eqn\prfe{b=l^{(0)}_{1+j_0}=1+j_0,~~~0\leq j_0 \leq \nu_0.}
As long as $1\leq j_0 \leq \nu_0$  $F_{1,s}(L,b)$ satisfies the
first equation of ~\brecrelt~because $r=1$ does not change. 
Comparing this recursion relation with the fermionic recursive
property ~\freczero~and using the values for $F_{1,s}(L,1)$ and
$F_{1,s}(L,2)$ of ~\prfb~and ~\prfd~we conclude that for an $n+1$
zone problem with

{\bf b in zone 0}
\eqn\prff{F_{1,s}(L,b)=f_s(L,{\bf e}_{j_0}-\E{1}{n}),
~~~b=l^{(0)}_{1+j_0}=1+j_0~~ {\rm and}~~ 0\leq j_0 \leq \nu_0}

When $b=2+\nu_0=l^{(1)}_{2+\nu_0}$ we have entered into zone
one. However we still have $r=1$  and an
identical computation gives
\eqn\zoneone{F_{1,s}(L,2+\nu_0)=
f_s(L,{\bf e}_{1+\nu_0}-\E{2}{n}).}

To proceed further into zone 1 we must cross a  boundary where $r$
changes from 1 to 2. Here we use second bosonic recursion relation 
in ~\brecrelt~ for
$b=2+\nu_0$ which has the Takahashi expansion ~\becopm~ with $\beta
=1$ and $\m_1=1,~j_1=1+\nu_0$
\eqn\prfg{b=2+\nu_0=l_{2+\nu_0}^{(1)}}
and find
\eqn\prfh{q^{L-1\over
2}F_{2,s}(L-1,\nu_0+3)=F_{1,s}(L,\nu_0+2)-F_{1,s}(L-1,\nu_0+1)}
which, after using~\zoneone~and \prff~becomes
\eqn\prfi{q^{L-1\over 2}F_{2,s}(L-1,\nu_0+3)=f_s(L,{\bf
e}_{1+\nu_0}-\E{2}{n})-f_s(L-1,-\E{2}{n}).}

To reduce this we replace $f_s(L,{\bf e}_{1+\nu_0}-\E{1}{n})$ 
by an expression in terms of fermionic polynomials with arguments 
 $L-1$ and $L-2$ using relation ~\frecrelb~with $\m=1.$ 
In this case several of the terms in
the general expression ~\frecrelb~vanish and we have
\eqn\prfj{\eqalign{f_s(L&,{\bf e}_{1+\nu_0}-\E{2}{n})
=q^{{L-1\over 2}-{\nu_0-2\over 4}}f_s(L-1,{\bf
e}_{2+\nu_0}-\E{1}{n})+f_s(L-1,-\E{2}{n})\cr
&+f_s(L-1,{\bf e}_{-1+\nu_0}+{\bf
e}_{1+\nu_0}-\E{1}{n})+(q^{L-1}-1)f_s(L-2,{\bf e}_{1+\nu_0}-\E{2}{n}).\cr}}
Using this in ~\prfi~and setting $L\rightarrow L+1$ we find
\eqn\prfk{\eqalign{&F_{2,s}(L,\nu_0+3)=q^{-{\nu_0-2\over 4}}f_s(L,{\bf
e}_{2+\nu_0}-\E{1}{n})\cr
&+q^{-{L\over 2}}f_s(L,{\bf e}_{-1+\nu_0}+{\bf
e}_{1+\nu_0}-\E{1}{n})+(q^{L\over 2}-q^{-{L\over
2}})f_s(L-1,{\bf e}_{1+\nu_0}-\E{2}{n})\cr
&~~~~~~~~=q^{-{\nu_0-2\over 4}}
f_s(L,{\bf e}_{2+\nu_0}-\E{1}{n})+{\tilde f}_s(L,\e{-1+\nu_0}
+{\bf e}_{1+\nu_0}-\E{1}{n})\cr}}
where in the last line we have used the case $j_0=0$ in the definition 
\ftildedf.

In a similar fashion we consider
\eqn\prfl{b=3+\nu_0=l^{(1)}_{2+\nu_0}+l^{(0)}_1.}
We use the third recursion relation
in ~\brecrelt~written as
\eqn\prfm{F_{2,s}(L-1,\nu_0+4)=F_{2,s}(L,\nu_0+3)-q^{L\over
2}F_{1,s}(L-1,\nu_0+2)}
which, after using ~\zoneone~and \prfk~becomes
\eqn\prfn{\eqalign{F_{2,s}(L-1&,\nu_0+4)=q^{-{\nu_0-2\over 4}}f_s(L,{\bf
e}_{2+\nu_0}-\E{1}{n})\cr
&+q^{-{L\over 2}}[f_s(L,{\bf e}_{-1+\nu_0}+{\bf
e}_{1+\nu_0}-\E{1}{n})
-f_s(L-1,{\bf e}_{1+\nu_0}-\E{2}{n})].\cr}}
We now reduce this by using the first equation in ~\freczero, the first
line in the definition ~\ftildedf~and setting $L\rightarrow L+1$ to
obtain the result
\eqn\prfo{F_{2,s}(L,\nu_0+4)=q^{-{\nu_0-2\over 4}}f_s(L,{\bf e}_{1}+{\bf
e}_{2+\nu_0}-\E{1}{n})+{\tilde f}_s(L,\e{\nu_0-2}
+{\bf e}_{1+\nu_0}-\E{1}{n}).}

Let us review what has been done. We started with $F_{1,s}(L,1)$ and
constructed all the values of $F_{1,s}(L,b)$ with $1\leq b \leq 2+\nu_0$
to satisfy the first bosonic recursion relation in ~\brecrelt~where
$r$ does not change. We refer to these recursion relations where $r$
does not change as ``moving $b$ on the plateau''. In this process we
{\it did not create any new terms} in the linear combination.

We then constructed $F_{2,s}(L,3+\nu_0)$ by using the second recursion
relation in ~\brecrelt~in which the term with $b+1$ has $r+1.$ In this
step we created the new term $\tilde f.$ We refer to this new term as
a reflected term (much as there is a reflected wave at a boundary in an
optics problem). We then created $F_{2,s}(L,4+\nu_0)$ by using the third
recursion relation in ~\brecrelt~where $b-1$ has $r-1.$ This process
did not create any additional terms. We refer to this process of using
the two equations in ~\brecrelt~ where $r$ changes by one as ``transiting an
$r$ boundary''. 
In most cases when $b$ moves on the plateau we apply \freczero~to the
$f_s$--terms and \reccft~to the $\tilde{f}_s$--terms in \polyform. Note
that the fermionic recurrences we employ may still vary from term to term
in \polyform. Again in most cases while transiting the $r$--boundary
we use \frecrelb--\frecrele. However, there are important exceptional
cases (related to the so--called dissynchronization effect  discussed
in sec. 8 and appendix B) where this rule breaks down and, as a result,
we are forced to apply \freczero~to some terms in \polyform~and \frecrelb--
\frecrele~to others (for examples of this phenomena see (B.2), (B.8) and 
(B.13)).


With this overview in mind we can continue the construction process as
far as we like. To be explicit we present the results of the
construction for $b$ in zones 1 and 2 in a problem with four or more zones

{\bf b in zone 1}

\eqn\tzcb{\eqalign{1:~~&b=l^{(1)}_{1+j_1}=(j_1-t_1)(\nu_0+1)+1,
~~t_1+1\leq j_1\leq t_2\cr
&r={\tilde l}^{(1)}_{1+j_1}=j_1-\nu_0\cr
&F_{r,s}(L,b)=q^{c(j_1)}f_s(L,\e{j_1}-\E{2}{n})\cr
&~~\cr
2:~~&b=l^{(1)}_{1+j_1}+l^{(0)}_{1+j_0}=l^{(1)}_{1+j_1}+1+j_0,
~~~0\leq j_0\leq t_1-1,~~t_1+1\leq
j_1\leq t_2-1\cr
&r={\tilde l}^{(1)}_{1+j_1}+1=j_1-\nu_0+1\cr
&F_{r,s}(L,b)=q^{c(j_1)-{\nu_0-2\over 4}}f_s(L,\e{j_0}+\e{j_1+1}-\E{1}{n})
+q^{c(j_1)}{\tilde f}_s(L,\e{\nu_0-j_0-1}+\e{j_1}-\E{1}{n})\cr
&~~\cr
3:~~&b=l^{(1)}_{1+t_2}+l^{(0)}_{1+j_0}=\nu_1(\nu_0+1)+2+j_0,~~0\leq j_0\leq t_1-1\cr
&r={\tilde l}^{(1)}_{1+t_2}+1=\nu_1+2\cr
&F_{r,s}(L,b)=q^{c(t_2)-{\nu_0-1\over 4}}f_s(L,\e{j_0}
-\e{t_1}+\e{1+t_2}-\E{3}{n})\cr
&~~~~~~~~+q^{c(t_2)}\tilde{f}_s(L,\e{\nu_0-j_0-1}+\e{t_2}-\E{1}{n})}}
where 
\eqn\prfq{c(j_1)=-{1\over 4}(\nu_0-2)(j_1-\nu_0-1)~~{\rm
for}~~t_1+1\leq j_1\leq t_2.}

{\bf b in zone 2}

Here we distinguish separate cases depending on whether $1\leq j_1\leq t_2-2$
or $j_1=t_2-1.$ The restriction~\restriction~
says that $j_1=t_2$ does
not occur. We also do not consider the cases $b>y_3$
\eqn\tzcc{\eqalign{1:~~&b=l^{(2)}_{1+j_2}=y_1+(j_2-t_2)y_2,~~1+t_2\leq j_2
\leq t_3,\cr
&r={\tilde l}^{(2)}_{1+j_2}=1+(j_2-t_2)\nu_1\cr
&F_{r,s}(L,b)=q^{c(j_2)}f_s(L,\e{j_2}-\E{3}{n}) \cr
&~~\cr
2:~~&b=l^{(2)}_{1+j_2}+l^{(1)}_{1+j_1}+l^{(0)}_{1+j_0}=l^{(2)}_{1+j_2}
+(j_1-t_1)y_1+2+j_0\cr
&0\leq j_0\leq t_1-1,~t_1+1\leq j_1\leq t_2-2,~t_2+1\leq j_2\leq t_3-1,\cr
&r={\tilde l}^{(2)}_{1+j_2}+{\tilde l}^{(1)}_{1+j_1}+1=
{\tilde l}^{(2)}_{1+j_2}+(j_1-t_1)+1\cr
&F_{r,s}(L,b)=q^{c(j_2)+c(j_1)-{\nu_0\over 2}+{1\over 2}}
f_s(L,\e{1+j_0}+\e{t_2-(j_1-t_1)-1}+\e{j_2}-\E{1}{n})\cr
&~~+q^{c(j_2)+c(j_1)-{\nu_0-2\over 4}}\tilde{f}_s(L,\e{\nu_0-j_0-2}
+\e{t_2-(j_1-t_1)}+\e{j_2}-\E{1}{n})\cr
&~~+q^{c(j_2)+c(j_1)-{\nu_0\over
2}+{1\over 4}}f_s(L,\e{j_0}+\e{j_1+1}+\e{j_2+1}-\E{1}{n})\cr
&~~+q^{c(j_2)+c(j_1)-{\nu_0+1\over 4}}
\tilde{f}_s(L,\e{\nu_0-j_0-1}+\e{j_1}+\e{j_2+1}-\E{1}{n})\cr
&~~\cr
3:~~&b=l^{(2)}_{1+j_2}+l^{(0)}_{1+j_0}=l^{(2)}_{1+j_2}+1+j_0,
~0\leq j_0\leq t_1,~~t_2+1\leq j_2
\leq t_3-1\cr
&r={\tilde l}^{(2)}_{1+j_2}+1\cr
&F_{r,s}(L,b)=q^{c(j_2)-{\nu_0\over 4}}f_s(L,\e{j_0}+\e{t_2-1}
+\e{j_2}-\E{1}{n})\cr
&~~+q^{c(j_2)}\tilde{f}_s(L,\e{\nu_0-j_0-1}-\e{t_1}+\e{j_2}-\E{3}{n})\cr
&~~+q^{c(j_2)-{\nu_0+1\over 4}}
f_s(L,\e{j_0}+\e{j_2+1}-\E{1}{n})\cr
&~~\cr
4:~~&b=l^{(2)}_{1+j_2}+l^{(1)}_{1+j_1}=
l^{(2)}_{1+j_2}+(j_1-t_1)y_1+1,
t_1+1\leq j_1\leq t_2-2,~t_2+1\leq j_2 \leq t_3-1\cr
&r={\tilde l}^{(2)}_{1+j_2}+{\tilde l}^{(1)}_{1+j_1}=
{\tilde l}^{(2)}_{1+j_2}+(j_1-t_1)\cr
&F_{r,s}(L,b)=q^{c(j_2)+c(j_1)-{\nu_0\over 4}}
f_s(L,\e{t_1-1}+\e{t_2-(j_1-t_1)}+\e{j_2}-\E{1}{n})\cr
&~~+q^{c(j_2)+c(j_1)-{\nu_0\over 2}+{L\over 2}}
f_s(L,\e{t_2-(j_1-t_1)-1}+\e{j_2}-\E{1}{n})\cr
&~~+q^{c(j_2)+c(j_1)-{\nu_0+1\over 4}}
f_s(L,\e{j_1}+\e{j_2+1}-\E{2}{n})}}

\eqn\tzcd{\eqalign{
5:~~&b=l^{(2)}_{1+j_2}+l_{t_2}^{(1)}=l^{(2)}_{1+j_2}+(\nu_1-1)y_1+1
,~~t_2+1\leq j_2 \leq t_3-1\cr
&r={\tilde l}^{(2)}_{1+j_2}+{\tilde l}^{(1)}_{t_2}={\tilde
l}_{1+j_2}^{(2)}+\nu_1-1\cr
&F_{r,s}(L,b)=q^{c(j_2)+c(t_2-1)-{\nu_0\over 4}}
f_s(L,\e{j_2}-\E{2}{n})\cr
&~~+q^{c(j_2)+c(t_2-1)-{\nu_0+1\over 4}}f_s(L,\e{-1+t_2}+\e{1+j_2}-\E{2}{n})\cr
&~~\cr
6:~~&l^{(2)}_{1+j_2}+l^{(1)}_{t_2}+l^{(0)}_{1+j_0}=l^{(2)}_{1+j_2}+
(\nu_1-1)y_1+2+j_0,\cr
&0\leq j_0 \leq t_1-1,~t_2+1\leq j_2 \leq t_3-1\cr
&r={\tilde l}^{(2)}_{1+j_2}+{\tilde l}^{(1)}_{t_2}+1=
{\tilde l} ^{(2)}_{1+j_2}+\nu_1\cr
&F_{r,s}(L,b)=q^{c(j_2)+c(t_2-1)-{\nu_0\over 4}}
\tilde{f}_s(L,\e{\nu_0-j_0-1}+\e{j_2}-\E{1}{n})\cr
&~~+q^{c(j_2)+c(t_2-1)-{\nu_0\over 2}+{1\over 4}}
f_s(L,\e{j_0}+\e{t_2}+\e{j_2+1}-\E{1}{n})\cr
&~~+q^{c(j_2)+c(t_2-1)-{\nu_0+1\over 4}}
\tilde{f}_s(L,\e{\nu_0-j_0-1}+\e{t_2-1}+\e{j_2+1}-\E{1}{n})\cr
&\cr}}

where
\eqn\ctdefn{c(j_2)={1\over 2}-{3+\nu_1(\nu_0-2)\over 4}(j_2-t_2)~~~{\rm
for}~~t_2+1\leq j_2 \leq t_3.}

Thus far we have used our constructive procedure to generate all
polynomials $F_{r(b),s}(L,b)$ for $b$ in zones 0, 1 and 
2 where the total number of
zones is 4 or greater. However, to complete the process we must carry
out the construction for $b$ in the final zone and show that the closing 
relation in \brecreltd~is satisfied. The construction for 
$b> l^{(n)}_{1+t_{1+n}}$ has two
features not present in any other zone. The first is that the map
$b\rightarrow r$ of ~\rbgen~has changed and the second is that the
parity restriction on $m_{t_{1+n}},$ which was even before, 
is now sometimes allowed to be odd. We recall that this parity
is specified in our notation by the parity of the $1+t_{1+n}$
component of the vector $\bf u.$ Thus we compute the
following results for the final zone.

{\bf b in zone 1 in a two zone problem}

1:~~The first equation of ~\tzcb~now may be extended to
all $1+t_1\leq j_1\leq 1+t_2$
with $-\E{2}{n}$ replaced by zero and $c(t_2+1)$ given by ~\prfq.

2:~~The second  equation of ~\tzcb~now holds for all $1+t_1\leq j_1\leq t_2$
with $-\E{2}{n}$ replaced by zero.

3:~~Equation three of ~\tzcb~is omitted and replaced by

\eqn\tzcl{\eqalign{b&=
l^{(1)}_{2+t_2}+l^{(0)}_{1+j_0}=
(\nu_1+1)(\nu_0+1)+2+j_0,~~0\leq j_0\leq t_1-1\cr
r&={\tilde l}^{(1)}_{2+t_2}= \nu_1+1=p-1\cr
F_{r,s}(L,b)&=q^{{ c}(t_2+1)}f_s(L,{\bf e}_{t_1-j_0-1}-
{\bf e}_{t_1}+{\bf e}_{1+t_2})\cr}}
where ${c}(t_2+1)$ is given by \prfq~instead of~\ctdefn. 

{\bf b in zone 2 in a three zone problem}

1:~~Equation 1 in ~\tzcc~ is now valid for $1+t_2 \leq j_2 \leq t_3+1$ and 
equations 2-6 in ~\tzcc~ and ~\tzcd~are now valid for
$t_2+1\leq j_2 \leq t_3+1$ with the convention that 
wherever ${\bf e}_{2+t_3}$
appears in the argument of some $f_s(L,\vu)$ or 
${\tilde f}_s(L,{\tilde\vu})$ the term is omitted~
and $c(t_3+1)$ is given by \ctdefn.

2:~~We have the following additional closing equation

\eqn\closing{\eqalign{&b=l^{(2)}_{2+t_3}+l^{(1)}_{t_2}+l^{(0)}_{1+j_0}=
p'-\nu_0+j_0-1\cr
&0\leq j_0\leq t_1\cr
&r={\tilde l}^{(2)}_{2+t_3}+{\tilde l}^{(1)}_{t_2}=\nu_1(\nu_2+2)=p-1\cr
&F_{r,s}(L,b)=q^{c(t_3+1)+c(t_2)-{1\over 2}}
f_s(L,\e{\nu_0-j_0}-\e{t_1}-\e{t_2}+\e{1+t_3})\cr}}
where $c(t_3+1)$ is given by ~\ctdefn.

For the problem with 2 and 3 zones we have now constructed a complete
set of fermionic polynomials for all $1\leq b \leq p'-1$ which satisfy
the bosonic recursion relations~\brecrelt, 
\brecreli~and the first equation in~\brecreltd~by construction and the
second equation in ~\brecreltd~by use of the first equation
of ~\freczero. In order to complete the
proof of the bose/fermi identities it remains to show that the
fermionic polynomials satisfy the boundary conditions for
$B_{r(b),s}(L,b)$ at $L=0$~\binitial~and to compute the normalization
constant in \eidbf. This is easily done and thus we obtain the final
result that for the cases of 2 and 3 zones with $s=l^{(\m_s)}_{1+j_{s}}$ 
\eqn\finalres{F_{r(b),s}(L,b)=q^{{1\over
2}({\phi}_{r(b),s}-\phi_{r(s),s})+c(j_s)} B_{r(b),s}(L,b)}
where $F_{r(b),s}(L,b)$ is given by the formulae of this section.

\newsec{The inductive analysis of evolution}

In the previous section we constructed all fermionic polynomials
$F_{r(b),s}(L,b),~1\leq b \leq p'-1$ in the problem with 1,2 and 3 zones. 
The procedure we used is completely general, but
becomes somewhat tedious to execute when  the number of
zones increases beyond three. This is because the results
are sensitive to the details of the Takahashi decomposition ~\becopm.
On the other hand, there are certain cases of ~\becopm~such as
$b=l^{(\mu)}_{1+j_{\mu}},~1+t_{\mu}\leq j_{\mu}\leq
t_{1+\mu}+\delta_{\mu,n},~0\leq \mu \leq n$ where the form of
$F_{r(b),s}(L,b)$ remains very simple for any $\mu.$ The question
arises if one can treat certain classes of decompositions
~\becopm~without having explicit formulas for all fermionic polynomials  
$F_{r(b),s}(L,b),~1\leq b\leq p'-1.$

In this section we shall provide a positive answer to this question by
proving inductively a set of explicit formulas for certain flows of
length $y_{\mu}-1.$ This inductive analysis is possible because as can
be seen from ~\rpropiv~the construction of
$F_{r(b),s}(L,b),~b=l^{(\mu)}_{2+j_{\mu}},~1+t_{\mu}\leq j_{\mu}\leq
t_{1+\mu}-1+\delta_{\mu,n},\mu\geq 2$ from the pair
$\{F_{r(b')+1,s}(L,b'+1),F_{r(b')+1,s}(L,b'+2)\},~b'=l^{(\mu)}_{1+j_{\mu}}$
involves exactly the same steps as that of
$F_{z_{\mu-1},s}(L,y_{\mu})$ from the pair $\{F_{1,s}(L,1),F_{1,s}(L,2)\}.$
Furthermore, recalling that $y_{1+\mu}=y_{\mu-1}+\nu_{\mu}y_{\mu},$ it
is natural to decompose the flow of length $y_{\mu+1}-1$ into flows of
smaller length and to take this decomposition as a starting point of
our inductive analysis of evolution. In this direction, we first
discuss the notation for the flows in terms of the $b\rightarrow r$
map, whose properties are summarized in ~\rpropi-\rpropv.


{\bf Notation}

{\it

The flow $\aw{x}$ of length $x+1$ denotes the sequence 
of steps corresponding to the
$b\rightarrow r$ map as defined in \rbgen~(or equivalently \alterrnew)
between $b=1$ and $b=2+x$.
\eqn\pprba{\psfig{file=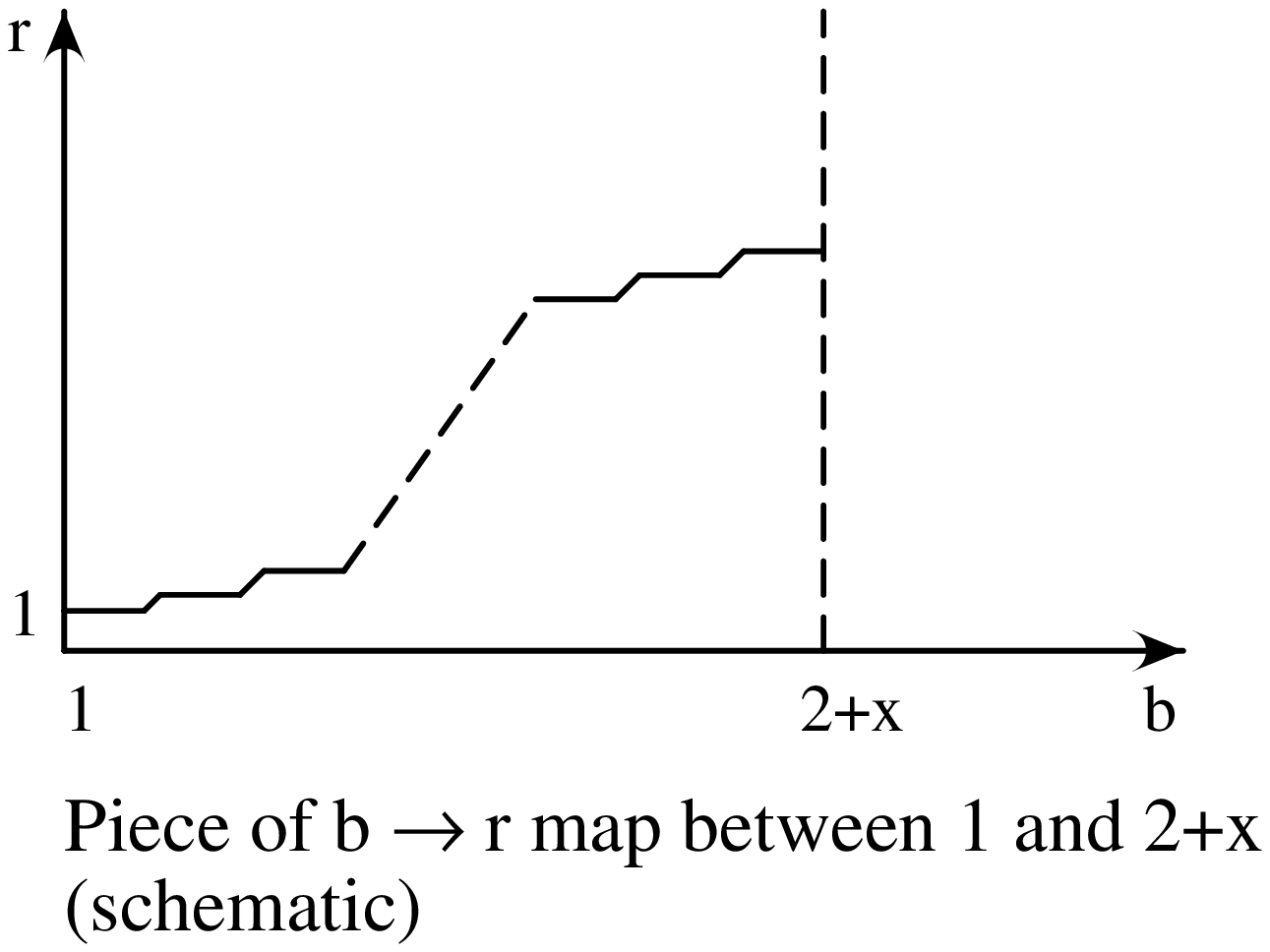,width=2.5in}}
According to ~\rpropi~and ~\rpropii~
this sequence is made up of three steps:
\eqn\pprbb{\psfig{file=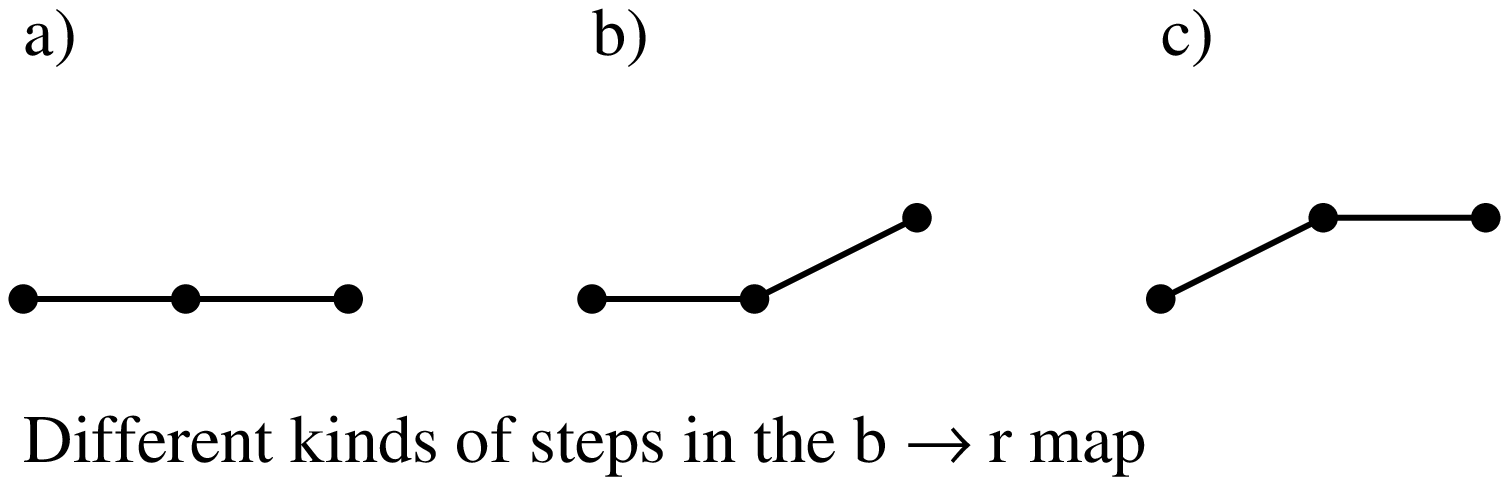,width=3in}}
(where we note that steps like
\eqn\notallow{\psfig{file=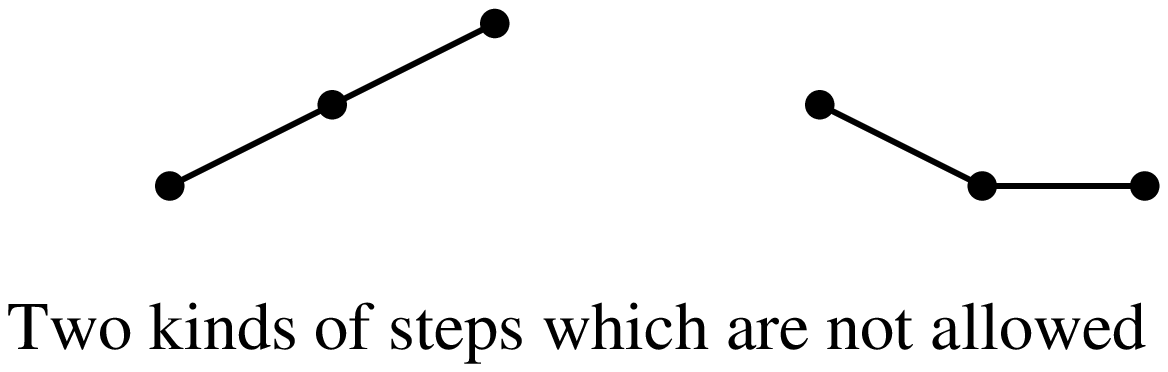,width=3in}}
are not allowed). 
The $b\rightarrow r$ map specifies exactly which bosonic recursion
relation is being used to define the (in this case fermionic)
polynomials.
Case $a)$ in \pprbb~denotes that $r$ is the same for all three objects
$F_{r,s}(L,b-1)$, $F_{r,s}(L,b)$ and $F_{r,s}(L,b+1)$ involved and 
hence one uses
the first recursion relation \brecrelt. Case $b)$ indicates the use
of the second and case $c)$ the use of the last recursion relation in 
\brecrelt.

Let us further denote by $\evs{2}$ the flow according to
the steps as defined by the following graph:
\eqn\pprbs{\psfig{file=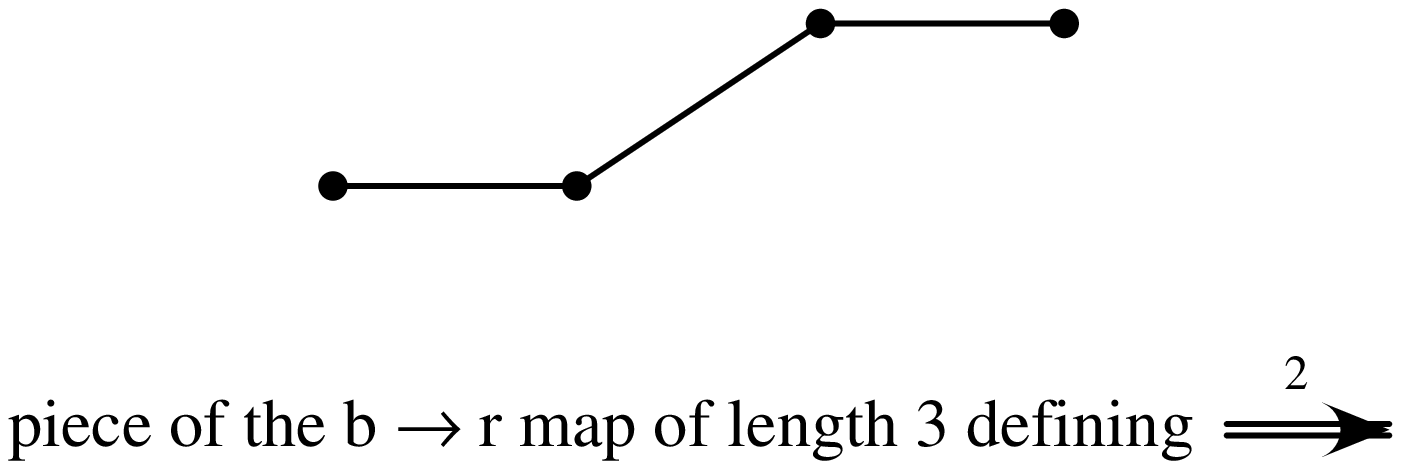,width=3in}}
} 
\bigskip

The $b\rightarrow r$ map as defined by \rbgen~(or \alterrnew)~ 
has an unambiguous initial point
$(r_1,b_1)=(1,1)$. By $\aw{x}$ we mean only the sequence of $x$ steps
found on the $b \rightarrow r$ map as specified in \pprba. 
The notation $\aw{x}$ 
does not fix the initial point $(r_1,b_1)$. The initial point
$(r_1,b_1)$ can be placed anywhere . Having agreed on that, we can now
piece together flow 1 and flow 2 such that the last segment of the
first flow is identical to the first segment of the second flow by
identifying the two final points of the first flow as the two initial
points of the second flow. For example adjoining flow c of \pprbb~ to flow b
of \pprbb~ gives flow \pprbs~ and adjoining a of ~\pprbb~ to c of
~\pprbb~gives
\eqn\newgra{\psfig{file=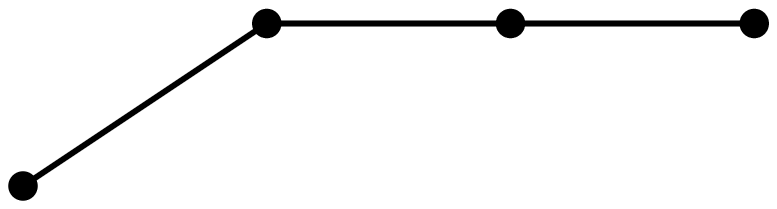,width=1.5in}}

\bigskip

We then denote by $\aw{x_{(1)}}\aw{x_{(2)}}$ the flow given by piecing
together the flow $\aw{x_{(1)}}$ and $\aw{x_{(2)}}$. Note that in
general it is not true that
$\aw{x_{(1)}+x_{(2)}}=\aw{x_{(1)}}\aw{x_{(2)}}.$

With these conventions we now show that there is a 
decomposition of $\aw{y_{\m+1}-2}$
in terms of $\aw{y_{\m-1}-2},~\aw{y_{\m}-2}$ and $\evs{2}$ given 
by the following graph
\eqn\ppdeca{\psfig{file=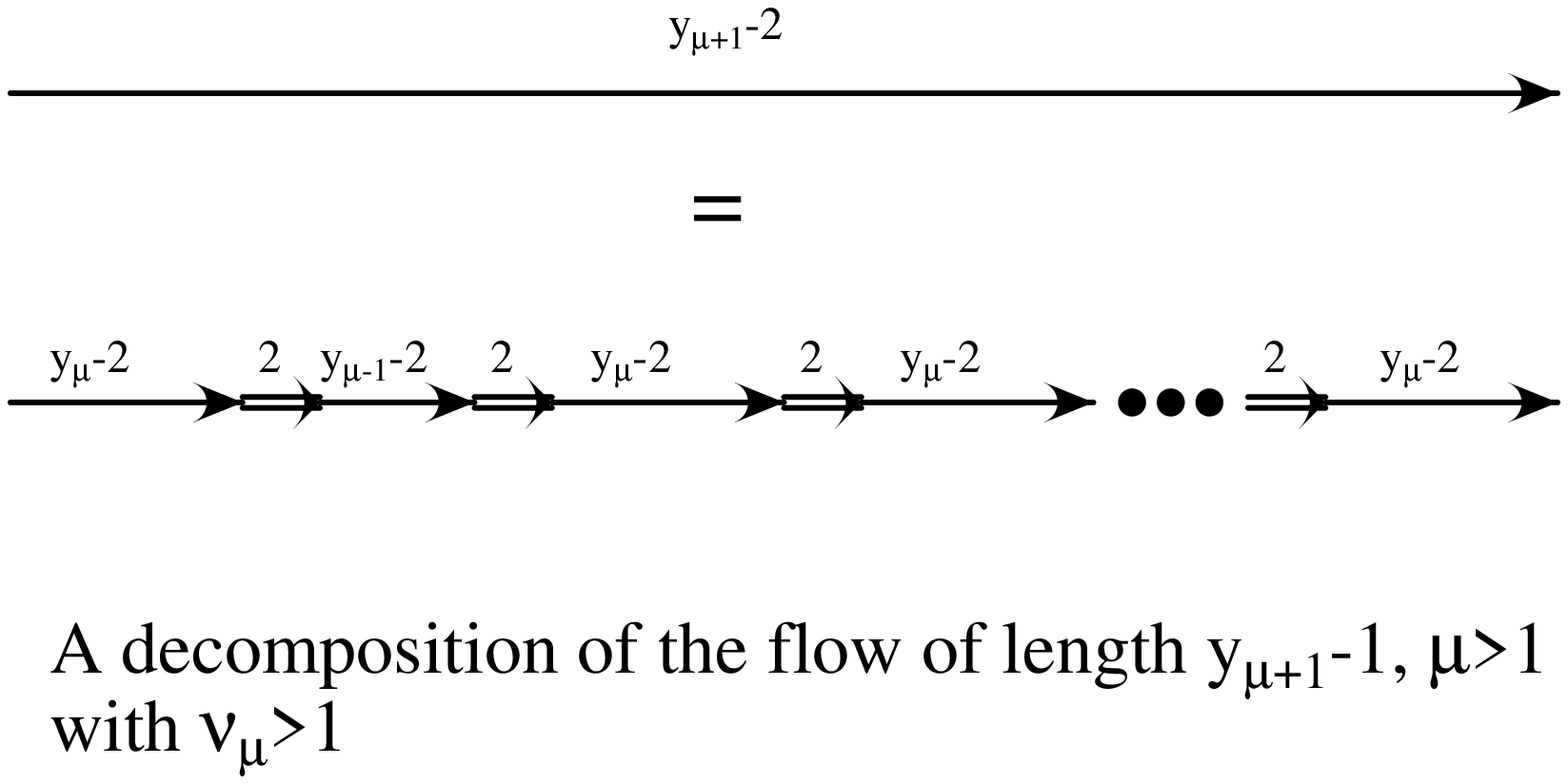,height=1.4in}}

In other words we need to show that piecing 
together $\aw{y_{\m}-2},~\aw{y_{\m-1}-2}$ and $\evs{2}$
as shown in \ppdeca~amounts to $\aw{y_{\m+1}-2}.$ This 
can be easily done by recalling ~\rpropiii~ with $b=y_{\mu}$ 
and  $b=y_{\mu-1}+ky_{\mu}$ for $1\leq k \leq \nu_{\mu}-1,~2\leq \mu$
and proving two additional results: 
\eqn\pcha{r(b+y_{\m})-r(b)=r(y_{\m})=z_{\m-1},~{\rm for}~1\leq b\leq y_{\m-1}}
and
\eqn\pchb{r(b+y_{\m-1}+ky_{\m})-r(b)=r(y_{\m-1}+ky_{\m})=
\tilde{l}^{(\m)}_{1+t_{\m}+k}~{\rm for}~1\leq b\leq y_{\m} }

To check \pcha~we use the Takahashi decomposition of $b$ ~\becopm~with
$\m_\beta \leq \m-2.$
 If $b\not= y_{\m-1}$ then
$b+y_{\m}=\sum_{i=1}^{\beta}l^{(\m_i)}_{1+j_{\m_i}}+l^{(\m-1)}_{1+t_{\m}}$ is 
a valid decomposition in Takahashi lengths and according to \rbgen~we have
$r(y_{\m}+b)-r(b)=\tilde{l}^{(\m-1)}_{1+t_{\m}}=
z_{\m-1}$. Similarly for $b=y_{\m-1}$
we have $r(y_{\m-1}+y_{\m})-r(y_{\m-1})=\tilde{l}^{(\m)}_{2+t_{\m}}
-\tilde{l}^{(\m-2)}_{1+t_{\m-1}}=z_{\m-1}$.
\pchb~may be checked in a similar fashion.

Notice however that the order of the arrows in the
decomposition as shown in the previous figure is crucial. If we 
moved for example $\aw{y_{\mu-1}-2}\evs{2}$ to the one before last
position as shown in the next figure the decomposition does not agree
with all the steps as defined by $\aw{y_{\m+1}-2}$.
\eqn\ppdecb{\psfig{file=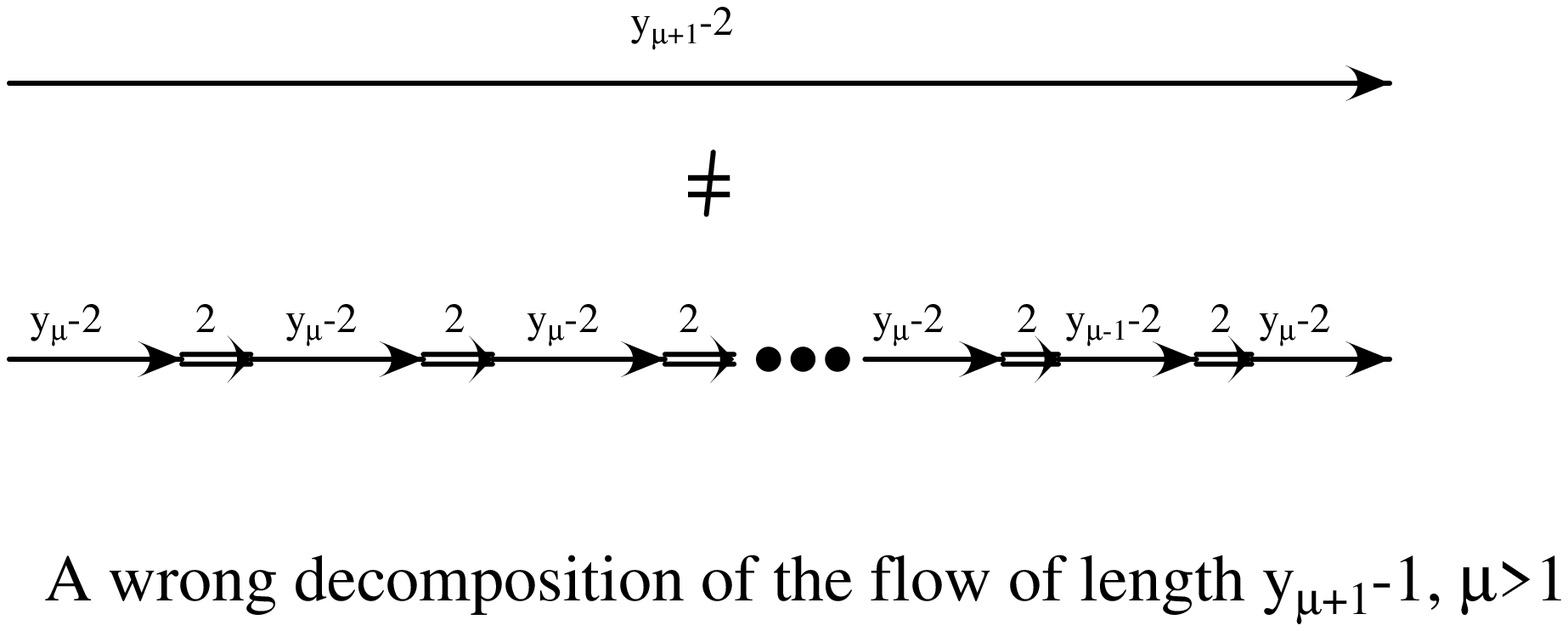,height=1.4in}}

This can be easily seen from  \rpropiii~with
$b=y_{\mu-1}+(\nu_{\mu}-1)y_{\mu}$ and $b=ky_{\mu},
~1\leq k \leq \nu_{\mu}-1,~2\leq \mu$
and the following lemma:

{\bf Lemma 1.1}

{\it
{\bf a.} For $1\leq k\leq \nu_{\m}-2$, $\m\geq 2$ we have 
\eqn\prbt{r(b+ky_{\m})-r(b)=kz_{\m-1}}
for all $1\leq b\leq y_{\m}$ except for 
$b=b_{\m,-1}=\sum_{i=0}^{\m-1}y_i$ for which
\eqn\prbts{r(b_{\m,-1}+ky_{\m})-r(b_{\m,-1})=kz_{\m-1}+(-1)^{\m},}

{\bf b.} for $1\leq b\leq y_{\mu-1}$ we have
\eqn\lemmanewa{r(b+(\nu_{\mu}-1)y_{\mu})-r(b)=(\nu_{\mu}-1)z_{\mu-1},}

{\bf c.} for $1\leq b \leq y_{\mu}$
\eqn\lemmanewb{r(b+l^{(\mu)}_{t_{1+\mu}})-r(b)={\tilde l}^{(\mu)}_{t_{1+\mu}}}
}

We shall also require the companion of ~\prbt~and ~\prbts

{\bf Lemma 1.2}

{\it 
For $\m\geq 3$ we have
\eqn\prbo{r(b+y_{\m-1})-r(b)=z_{\m-2}}
for all $1\leq b\leq y_{\m}$ except for $b=b_{\m,1}=\sum_{i=0}^{\m-2}y_i$ 
for which
\eqn\prbos{r(b_{\m,1}+y_{\m-1})-r(b_{\m,1})=z_{\m-2}+(-1)^{\m-1}}
}

These results can be checked in the same fashion as equations \pcha~and \pchb.

The decomposition ~\ppdecb~fails precisely because there is a special
value $b=b_{\mu,-1}$ for which ~\prbt~is not
valid. However, one can find a slightly modified version of
~\ppdecb~which does indeed hold. For this we need to define two
further types of arrows.

{\bf Definition}

{\it
We denote by $\ev{y_{\m}-2}{-1}$ the flow according to 
the steps defined by the $b\rightarrow r$ map \rbgen~between $1+y_{\m}$ and 
$2y_{\m}$.
\eqn\pprbc{\psfig{file=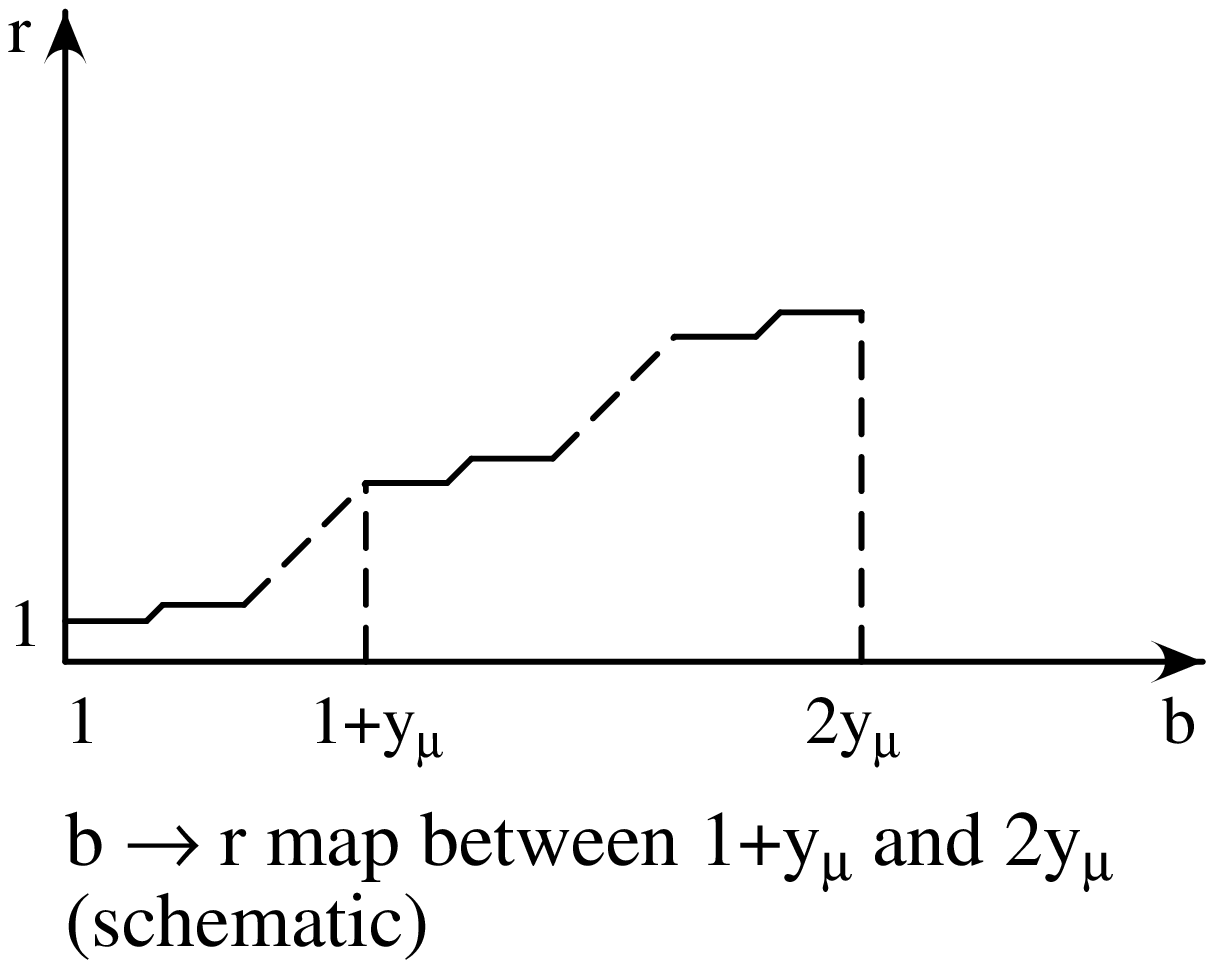,height=2.5in}}

\bigskip

We denote by $\ev{y_{\m}-2}{1}$ the flow according to
the steps defined by the $b\rightarrow r$ map \rbgen~between $1+y_{\m-1}$ and 
$y_{\m-1}+y_{\m}$.
\eqn\pprbd{\psfig{file=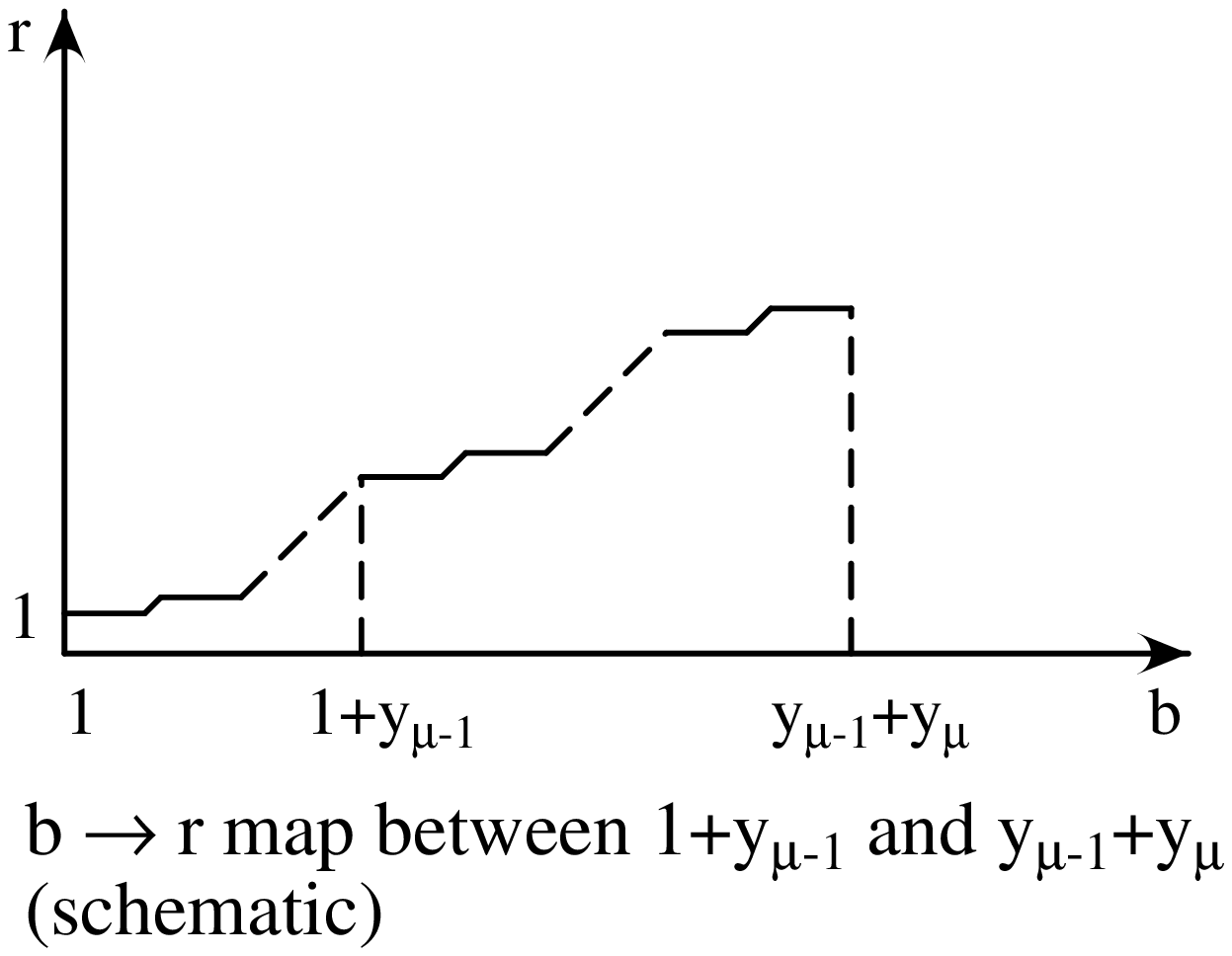,height=2.5in}}

\bigskip

We further identify 
\eqn\onemore{\ev{y_{\m}-2}{0}=\aw{y_{\m}-2}}
}

Slightly generalizing the above discussion we recapitulate.
Let $O_i(L,b_i,q)$ for $i=1,2$ be polynomials (in $q^{1/4})$ 
depending on $L$ where $L$
is a nonnegative integer. Let $O_i(L,b_i,q)$ for $i=3,4$ 
be polynomials 
depending on $L$ obtained recursively from
$O_1(L,b_1,q),~O_2(L,b_2,q)$ 
by the flow
$\ev{x}{t}$ of length $1+x$. We denote this as
\eqn\absflow{\left\{ {O}_1(L,b_1,q)~,{O}_2(L,b_2,q) \right\} \ev{x}{t} 
\left\{ {O}_3(L,b_3,q),~{O}_4(L,b_4,q) \right\}}
where the symbol $x$ above the arrow denotes that
$O_4(L,b_4,q)$ follows from $O_1(L,b_1,q)$ and $O_2(L,b_2,q)$ after $x$ steps. 
Parameters $b_i,~i=1,2,3,4$ with $b_2=b_1+1,~b_3=b_1+x,~b_4=b_1+x+1$
and $q$
associated with $O_i$ will often be suppressed. The symbol $t$
below the arrow  denotes which sequences of
the recursion relations in ~\brecrelt~ are being used. 
The sequence $t$ can be thought of
pictorially as a continuous graph made up of horizontal and diagonal
segments where by horizontal segment we mean only the segment 
between $(r,b)$ and
$(r,b+1)$ and by diagonal segment we mean the segment between $(r,b)$ and
$(r+1,b+1).$ The only restriction on the sequence t is that diagonal
segment must be preceded and followed by horizontal segments unless the
diagonal segment is the first or last segment. As in the case of the
notation $\aw{x},$ the notation $\ev{x}{t}$ does
not fix the initial point. Each pair of adjacent segments, called a
step, in the flow $\ev{x}{t}$ represents one of the three recursion
relations in ~\brecrelt~ in a manner exactly analogous to the discussion
following ~\pprbb. More precisely,
\eqn\miwaa{\{O_1(L,b_1,q),O_2(L,b_1+1,q)\}\ev{1}{t}
\{O_3(L,b_1+1,q),O_4(L,b_1+2,q)\}}
where
\eqn\miwab{O_3(L,b,q)=O_2(L,b,q)}
and
\eqn\miwac{O_4(L,b_1+2,q)=\cases{O_2(L+1,b_1+1,q)-O_1(L,b_1,q)
+(1-q^L)O_2(L-1,b_1+1,q)&~~~\cr
~~~{\rm if}~
t~ {\rm is~ a~ of~ (8.2)}&~~~\cr
q^{-L/2}[O_2(L+1,b_1+1,q)-O_1(L,b_1,q)]&~~~\cr
~~~{\rm if}~ t~ {\rm is~ b~ of~ (8.2)}&~~~\cr
O_2(L+1,b_1+1,q)-q^{(L+1)/2}O_1(L,b_1,q)&~~~\cr
~~~~ {\rm if} ~t~ {\rm is~ c~ of~ (8.2).}&~~~\cr}}
Note that for $L=0$ the last term in the top line of ~\miwac~
vanishes. Therefore one does not need to know the polynomials
$O_1(L,b_1,q),~O_2(L,b_1+1,q)$ for $L<0$ to determine $O_4(L,b_1+2,q)$
for $L \geq 0.$

To compare three different flows~\pprbc--\onemore, let us set 
the initial points of the
flows $\ev{y_{\mu}-2}{\pm 1}$, $\ev{y_{\mu}-2}{0}$ 
to be $(r_1,b_1)= (1,1).$ Then
according to \prbt, \prbts~and \prbo~,\prbos~the pieces of the 
$b\rightarrow r$ map used to define $\ev{y_{\m}-2}{a},~a=\pm 1$ 
differ from the one used  to define
$\ev{y_{\m}-2}{0}$ at exactly one point.
We refer to this phenomena as the dissynchronization effect.
Recalling the definition for $a=\pm 1$ that
$b_{\m,a}=\sum_{i=0}^{\m-1-\delta_{a,1}}y_i$, 
$\m\geq 2+\delta_{a,1}$ we find that the difference 
between
$\ev{y_{\m}-2}{a},~a=\pm 1$ and $\ev{y_{\m}-2}{0}$ can be illustrated by the
following two figures. Here the solid line denotes the piece of 
the $b\to r$ map for 
$\ev{y_{\m}-2}{a},~a=\pm 1$ and the dashed line that for $\ev{y_{\m}-2}{0}$.
\eqn\ppdis{\psfig{file=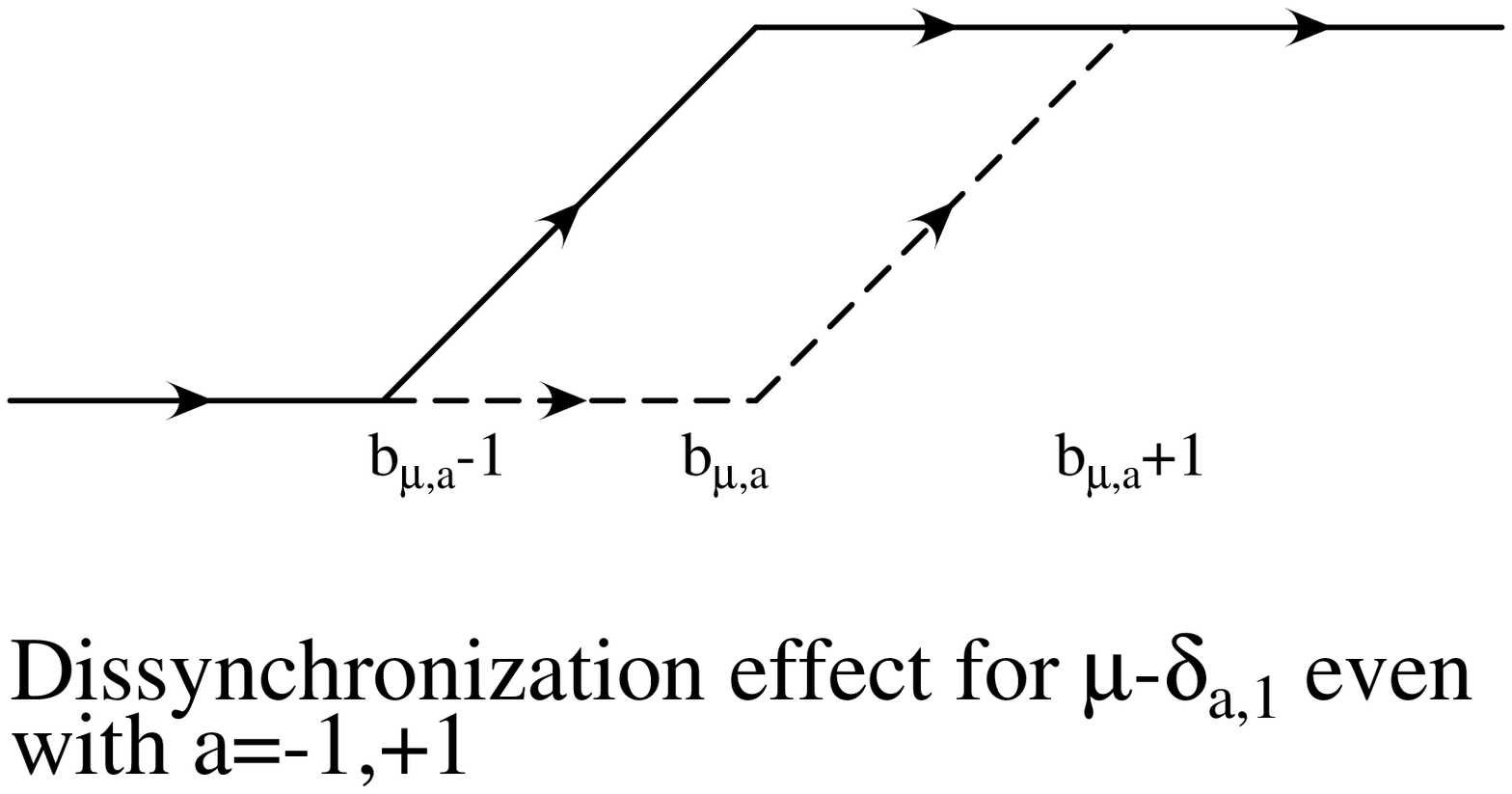,width=2.5in} \psfig{file=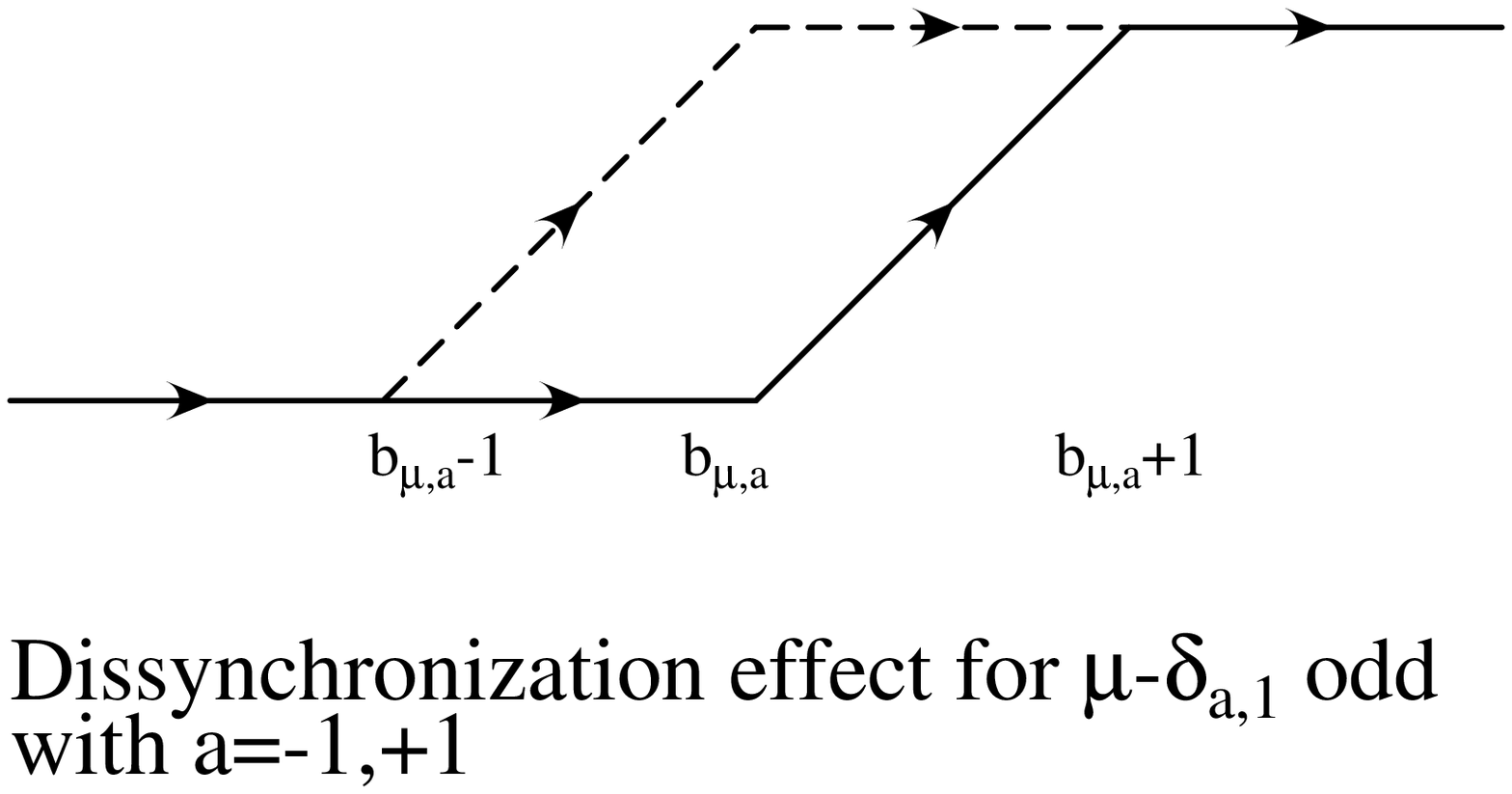,width=2.5in}}

More precisely, all elementary segments in these two flows are
identical   except in the interval $(b_{\mu,a}-1,b_{\mu,a}+1).$ For
$\mu-\delta_{a,1} $
even (respect. odd) the flow $\ev{y_{\mu}-2}{a},~a=\pm 1$ restricted to the
interval $[b_{\mu,a}-1,b_{\mu,a}+1]$ is given by c of ~\pprbb~(b of
~\pprbb). The flow, $\ev{y_{\mu}-2}{0},$ restricted to the same
interval is given by b of ~\pprbb~(c of ~\pprbb). 

We now may prove  the following two decompositions of
$\ev{y_{\m+1}-2}{a}$ where $a=-1,0,1.$ 
\eqn\pppropa{\psfig{file=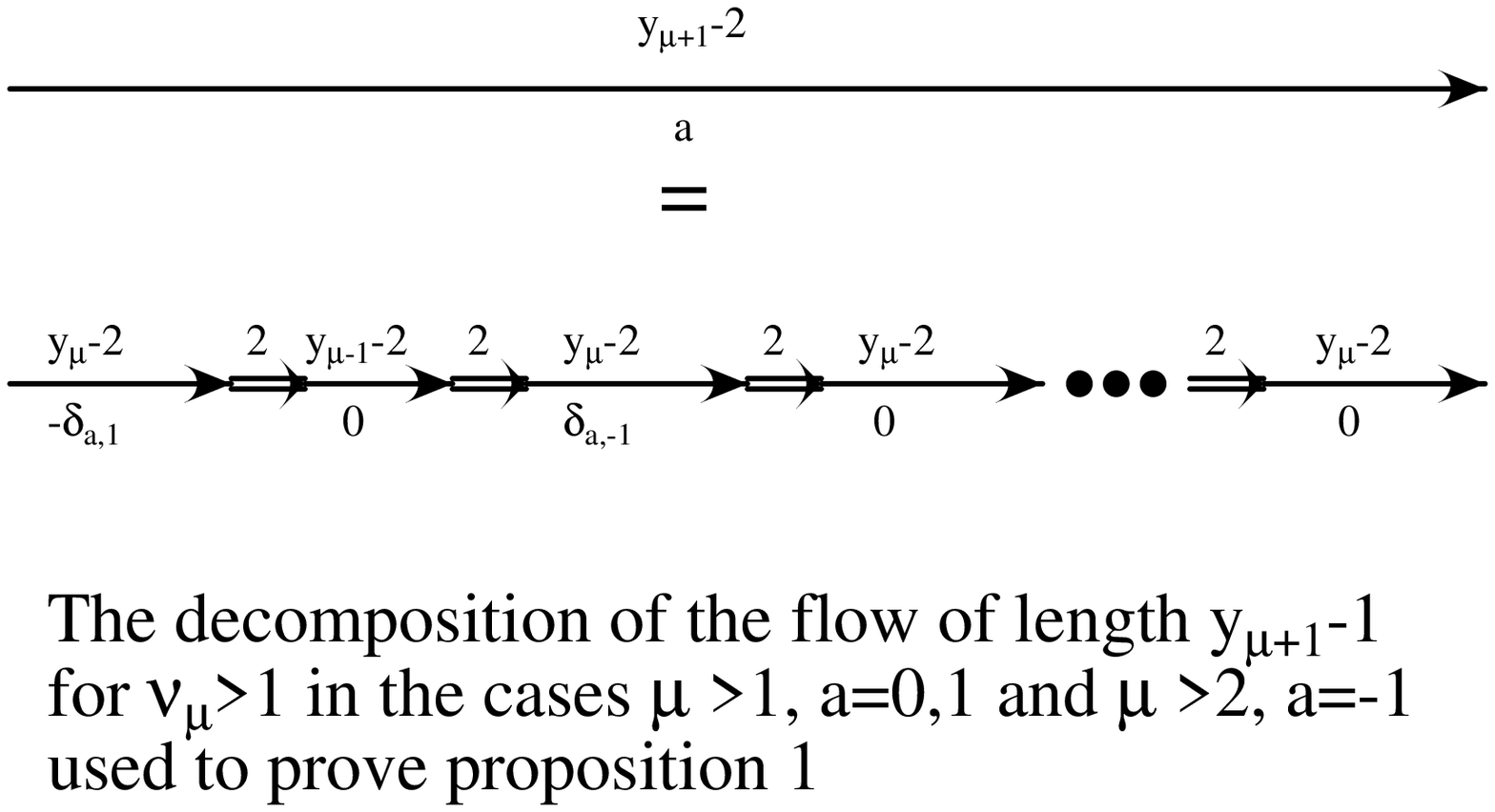,height=2.0in}}
and 
\eqn\pppropb{\psfig{file=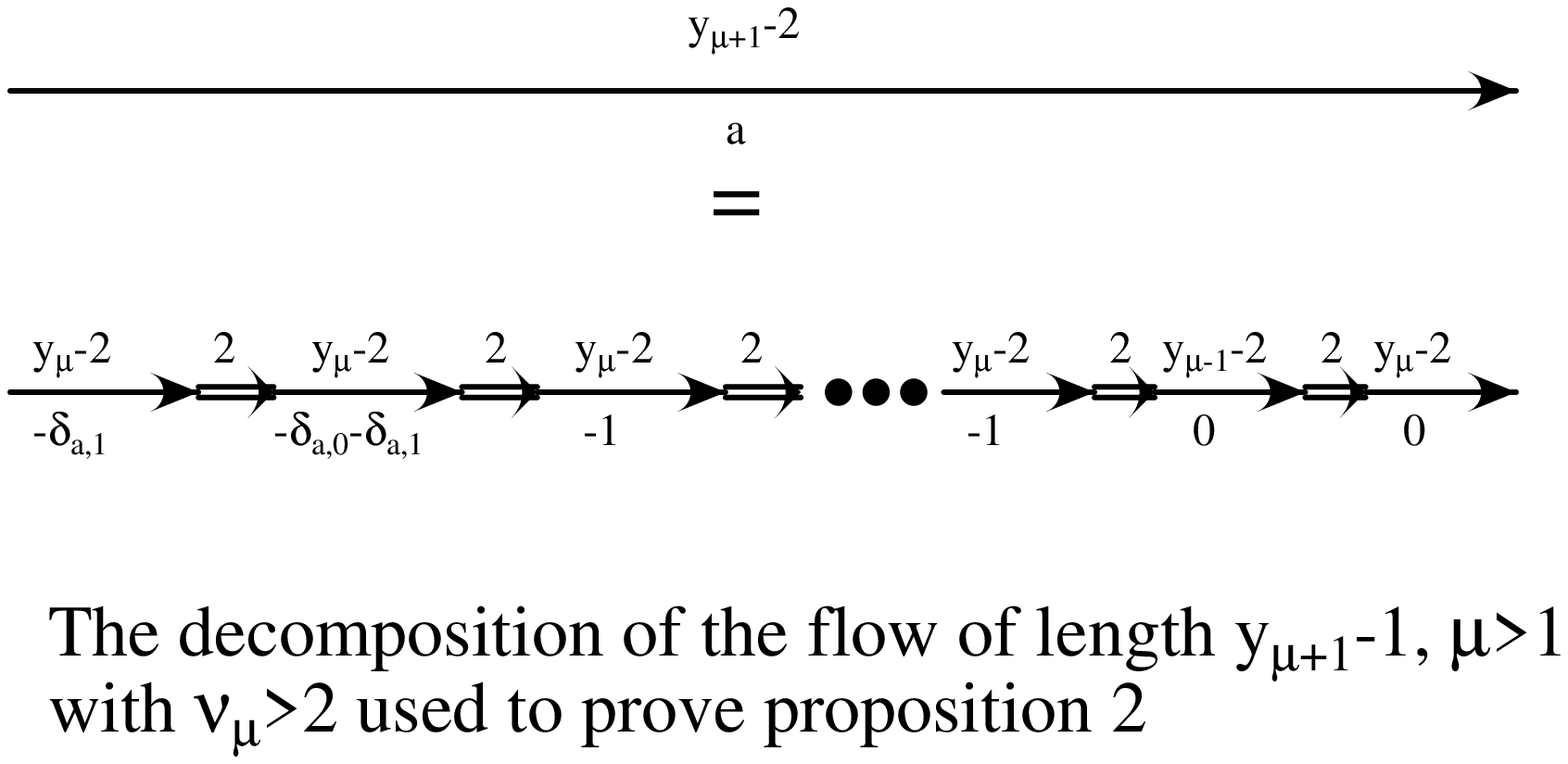,height=1.4in}}
In ~\pppropa~the arrows $\evs{2}\ev{y_{\mu}-2}{0}$ appear
$\nu_{\mu}-2$ times and in ~\pppropb~the arrows
$\evs{2}\ev{y_{\mu}-2}{-1}$
appear $\nu_{\mu}-3$ times.

The decomposition ~\pppropa~for $a=0$  is ~\ppdeca. To
prove~\pppropa~for $a=\pm 1$ we use~\pppropa~  with $a=0$ and~\ppdis~  with
$a=\pm 1,$ $\mu$ replaced by $\mu+1$ to convert the top graph of
~\pppropa~with $a=0$ into the top graph of~\pppropa~  with $a=\pm 1,$ by
altering just the single step centered at $b_{\mu+1,a},~a=\pm 1$ as
prescribed by~\ppdis. It is easy to see from~\ppdis~ that this 
procedure converts
the first (fifth) arrow of the bottom graph~\pppropa~  with $a=0$ into
$\ev{y_{{\mu}-2}}{-1} (\ev{y_{{\mu}-2}}{1})$ for $a=1(-1)$ and has no
effect on other arrows of the bottom graph ~\pppropa. This completes the
proof of~\pppropa. 

The decomposition ~\pppropb~is  the correct modification of
~\ppdecb. When $a=0$ it follows from \prbt-\lemmanewb~and 
~\ppdis~with $a=-1$. To prove ~\pppropb~for
$a=1(-1)$ we again replace one single step of the top graph of
~\pppropb~with $a=0.$ This converts the first (third) arrow of the
bottom graph into $\ev{y_{\mu}-2}{-1} (\ev {y_{\mu}-2}{0})$ and thus
completes the proof.

Let us finally define some useful sets of terms

{\bf Definition}

{\it For $1\leq \m\leq n,$
$P_{-1}(L,\m,\u{\m})$, $P_0(L,\m,\u{\m})$ and $P_1(L,\m,\u{\m})$
are defined as
follows
\eqn\ppmsl{\eqalign{P_{-1}(L,\m,\u{\m})=&
\sum_{i=2}^{\m}q^{{L\over 2}-{\nu_0-\kio{i}
\over 4}} f_s(L,-\E{1}{i}+\e{-1+t_{i}}+\u{\m})
+f_s(L,\e{-1+t_1}-\e{t_1}+\u{\m})\cr
=&\sum_{i=2}^{\m} q^{{\nu_0-1-\kie{i}
\over 4}} \tilde{f}_s(L,-\E{1}{i}+\e{-1+t_{i}}+\u{\m})
+f_s(L,\e{-1+t_1}-\e{t_1}+\u{\m})}}
\eqn\ppms{\eqalign{P_0(L,\m,\u{\m})
=\sum_{i=2}^{\m} q^{-{\nu_0-\kio{i}\over
4}} f_s(L,-\E{1}{i}+\e{-1+t_i}+\u{\m})
+\tilde{f}_s(L,\e{\nu_0-1}-\e{\nu_0}+\u{\m})}}
and
\eqn\ppmsh{\eqalign{P_1(L,\m,\u{\m})
=\sum_{i=2}^{\m} q^{-{\nu_0-\kio{i}\over
4}} f_s(L,\e{1}-\E{1}{i}+\e{-1+t_i}+\u{\m})
+\tilde{f}_s(L,\e{\nu_0-2}-\e{\nu_0}+\u{\m})}}
where $\u{\m}$ is any $1+t_{n+1}$-dimensional vector with 
non-zero entries only in zone $\m$ or higher, i.e. $(\u{\m})_i=0$
for $i\leq t_{\m}$. 
(We consider here the case that $\nu_i\geq 2,~1\leq i\leq 
n-1$).
}

Equipped with the above notations and explanations we are now in the
position to formulate two propositions which will be important in the sequel.

{\bf Proposition 1}

{\it 
Let $a=-1,0,1$ and $P_{-1}(L,\m,\u{\m})$ be defined as in \ppmsl.
Then we have for $a=0$, $\m \geq 1$ and $a=-1$, 
$\m \geq 2$ and $a=1$, $\m\geq 3$
\eqn\pprop{\eqalign{
\left\{ f_s(L,-\E{1}{\m}+\u{\m}),f_s(L,\e{1}-\E{1}{\m}+\u{\m})
 \right\}
\ev{y_{\m}-2}{a} q^{c(t_{\m})+{1\over 2}a(-1)^{\m}}  \left\{
P_{-1}(L,\m,\u{\m}),f_s(L,\u{\m}) \right\}
}}
}

{\bf Proposition 2}

{\it
Let $a=-1,0,1$ and $P_0(L,\m,\u{\m})$ and $P_1(L,\m,\u{\m})$ 
be defined as in \ppms~ and \ppmsh. Then we have for $a=0$,
 $\m\geq 1$ and
$a=-1$, $\m\geq 2$ and $a=1$, $\m\geq 3$
\eqn\ppmsflow{\eqalign{ \left\{ P_0(L,\m,\u{\m}),
P_1(L,\m,\u{\m}) \right\}&\ev{y_{\m}-2}{a} 
\left\{ q^{{L-\nu_0+\kio{\m}\over 2}+
c(t_{\m})+{1\over 2}a(-1)^{\m}}f_s(L,-\E{1}{\m}+\u{\m}),0 \right\}\cr
&= q^{c(t_{\m})-{1\over 2}\kie{\m}+{1\over 2}a(-1)^{\m}}
\left\{ \tilde{f}_s(L,-\E{1}{\m}+\u{\m}),0 \right\}
}}
}

Here $c(t_{\m})$ is defined recursively by
\eqn\pct{c(t_{\m+1})=c(t_{\m-1})
+\nu_{\m} c(t_{\m})-{\nu_0-\kie{\m}\over 4}
+(\nu_{\m}-1)(-{\nu_0+1\over 4}+{3\over 4}\kio{\m})}
where
\eqn\pcin{c(t_0)={\nu_0\over 4},~~~~c(t_1)=0.}
We will also need 
\eqn\pcj{\eqalign{c(j_0)&=0~~{\rm for}~1\leq j_0\leq t_1\cr
c(j_{\m})&={1\over 2}(-)^{\m}+(j_{\m}-t_{\m})
\left\{ -{\nu_0+1\over 4}
+{3\over 4}
\kio{\m}+c(t_{\m})\right\}+c(t_{\m-1})\cr
&~~~{\rm for}~~t_{\m}+1\leq j_{\m} \leq t_{\m+1}+2\delta_{\m,n}
~~{\rm and}~~1\leq \m\leq n.}}

Propositions 1 and 2 are very important because they enable us to prove
many identities without explicitly constructing all fermionic polynomials
in a one by one fashion. From them we will prove Theorem 1 (8.50)
and as a result obtain the 
Rogers-Schur-Ramanujan identities for $M(p,p')$ models at 
$b=l^{(\m)}_{1+j_{\m}}.$  
The reason for the introduction of the several sets of terms
$P_0$ and $P_{\pm 1}$ is that (as can be seen quite explicitly in
(10.5) with $j_0=0,1$) the
polynomials $F_{r(b+a),s}(L,b+a)$ for 
$b=l^{(\mu)}_{1+j_{\mu}},~a=1,2,~1\leq \mu \leq n,~1+t_{\mu}\leq j_{\mu}
\leq t_{\mu+1}-1+\delta_{\mu,n}$ can be written as 
\eqn\newthing{\eqalign{& \{F_{r(b+1),s}(L,b+1),F_{r(b+2),s}(L,b+2)\} \cr
&= q^{c(j_{\mu})-{\nu_0+1\over 4}+{3\over 4}\kio{\mu}}
\{f_s(L, {\bf e}_{1+j_{\mu}}-{\bf E}^{(t)}_{1,n}),
 f_s(L, {\bf e}_1+{\bf e}_{1+j_{\mu}}-{\bf E}^{(t)}_{1,n})\} \cr
&+q^{c(j_{\mu})}\{ P_0(L, \mu, {\bf e}_{j_{\mu}}-{\bf
E}^{(t)}_{\mu+1,n}),P_1(L, \mu, {\bf e}_{j_{\mu}}-{\bf
E}^{(t)}_{\mu+1,n})\}}}
and then we see from ~\pprop~and ~\ppmsflow~that the first pair 
on the right hand
side of~\newthing~ may be studied independently from the second pair under
the flow $\aw{b-2}, b\leq y_{\mu}$ (but we note this independence is
not true when $b >y_\mu).$ We note too that the
proofs given below of the two propositions are also quite 
independent. In particular, 
for~\pprop~ we use decomposition~\pppropa~ while for
~\ppmsflow~we use~\pppropb.

The final tool we need for our inductive proof are the following lemmas
which we use to treat the evolution along $\evs{2}$ in the
decompositions ~\pppropa~and~\pppropb:

{\bf Lemma 2.1}

{\it
Define $g(\m,n,j_{\m})=-{\nu_0+1\over 4}+{3\over 4}\kio{\m}
+{(-1)^{\m}\over 4}\theta(n>\m)\delta_{j_{\m},t_{\m+1}},~
~1\leq \m\leq n,~~1+\delta_{\m,1}+t_{\m}
\leq j_{\m} \leq t_{\m+1}+\delta_{\m,n}$,
and $f_s(L,{\bf e}_a+\e{2+t_{n+1}}-{\bf E}^{(t)}_{1,n}
+{\tilde p}{\bf e}_{1+t_{1+n}})=0$ 
for $a=0,1$ and ${\tilde p} \in Z.$ Then we have
\eqn\newlemmaa{\eqalign{
& \left\{
P_{-1}(L,\m-\delta_{1+t_{\m},j_{\m}},\e{j_{\m}}-
\theta(n>\mu){\bf e}_{t_{1+\mu}}+{\bf u}''(j_{\mu}))
\right. \cr
&~~~~~~~~ +q^{{L\over 2}-{\nu_0-\kie{\m}\over 4}} 
f_s(L,\e{-1+j_{\m}}-\theta(n>\mu){\bf e}_{t_{1+\mu}}
-{\bf E}^{(t)}_{1,\mu}+{\bf u}''(j_{\mu})),\cr
&\left.~~~~~~~~~~f_s(L,\e{j_{\m}}-\theta(n>\mu){\bf e}_{t_{1+\mu}}+
{\bf u}''(j_{\mu})) \right\}\cr
\evs{2} &  q^{g(\m,n,j_{\m})}\left\{ 
f_s(L,\e{1+j_{\m}}+\theta(n>\mu)\delta_{t_{1+\m},j_{\m}}
\e{t_{\m+1}}-\E{1}{\mu+\theta(n>\mu)}+{\bf u}''(j_{\mu})),\right. \cr
&\left.  
f_s(L,\e{1}+\e{1+j_{\m}}+\delta_{t_{1+\m},j_{\m}}
\theta(n>\m)\e{t_{\m+1}}
-\E{1}{\mu  +\theta(n>\mu)}+{\bf u}''(j_{\mu})) \right\} \cr
&~~~~+\left\{P_{0}(L,\m,\e{j_{\m}}-\theta(n>\mu){\bf e}_{t_{1+\mu}}
+{\bf u}''(j_{\mu})),\right. \cr
&~~~~~~~~~~~\left. P_{1}(L,\m,\e{j_{\m}}-\theta(n>\mu){\bf
e}_{t_{1+\mu}}+{\bf u}''(j_{\mu})) \right\}
}}

where
\eqn\udouble{{\bf u}''(j_{\mu})=\cases{{\bf
u}'_{\mu+1}-\delta_{j_{\mu},t_{1+\mu}}({\bf u}'_{\mu+1})_{1+t_{1+n}}{\bf
e}_{1+t_{1+\mu}}&if~$\mu<n$\cr
{\tilde p}{\bf e}_{1+t_{1+n}},~{\tilde p}\in Z&if~$\mu=n.$\cr}}
}

{\bf Lemma 2.2}

{\it
For $1+t_{\mu}\leq j_{\mu} \leq t_{1+\mu}-1-\delta_{\mu,1},~~1\leq \mu
\leq n-1$ we have
\eqn\newlemmab{\eqalign{&\{P_{-1}(L,\mu,{\bf
e}_{t_{1+\mu}-(j_{\mu}-t_{\mu})}-{\bf e}_{t_{1+{\mu}}}+{\bf
u}'_{\mu+1})+\cr
&~~~~~~q^{{L\over 2}-{\nu_0+1\over 4}+{3\over 4}\kio{\mu}}
f_s(L,{\bf e}_{t_{1+\mu}-(j_{\mu}-t_{\mu})+1}-{\bf E}^{(t)}_{1,\mu+1}+{\bf
u}'_{\mu+1}),\cr
&~~~~~~f_s(L,{\bf e}_{t_{1+\mu}-(j_{\mu}-t_{\mu})}-
{\bf e}_{t_{1+\mu}}+{\bf u}'_{\mu+1})\}\cr
&\evs{2}q^{-{\nu_0-\kie{\mu}\over 4}}\{f_s(L,{\bf
e}_{t_{1+\mu}-(j_{\mu}-t_{\mu})-1}-{\bf E}^{(t)}_{1,\mu+1}+{\bf
u}'_{\mu+1}),\cr
&~~~~~~~~~~~~~~f_s(L,{\bf e}_1+{\bf e}_{t_{1+\mu}-(j_{\mu}-t_{\mu})-1}-{\bf
E}^{(t)}_{1,\mu+1}+ {\bf u}'_{\mu+1})\}\cr
&+\{P_0(L,\mu-\delta_{j_{\mu},t_{1+\mu}-1},{\bf
e}_{t_{1+\mu}-(j_{\mu}-t_{\mu})}-{\bf e}_{t_{1+\mu}}+{\bf
u}'_{\mu+1}),\cr
&~~~~~~~~P_1(L,\mu-\delta_{j_{\mu},t_{1+\mu}-1},{\bf
e}_{t_{1+\mu}-(j_{\mu}-t_{\mu})}-{\bf e}_{t_{1+\mu}}+{\bf
u}'_{\mu+1})\}\cr}}
}
These lemmas follow immediately from the recursive properties 
~\freczero-\frecrele~of $f_s$ and the definition ~\ftildedf~of ${\tilde f}_s$.

{\bf Proof of Propositions 1 and 2}

We prove propositions 1 and 2 by induction on $\m$. For the proof of 
proposition 1 we use the decomposition of $\ev{y_{\m+1}-2}{a}$ as given
in \pppropa~and for the proof of proposition 2 the decomposition \pppropb.
Since in the decompositions \pppropa~and \pppropb~one uses  flows
$\ev{y_{\mu}-2}{\pm 1,0}$ and $\ev{y_{\mu-1}-2}{0}$
but never $\ev{y_{\mu-1}-2}{\pm 1}$
propositions 1 and 2 for $\m+1$ and all $a$ will follow if propositions 1 and 2
are true for $a=0$ and $\m, \m-1$ and $a=-1,1$ and $\m$.

One may check that it is sufficient to prove the following initial conditions
for propositions 1 and 2:

a) propositions 1 and 2 for $a=0$ and $\m=1,2$

b) propositions 1 and 2 for $a=-1$ and $\m=2$

c) proposition 1 for $a=-1$ and $\m=3$

With these initial conditions propositions 1 and 2 with $a=0$ follow for all
$\m\geq 1$, propositions 1 and 2 with $a=-1$ follow for
all $\m\geq 2$ and 
propositions 1 and 2 with $a=1$ follow for all $\m\geq 3$. 

{\it Proof of the initial conditions}
 
Here we prove point a) only. Points b) and c) are treated in appendix
B. 

We first consider proposition 1 for $\mu=1,2$ and $a=0$. 
Taking into account ~\prff~with
$j_0=0,1,\nu_0-1,\nu_0$ and case 2 of ~\tzcb~ with
$j_0=-1+t_1,~j_1=-1+t_2$ along with case 1 of ~\tzcb~ with $j_1=t_2$ and
recalling definition of $P_{-1}$~\ppmsl~ we see that
\eqn\initiala{\eqalign{&\{ f_s(L,-{\bf E}^{(t)}_{1,i}-{\bf
E}^{(t)}_{1+i,n}),  f_s(L,{\bf e}_1-{\bf E}^{(t)}_{1,i}-{\bf
E}^{(t)}_{1+i,n})\cr
&~~~~~~~~~\ev{y_i-2}{0}q^{c(t_{i})}\{P_{-1}(L,i,-{\bf
E}^{(t)}_{1+i,n}),f_s(L,-{\bf E}^{(t)}_{1+i,n})\}}}
with $i=1,2.$ Now if we replace $-{\bf E}^{(t)}_{1+i,n}$ by ${\bf
u}'_{i}$ in ~\initiala~we obtain proposition 1 for $\mu=1,2.$ Note
that the above replacement is legitimate because the recursive
properties of secs. 5 and 6 allow us to repeat the constructions of
sec. 7 for any vector ${\bf u}'_{i}.$

The proof of proposition 2 for $\mu=1,2$ and $a=0$ is only slightly more
involved. Making use of case 2 of ~\tzcb~with $j_0=0,1, \nu_0-1,$ 
case 1 of ~\tzcb~with $j_1\rightarrow j_1+1$ on one hand and
case 3 of ~\tzcc~with $j_0=0,1,~~1+t_2 \leq j_2 \leq t_3-2,$ case 6 of
~\tzcd~ with $j_0=-1+t_1$ and case 1 of ~\tzcc~ with $j_2\rightarrow j_2+1$ on
the other hand we obtain upon recalling the definitions
~\ppmsl-\ppmsh~and property ~\rpropiv~with $b=l^{(i)}_{1+j_i}<y_{n}$
\eqn\initialb{\eqalign{&q^{-{\nu_0+1\over 4}+{3\over 4}\delta_{1,i}}
\left\{f_s(L,-{\bf E}^{(t)}_{1,i}+{\bf e}_{1+j_i}-{\bf E}^{(t)}_{1+i,n}),
f_s(L,{\bf e}_1-{\bf E}^{(t)}_{1,i}+{\bf e}_{1+j_i}-{\bf
E}^{(t)}_{1+i,n})\right\}\cr
&~~~~~~~~~~~~~~~~~~+\left\{P_0(L,i,{\bf e}_{j_i}-{\bf E}^{(t)}_{1+i,n}),
P_1(L,i,{\bf e}_{j_i}-{\bf E}^{(t)}_{1+i,n})\right\}\cr
&\ev{y_i-2}{0}q^{c(t_i)-{\nu_0+1\over 4}+{3\over 4}\delta_{i,1}}\left\{
P_{-1}(L,i,{\bf e}_{1+j_i}-{\bf E}^{(t)}_{1+i,n}), f_s(L,{\bf
e}_{1+j_i}-{\bf  E}_{1+i,n}^{(t)})\right\}\cr
&~~~~~~~~~~~~~~~~~+q^{c(t_i)-{1\over 2}\delta_{i,2}}\left\{{\tilde f}_s(L,{\bf
e}_{j_i}-{\bf E}^{(t)}_{1,n}),0\right\}\cr}}
with $i=1,2.$ Next we use proposition 1 with $\mu=i,~{\bf u}'_{i}={\bf
e}_{1+j_i}-{\bf E}^{(t)}_{1+i,n}$ to derive from ~\initialb
\eqn\initialc{\left\{P_0(L,i,{\bf e}_{j_i}-{\bf E}^{(t)}_{1+i,n}),
P_1(L,i,{\bf e}_{j_i}-{\bf E}_{1+i,n}^{(t)})\right\}\ev{y_i-2}{0}
q^{c(t_i)-{1\over 2}\delta_{i,2}}\left\{{\tilde f}_s(L,{\bf
e}_{j_i}-{\bf E}^{(t)}_{1,n}),0\right\}.}
Thus, replacing ${\bf e}_{j_i}-{\bf E}^{(t)}_{1+i,n}$ by ${\bf u}'_i$
in ~\initialc~ we obtain proposition 2 with $\mu=1,2.$ 

We conclude this subsection with the following comment. In deriving
~\initialc~from ~\initialb~we succeeded in separating the evolution of
a $\{P_0,P_1\}$ pair from that of a $\{f_s,f_s\}$ pair. This
separation can be made in the formulas ~\tzcc~and \tzcd~as long 
as $l^{(2)}_{1+j_2}+1\leq b\leq l^{(2)}_{2+j_2},~1+t_2\leq j_2 \leq
t_3-1.$ All descendents of the $\{f_s,f_s\}$ pair will have ${\bf e}_{1+j_2}$
in their arguments and all descendants of the pair $\{P_0,P_1\}$ will have
${\bf e}_{j_2}$ in their arguments instead. Using this identification
principle we easily obtain for $\nu_1 \geq 3,~\nu_0 \geq 2$
\eqn\initiald{\eqalign{&\{P_0(L,2,{\bf u}'_2), P_1(L,2,{\bf u}'_2)\}\cr
&\aw{y_1+1}q^{-{\nu_0-1\over 2}}\{f_s(L,{\bf e}_1+{\bf
e}_{-2+t_2}-{\bf E}^{(t)}_{1,2}+{\bf u}'_2),
f_s(L,{\bf e}_2+{\bf e}_{-2+t_2}-{\bf E}^{(t)}_{1,2}+{\bf u}'_2)\}\cr
&+q^{-{\nu_0-2\over 4}}\{{\tilde f}_s(L,{\bf e}_{-2+t_1}+{\bf
e}_{-1+t_2}-{\bf E}_{1,2}^{(t)}+{\bf u}'_{2}),
{\tilde f}_s(L,{\bf e}_{-3+t_1}+{\bf
e}_{-1+t_2}-{\bf E}_{1,2}^{(t)}+{\bf u}'_{2})\}\cr}}
by comparing case 2 of ~\tzcc~ with $j_0=0,1;~j_1=1+t_1$ and case 3 of
~\tzcc with $j_0=0,1$ and replacing ${\bf  e}_{j_2}-{\bf
E}^{(t)}_{3,n}$ by ${\bf u}'_2.$ This result will be used in appendix B
 to prove proposition 2 with $a=-1,~\mu=2.$ 

{\it Proof of proposition 1 for $ \m-1,\m\rightarrow \m+1$}

Let us now show that proposition 1 with $a=-1,0,1$ follows inductively for 
$\m+1$. We will decompose $\ev{y_{\m+1}-2}{a}$ according to \pppropa~which
allows us to use proposition 1 and 2 for $\m$ and $\m-1$
for which they are true
by assumption.
Let us start by applying proposition 1 with $\m$,
$a'=-\delta_{a,1}$ and $\u{\m}
=\u{\m+1}-\e{t_{\m+1}}$ to
$\left\{ f_s(L,-\E{1}{\m+1}+\u{\m+1}),
f_s(L,\e{1}-\E{1}{\m+1}+\u{\m+1}) \right\}$
and subsequently using lemma 2.1 with $\mu\rightarrow
\mu-1,~~j_{\mu-1}=t_{\mu}$
 to obtain
\eqn\ppa{\eqalign{&
\left\{ f_s(L,-\E{1}{\m+1}+\u{\m+1}),
f_s(L,\e{1}-\E{1}{\m+1}+\u{\m+1}) \right\} \cr
&\ev{y_{\m}-2}{-\delta_{a,1}} q^{c(t_{\m})
+{1\over 2}\delta_{a,1}(-1)^{\m+1}}  \left\{
P_{-1}(L,\m,-\e{t_{\m+1}}+\u{\m+1}),
f_s(L,-\e{t_{\m+1}}+\u{\m+1}) \right\}\cr
&\evs{2} q^{c(t_{\m})+{1\over 2}\delta_{a,1}(-1)^{\m+1}} \left( 
q^{-{\nu_0-\kie{\m}\over 4}}\left\{
f_s(L,-\E{1}{\m-1}+\e{1+t_{\m}}-\e{t_{\m+1}}+\u{\m+1}),\right. \right.\cr
&~~~~~~~~~~~~~~~~~~~~~~~~~~~~~~~~~~\left. f_s(L,\e{1}-\E{1}{\m-1}+\e{1+t_{\m}}-
\e{t_{\m+1}}+\u{\m+1}) \right\}\cr
&+\left.\left\{P_0(L,\m-1,-\e{t_{\m+1}}+\u{\m+1}),
P_1(L,\m-1,-\e{t_{\m+1}}+\u{\m+1})\right\}\right).
}}
where we have noticed  from~\ppmsl~that
\eqn\anotherthing{\eqalign{P_{-1}(L,\mu,{\bf u}'_{\mu})&=P_{-1}(L,\mu-1,{\bf
e}_{t_\mu}-{\bf e}_{t_{\mu}}+{\bf u}'_{\mu})\cr
&+q^{{L\over 2}-{\nu_0-\kie{\mu-1}\over 4}}f_s(L,{
\bf e}_{-1+t_{\mu}}-{\bf e}_{t_{\mu}}-\E{1}{\mu-1}+{\bf u}'_{\mu}).\cr}}
We point out that the appearance of $\mu-1$ instead of $\mu$ in the
right hand side of ~\ppa~after $y_{\mu}$ steps explains the use of
decomposition ~\pppropa~instead of~\pppropb~ for the proof of proposition 1. 

We now apply proposition 1 with $\mu\rightarrow \mu-1,~a'=0$ to the
first pair and proposition 2 with $\mu \rightarrow \mu-1~a'=0$ to 
the second pair in the right hand side of ~\ppa~ 
which yields
\eqn\ppb{\eqalign{
&\left\{ f_s(L,-\E{1}{\m+1}+\u{\m+1}),
f_s(L,\e{1}-\E{1}{\m+1}+\u{\m+1}) \right\} \cr
\ev{y_{\m-1}+y_{\m}-2}{\pppropa} &  q^{c(1+t_{\m})
+{1\over 2}\delta_{a,1}(-1)^{\m +1}}
\left\{P_{-1}(L,\m-1,\e{1+t_{\m}}-\e{t_{\m+1}}+\u{\m+1}) \right.\cr
&\left.+q^{{L\over 2}-{\nu_0-\kie{\mu}\over 4}}
f_s(L,-\E{1}{\m-1}-\e{t_{\m+1}}+\u{\m+1}), 
f_s(L,\e{1+t_{\m}}-\e{t_{\m+1}}+\u{\m+1}) \right\}.
}}
Here we have used $c(1+t_{\m})=-{\nu_0-\kie{\m}\over
4}+c(t_{\m})+c(t_{\m-1})$
which follows from \pcj. The symbol \pppropa~under the arrow in ~\ppb~
means  that we have evolved
the initial state according to the first $y_{\m-1}+y_{\m}-2$ 
steps of the decomposition given by \pppropa.

Next applying lemma 2.1 with $\mu$ and $j_{\mu}=1+t_{\mu}$ we obtain
\eqn\ppc{\eqalign{
&\left\{ f_s(L,-\E{1}{\m+1}+\u{\m+1}),
f_s(L,\e{1}-\E{1}{\m+1}+\u{\m+1}) \right\} \cr
&\ev{y_{\m-1}+y_{\m}}{\pppropa} q^{c(1+t_{\m})
+{1\over 2}\delta_{a,1}(-1)^{\m+1}}\left(
q^{-{\nu_0+1\over 4}+{3\over 4}\kio{\m}}\left\{
f_s(L,-\E{1}{\m}+\e{2+t_{\m}}-\e{t_{\m+1}}+\u{\m+1}), \right. \right.\cr
&~~~~~~~~~~~~~~~~~~~~~~~~~~~~~~~~~~~~~~~~~~~~~~~~~~
\left. f_s(L,\e{1}-\E{1}{\m}+\e{2+t_{\m}}-\e{t_{\m+1}}+\u{\m+1})\right\}\cr
&+\left.\left\{P_0(L,\m,\e{1+t_{\m}}-\e{t_{\m+1}}+\u{\m+1}),
P_1(L,\m,\e{1+t_{\m}}-\e{t_{\m+1}}+\u{\m+1}) \right\}\right).
}}
Now we may use propositions 1 and 2 again with $\m$ and $a'=\delta_{a,-1}$ 
which yields
\eqn\ppd{\eqalign{
&\left\{ f_s(L,-\E{1}{\m+1}+\u{\m+1}),
f_s(L,\e{1}-\E{1}{\m+1}+\u{\m+1}) \right\} \cr
&\ev{y_{\m-1}+2y_{\m}-2}{\pppropa} 
q^{c(2+t_{\m})+{1\over 2}a(-1)^{\m+1}} 
\left\{P_{-1}(L,\m,\e{2+t_{\m}}
-\e{t_{\m+1}}+\u{\m+1}) \right. \cr
&\left. + q^{{L\over 2}-{\nu_0-\kie{\m}\over 4}}
f_s(L,-\E{1}{\m}+\e{1+t_{\m}}-\e{t_{\m+1}}+\u{\m+1}),
f_s(L,\e{2+t_{\m}}-\e{t_{\m+1}}+\u{\m+1}) \right\}.
}}
where we used 
$c(2+t_{\m})=c(1+t_{\m})-{\nu_0+1\over 4}+{3\over 4}\kio{\m}
+c(t_{\m})$ 
which again follows from \pcj.
Notice that the final entry of \ppb~$q^{c(1+t_{\m})
+{1\over 2}(-1)^{\m+1}\delta_{a,1}}
f_s(L,\e{1+t_{\m}}-\e{t_{\m+1}}+\u{\m+1})$ 
and the final entry of \ppd~$q^{c(2+t_{\m})+{1\over 2}a(-1)^{\m+1}} 
f_s(L,\e{2+t_{\m}}-\e{t_{\m+1}}+\u{\m+1})$ only differ 
in that  $\e{1+t_{\m}}$ has become $\e{2+t_{\m}}$ 
and the phase factor has changed. 
Applying now repeatedly $\evs{2}$ and $\ev{y_{\m}-2}{0}$
according to the 
decomposition \pppropa~and using lemma 2.1 and
propositions 1 and 2 for $\m$ we obtain
\eqn\ppe{\eqalign{
&\left\{ f_s(L,-\E{1}{\m+1}+\u{\m+1}),
f_s(L,\e{1}-\E{1}{\m+1}+\u{\m+1}) \right\} \cr
&\ev{y_{\m-1}+(j_{\m}-t_{\m})y_{\m}-2}{\pppropa} q^{c(j_{\m})+
{1\over 2}a(-1)^{\m+1}} 
\left\{ P_{-1}(L,\m,\e{j_{\m}}-\e{t_{\m+1}}+\u{\m+1}) \right. \cr
&\;\left.+q^{{L\over 2}-{\nu_0-\kie{\m}\over 2}}
f_s(L,-\E{1}{\m}+\e{-1+j_{\m}}-\e{t_{\m+1}}+\u{\m+1}),
 f_s(L,\e{j_{\m}}-\e{t_{\m+1}}+\u{\m+1})\right\}.
}}
The arrow in \ppe~denotes the flow
after the first $y_{\m-1}+(j_{\m}-t_{\m})y_{\m}-2$ 
steps ($t_{\m}+1<j_{\m}\leq t_{\m+1}$) according
to the decomposition \pppropa.
Finally setting $j_{\m}=t_{\m+1}$ in ~\ppe~and using the easily
verifiable identity 
\eqn\veriden{c(-1+t_{1+\mu})+{L-\nu_0+\kio{\mu}\over
2}+c(t_{\mu})-c(t_{1+{\mu}})={L\over 2}-{\nu_0-\kie{\mu}\over 4}}
along with ~\anotherthing~with $\mu\rightarrow \mu+1$ we obtain
proposition 1 for $\mu+1.$
This concludes the proof of
proposition 1.

{\it Proof of proposition 2 for $\m-1, \m\rightarrow \m+1$}

Let us show that proposition 2 holds for $\m+1$ and $a=-1,0,1$
inductively. First we assume $\nu_{\mu}>2.$
Recalling the definition of ${\bf u}'_{\mu}$ 
we see from \ppms~and  \ppmsh~that
\eqn\pa{\eqalign{&\left\{P_0(L,\m+1,\u{\m+1}),P_1(L,\m+1,\u{\m+1})\right\}\cr
&=q^{-{\nu_0-\kie{\m}\over 4}}
\left\{f_s(L,-\E{1}{\m+1}
+\e{-1+t_{\m+1}}+\u{\m+1}),\right.\cr
&\left.f_s(L,\e{1}-\E{1}{\m+1}
+\e{-1+t_{\m+1}}+\u{\m+1})\right\}
+\left\{P_0(L,\m,\u{\m+1}),P_1(L,\m,\u{\m+1})\right\}.\cr}}
We evolve these polynomials according to the decomposition \pppropb.
Thus we first evolve the first pair on the rhs of \pa~using proposition 1
and the second pair using proposition 2 for $\m$, $a'=-\delta_{a,1}$
and $\u{\m}= \e{-1+t_{\m+1}}-\e{t_{\m+1}}+\u{\m+1}$ followed by lemma
2.2 with $j_{\mu}=1+t_{\mu}$
to obtain
\eqn\pc{\eqalign{& \left\{ P_0(L,\m+1,\u{\m+1}),
P_1(L,\m+1,\u{\m+1}) \right\}\cr
\ev{y_{\m}-2}{-\delta_{a,1}} & 
q^{c(t_{\m})-{\nu_0-\kie{\m}\over 4}+{1\over 2}\delta_{a,1}(-1)^{\m+1}}
\left\{ P_{-1}(L,\m,\e{-1+t_{\m+1}}
-\e{t_{\m+1}}+\u{\m+1}) \right.  \cr
&~~~~~~\left.+q^{{L\over 2}-{\nu_0+1\over 4}+{3\over 4}\kio{\mu}}
f_s(L,-\E{1}{\m}+\u{\m+1}), 
f_s(L,\e{-1+t_{\m+1}}-\e{t_{\m+1}}+\u{\m+1}) \right\} \cr
\evs{2}& q^{c(t_{\m})-{\nu_0-\kie{\m}\over 4}
+{1\over 2}\delta_{a,1}(-1)^{\m+1}}\cr
&\times \left(
q^{-{\nu_0-\kie{\m}\over 4}}\left\{f_s(L,-\E{1}{\m}+\e{-2+t_{\m+1}}
-\e{t_{\m+1}}+\u{\m+1}),\right. \right.\cr
&~~~~~~~~~~~~~~~~~~~~~~~~~~~~~\left.f_s(L,\e{1}-\E{1}{\m}
+\e{-2+t_{\m+1}}-\e{t_{\m+1}}+\u{\m+1})\right\}\cr
&\left.+\left\{
P_0(L,\m,\e{-1+t_{\m+1}}-\e{t_{\m+1}}+\u{\m+1}), 
P_1(L,\m,\e{-1+t_{\m+1}}-\e{t_{\m+1}}+\u{\m+1}) \right\}\right).
}}
In the next step we apply proposition 1 and proposition 2 with $\m$ and 
$a'=-\delta_{a,0}-\delta_{a,1}$ which yields
\eqn\pd{\eqalign{
& \left\{ P_0(L,\m+1,\u{\m+1}),P_1(L,\m+1,\u{\m+1}) \right\}\cr
\ev{2y_{\m}-2}{\pppropb} &
q^{c(2+t_\mu)-c(t_{\m-1})+{1\over 2}a(-1)^{\m+1}}
\left\{ P_{-1}(L,\m,\e{-2+t_{\m+1}}-\e{t_{\m+1}}+\u{\m+1}) \right. \cr
&+q^{{L\over 2}-{\nu_0+1\over 4}+{3\over 4}\kio{\mu}}
f_s(L,-\E{1}{\m}+\e{-1+t_{\m+1}}-\e{t_{\m+1}}+\u{\m+1}),\cr
&~~~~~~~~~~~~~~~~\left. f_s(L,\e{-2+t_{\m+1}}-\e{t_{\m+1}}+\u{\m+1}) \right\}
}}
where  $\ev{2y_{\m}-2}{\pppropb}$ denotes the flow according to the first
$2y_{\m}-2$ steps of the decomposition \pppropb~ and we used the
identity $2c(t_\mu)-{\nu_0-\kie{\mu}\over
2}+{\delta{a,0}+2\delta_{a,1}\over
2}(-1)^{\mu+1}=c(2+t_{\mu})-c(t_{\mu-1})+{1\over 2}a(-1)^{\mu+1}$
which follows from ~\pcj.

Applying now further $j_{\m}-t_{\m}-2$ $(t_{\m}+3\leq 
j_{\m}\leq t_{\m+1}-1)$ 
times the combination $\evs{2}\ev{y_{\m}-2}{-1}$
using lemma 2.2 
and proposition 1 and 2 
with $\m$ and $a'=-1$ we obtain
\eqn\pda{\eqalign{
& \left\{ P_0(L,\m+1,\u{\m+1}),P_1(L,\m+1,\u{\m+1}) \right\}\cr
\ev{(j_{\m}-t_{\m}) y_{\m}-2}{\pppropb} &
q^{c(j_{\m})-c(t_{\m-1})+{1\over 2}a(-1)^{\m+1}} 
\left\{ P_{-1}(L,\m,\e{t_{\m+1}-(j_{\m}-t_{\m})}-\e{t_{\m+1}}
+\u{\m+1}) \right.\cr
&~+q^{{L\over 2}-{\nu_0+1\over 4}+{3\over 4}\kio{\mu}} f_s(L,-\E{1}{\m}
+\e{t_{\m+1}-(j_{\m}-t_{\m})+1}-\e{t_{\m+1}}+\u{\m+1}),\cr
&~\left. f_s(L,\e{t_{\m+1}-(j_{\m}-t_{\m})}-\e{t_{\m+1}}+\u{\m+1}) \right\}.
}}
Setting $j_{\m}=t_{\m+1}-1$ in the last formula gives
\eqn\pe{\eqalign{
& \left\{ P_0(L,\m+1,\u{\m+1}),P_1(L,\m+1,\u{\m+1}) \right\}\cr
\ev{(\nu_{\m}-1)y_{\m}-2}{\pppropb} & 
q^{v}\left\{P_{-1}(L,\m,\e{1+t_{\m}}-\e{t_{\m+1}}+\u{\m+1})\right.\cr
&\;\;\;+q^{{L\over 2}-{\nu_0+1\over 4}+{3\over 4}\kio{\mu}}
f_s(L,-\E{1}{\m}+\e{2+t_{\m}}-\e{t_{\m+1}}+\u{\m+1})
,\cr
&\;\;\;\left.
 f_s(L,\e{1+t_{\m}}-\e{t_{\m+1}}+\u{\m+1}) \right\}
}}
where $v=c(-1+t_{\m+1})-c(t_{\m-1})+{1\over 2}a(-1)^{\m+1}$.

Doing the next step we enter zone $\m-1$. Using lemma 2.2 with
$j_{\mu}=t_{\mu+1}-1$
we get
\eqn\pf{\eqalign{
& \;\;\;\left\{ P_0(L,\m+1,\u{\m+1}),P_1(L,\m+1,\u{\m+1}) \right\}\cr
&\ev{(\nu_{\m}-1)y_{\m}}{\pppropb}q^{v-{\nu_0-\kie{\m}\over 4}} \left\{
f_s(L,-\E{1}{\m-1}-\e{t_{\m+1}}+\u{\m+1}),\right.\cr
&~~~~~~~~~~~~~~~~~~~~~~~~~\left.f_s(L,\e{1}-\E{1}{\m-1}
-\e{t_{\m+1}}+\u{\m+1})\right\}\cr
&~~~~~~~~+q^v\left\{P_0(L,\m-1,\e{1+t_{\m}}-\e{t_{\m+1}}+\u{\m+1}),
P_1(L,\m-1,\e{1+t_{\m}}-\e{t_{\m+1}}+\u{\m+1})\right\}
}}
The appearance of $\mu-1$ instead of $\mu$ on the right hand side of
~\pf~which occurred after $(\nu_{\mu}-1)y_{\mu}$ steps explains the
choice of decomposition ~\pppropb~instead of ~\pppropa~in the proof of
proposition 2. Now using proposition 1 and proposition 2 with
$\mu\rightarrow \mu-1$ and $a'=0$ we find
\eqn\pg{\eqalign{
& \;\;\;\left\{ P_0(L,\m+1,\u{\m+1}),P_1(L,\m+1,\u{\m+1}) \right\}\cr
&\ev{(\nu_{\m}-1)y_{\m}+y_{\m-1}-2}{\pppropb}
q^{v+c(t_{\m-1})-{\nu_0-\kie{\m}\over 4}}
\left\{P_{-1}(L,\m-1,-\e{t_{\m+1}}+\u{\m+1}) \right. \cr
&\;\;\;\;\;\left.+q^{{L\over 2}-{\nu_0-\kie{\mu}\over 4}}
f_s(L,-\E{1}{\m-1}+\e{1+t_{\m}}-\e{t_{\m+1}}+\u{\m+1}), 
f_s(L,-\e{t_{\m+1}}+\u{\m+1})\right\}
}}

To complete the proof, we evolve \pg~ according to $\evs{2}$ and
use the fermionic recursive property \frecreld, the first equation of
~\freczero~and definition ~\ftildedf~ to obtain
\eqn\ph{\eqalign{
& \left\{ P_0(L,\m+1,\u{\m+1}),P_1(L,\m+1,\u{\m+1}) \right\}\cr
\ev{(\nu_{\m}-1)y_{\m}+y_{\m-1}}{\pppropb} &
q^{v-{\nu_0-\kie{\m}\over 4}+c(t_{\m-1})} \left\{
P_0(L,\m,-\e{t_{\m+1}}+\u{\m+1}),
P_1(L,\m,-\e{t_{\m+1}}+\u{\m+1})\right\}
}}
and then apply proposition 2 with $\m$ and $a'=0$ to arrive at
\eqn\pia{\eqalign{
& \left\{ P_0(L,\m+1,\u{\m+1}),P_1(L,\m+1,\u{\m+1}) \right\}\cr
\ev{y_{\m+1}-2}{\pppropb} & \left\{
q^{v-{\nu_0-\kie{\m}\over 4}+c(t_{\m-1})+
{L-\nu_0+\kio{\m}\over 2}+c(t_{\m})}
f_s(L,-\E{1}{\m+1}+\u{\m+1}),0\right\}
}}
where we used 
$(\nu_{\m}-1)y_{\m}+y_{\m-1}+y_{\m}-2=y_{\m+1}-2$. 
Finally from ~\pcj~we derive
\eqn\pj{\eqalign{
&v-{\nu_0-\kie{\m}\over 4}+c(t_{\m-1})
+{L-\nu_0+\kio{\m}\over 2}+c(t_{\m})\cr
=&{L-\nu_0+\kio{\m+1}\over 2}+c(t_{\m+1})+{1\over 2}a(-1)^{\m+1}.
}}
Combining~\pia~and ~\pj~we derive proposition 2 for $\mu+1$ 
This concludes the proof of proposition 2 for $\nu_{\mu}>2.$

When $\nu_{\mu}=2$ a separate treatment is needed because the
decomposition ~\pppropb~does not hold. However, in this case the
decomposition ~\pppropa~ with $\nu_{\mu}=2$ can be used and with this the
proof of proposition 2 for $\nu_{\mu}=2$ follows.

Recalling lemma 2.1, the equations ~\ppb,~\ppe~and the definition of
Takahashi length~\stringdef~
we organize the results proven in this section as 

{\bf Theorem 1}

{\it For $1\leq \m\leq n,~~1+\delta_{\m,1}+t_{\m}
\leq j_{\m} \leq t_{\m+1}+\delta_{\m,n}$, we have 
\eqn\ptheo{\eqalign{
& \left\{ f_s(L,-\E{1}{n}+P\e{1+t_{1+n}}),
f_s(L,\e{1}-\E{1}{n}+P\e{1+t_{1+n}}) \right\} \cr
&\aw{l^{(\m)}_{1+j_{\m}}-2}  q^{c(j_{\m})}
\left\{
P_{-1}(L,\m-\delta_{1+t_{\m},j_{\m}},\e{j_{\m}}-\E{\m+1}{n}+P\e{1+t_{1+n}})
\right. \cr
&\left. +q^{{L\over 2}-{\nu_0-\kie{\m}\over 4}} 
f_s(L,\e{-1+j_{\m}}-\E{1}{n}+P\e{1+t_{1+n}}),
f_s(L,\e{j_{\m}}-\E{\m+1}{n}+P\e{1+t_{1+n}}) \right\}\cr
&\evs{2}  q^{c(j_{\m})}\left( q^{g(\m,n,j_{\m})} 
\left\{f_s(L,\e{1+j_{\m}}+\delta_{t_{1+\m},j_{\m}}
\theta(n>\m)\e{t_{\m+1}}-\E{1}{n}+P\e{1+t_{1+n}}),\right. \right. \cr
&~~~~~~~~~~~~~~~~~~~~~~~\left.f_s(L,\e{1}+\e{1+j_{\m}}+\delta_{t_{1+\m},j_{\m}}
\theta(n>\m)\e{t_{\m+1}}
-\E{1}{n}+P\e{1+t_{1+n}})\right\}\cr
&~~~~\left.\left\{P_{0}(L,\m,\e{j_{\m}}-\E{\m+1}{n}+P\e{1+t_{1+n}}), 
P_{1}(L,\m,\e{j_{\m}}-\E{\m+1}{n}+P\e{1+t_{1+n}}) \right\}\right)\cr
}}
}
where $\E{n+1}{n}=0,~P=0,1$
and $g(\mu,n,j_{\mu})$ is defined just above ~\newlemmaa. 

Before we move on, we observe that the $r(b)$ graph of sec. 3
for $y_{\m}-y_1-1 \leq b \leq y_{\m}$ is (up to a shift) for $\mu \geq
2$
\eqn\ppgraph{\psfig{file=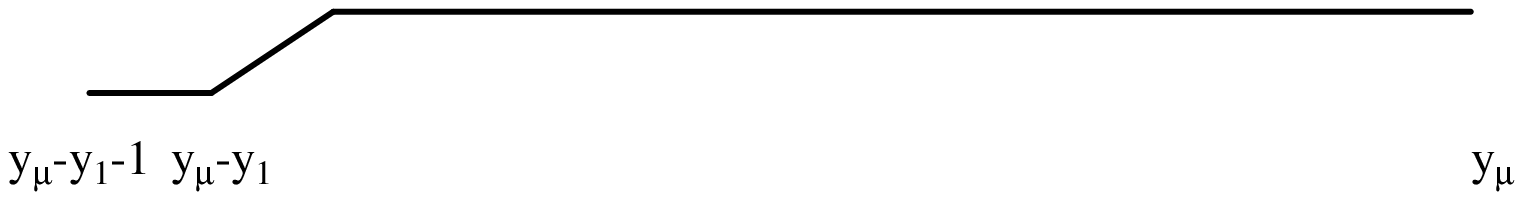,width=2.5in}}
Using the fermionic recursive properties of  sec. 5 for $f_s$ and~\reccft~ 
for ${\tilde f}_s$ as well as definition~\ftildedf~one can easily show that
\eqn\pfina{\eqalign{
\left\{ f_s(L,\e{1+t_1}-\E{2}{\m}+\u{\m}),
f_s(L,-\E{2}{\m}+\u{\m}) \right\}
\ev{y_1}{\ppgraph} \left\{ \tilde{f}_s(L,-\E{1}{\m}+\u{\m}),0 \right\}
}}
Notice that if we use
\eqn\ppgrapha{\psfig{file=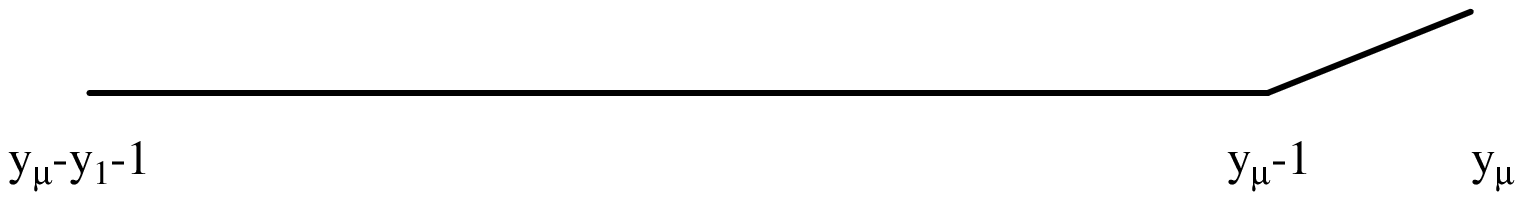,width=2.5in}}
instead, we obtain
\eqn\pfinb{\eqalign{
&\left\{ f_s(L,\e{1+t_1}-\E{2}{\m}+\u{\m}),f_s(L,-\E{2}{\m}+\u{\m}) \right\}\cr
\ev{y_1-1}{\ppgrapha} & \left\{ f_s(L,\e{1}-\E{1}{\m}+\u{\m}),
f_s(L,-\E{1}{\m}+\u{\m})\right\}\cr
\ev{1}{\pprbb~b} & \left\{ f_s(L,-\E{1}{\m}+\u{\m}),0 \right\}
}}
Comparing \pfina~with proposition 2 with $a=0$ we infer that
 for $\m\geq 2$
\eqn\pfinc{\eqalign{
& \left\{ P_{0}(L,\m,\u{\m}),P_{1}(L,\m,\u{\m}) \right\} \cr
\ev{y_{\m}-y_1-2}{~} & q^{c(t_{\m})-{1\over 2}\kie{\m}}
\left\{ f_s(L,\e{1+t_1}-\E{2}{\m}+\u{\m}),f_s(L,-\E{2}{\m}+\u{\m}) \right\}.
}}

We are now ready to discuss the evolution along the final
stretch of the
$r(b)$ map
\eqn\ppfin{\psfig{file=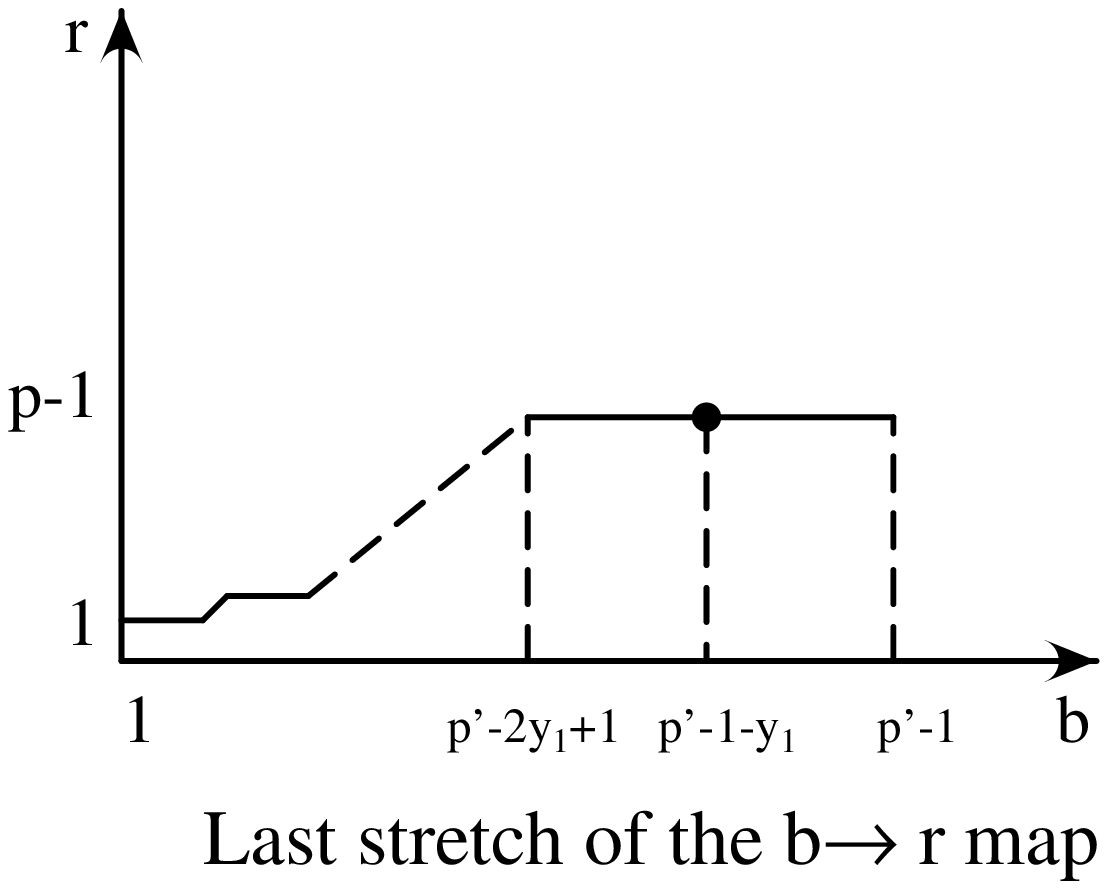,width=2.5in}}

Notice that the piece of the $b\to r$ map with
 $1+l^{(n)}_{2+t_{n+1}}\leq b
\leq p'-y_1-1$ is the same (up to a shift) as the map $b\to r$
of sec. 3
with $1\leq b\leq y_n-y_1-1$. The piece of the $b\to r$ map restricted to the interval $[p' -y_1 -1, p' -1]$ is graph \ppgrapha~with $\mu=n$ and the
last segment removed. Recalling theorem
1 (with $\m=n$ and $j_n=1+t_{n+1}$) we have
\eqn\pfind{\eqalign{
& \left\{ f_s(L,-\E{1}{n}+P\e{1+t_{1+n}}),
f_s(L,\e{1}-\E{1}{n}+P\e{1+t_{1+n}}) \right\} \cr
\aw{l^{(n)}_{2+t_{n+1}}=p'-y_n} & q^{c(1+t_{n+1})}
\left\{ P_{0}(L,n,(1-P)\e{1+t_{n+1}}),P_{1}(L,n,(1-P)
\e{1+t_{n+1}}) \right\}
}}
Applying \pfinc~and \pfinb~with $\m=n$, $\u{n}=(1-P)\e{1+t_{n+1}}$ 
to the right hand side of \pfind, we arrive at

{\bf Theorem 2}

\eqn\ptheof{\eqalign{
& \left\{ f_s(L,-\E{1}{n}+P\e{1+t_{1+n}}), 
f_s(L,\e{1}-\E{1}{n}+P\e{1+t_{1+n}}) \right\} \cr
\aw{p'-3} & q^{c(t_{n+1}+1)+c(t_n)-{1\over 2}\kie{n}}\cr
&~~~\times
\left\{ f_s(L,\e{1}-\E{1}{n}+(1-P)\e{1+t_{1+n}}),
f_s(L,-\E{1}{n}+(1-P)\e{1+t_{n+1}})\right\}
}}

If we now identify the fermionic polynomials generated in 
the evolution \ptheof~
with $F_{r(b),s}(L,b)$ (for $P=0$) and with $F_{r(b),s}^{(o)}(L,b)$ (for
$P=1$), then we have 
\eqn\Ffida{\eqalign{
&F_{1,s}(L,1)=f_s(L,-{\bf E}^{(t)}_{1,n})\cr
&F_{1,s}(L,2)=f_s(L,{\bf e}_1-{\bf E}^{(t)}_{1,n})\cr
&~~~~~~~~~~~\cdots\cr
&F_{p-1,s}(L,p'-2)=q^{c(1+t_{1+n})+c(t_n)-{1\over 2}\kie{n}}
f_s(L,\e{1}-\E{1}{n}+\e{1+t_{1+n}})\cr
&F_{p-1,s}(L,p'-1)=q^{c(1+t_{1+n})+c(t_n)-{1\over 2}\kie{n}}
f_s(L,-\E{1}{n}+\e{1+t_{1+n}})
}}
and
\eqn\Ffidb{\eqalign{
&F_{1,s}^{(o)}(L,1)=f_s(L,-\E{1}{n}+\e{1+t_{1+n}})\cr
&F_{1,s}^{(o)}(L,2)=f_s(L,\e{1}-\E{1}{n}+\e{1+t_{1+n}})\cr
&~~~~~~~~~~~\cdots\cr
&F_{p-1,s}^{(o)}(L,p'-2)=q^{c(1+t_{1+n})+c(t_n)-{1\over 2}\kie{n}}
f_s(L,\e{1}-\E{1}{n})\cr
&F_{p-1,s}^{(o)}(L,p'-1)=q^{c(1+t_{1+n})+c(t_n)-{1\over 2}\kie{n}}
f_s(L,-\E{1}{n})
}}

Finally we use the first equation in \freczero~to verify that
$F_{p-1,s}(L,p'-2), F_{p-1,s}(L,p'-1)$ ($F^{(o)}_{p-1,s}(L,p'-2), 
F^{(o)}_{p-1,s}(L,p'-1)$) satisfy the closing equation in
~\brecreltd. Thus we have  shown that
the constructive procedure defined in secs. 7 and 8 with the initial values
for $b=1$ and $b=2$ specified by the first two equations in \Ffida~or in
\Ffidb~gives rise to fermionic
polynomials which satisfy all bosonic recursion relations 
 ~\brecrelt--\brecreltd~ for $1\leq b \leq p'-1$

\newsec{Normalization constants and boundary conditions}

{}From the conclusion of the previous section it follows that when 
$s=l^{(\mu_s)}_{1+j_s},~1+t_{\mu_s}\leq j_s\leq
t_{1+\mu_s},~L+b+s\equiv 0 ({\rm mod} 2)$
the identity 
\eqn\bequ{F_{r(b),s}(L,b)=
\sum_{s'=1\atop s'\equiv s ({\rm mod} 2)}^{p'-1} k_{s,s'} 
{\tilde B}_{r(b),s'}(L,b)}
will hold for $L>0,$ provided constants $k_{s,s'}$ can be chosen so
that ~\bequ~holds for $L=0$. Using ~\mobose,~\binitial~it is trivial to verify that
for
\eqn\newnorm{k_{s,s'}=F_{r(s'),s}(0,s')}
with $s'\equiv s({\rm mod} 2),~1\leq s'\leq p'-1$ eqn.~\bequ~indeed holds
for $L=0.$ However,~\newnorm~is of very little use because
$F_{r(b),s}(L,b)$ have not been explicitly constructed for all $b\in
[1,p'-1].$ Fortunately, it turns out that the constants $k_{s,s'}$ can
be determined from ~\bequ~with $L=0,1,\cdots ,p'-1$ and $b=1,s,p'-1.$
In this direction we first calculate the threshold values of $L$,
i.e. the lowest values of $L$ such that $F_{r(b),s}(L,b){\not =} 0$ 
for $b=1,s,p'-1$
%

We know that $F_{r(b),s}(L,b)$ with $b=1,p'-1,s$
is given by 
\eqn\bFequ{F_{r(b),s}(L,b)=
q^{\cal N}\sum_{\vn \in \Z^{t_{1+n}} \atop m_{t_{1+n}}\equiv P (\mod 2)}
 q^{Q({\bf n},{\bf m})
+\vA^T \tilde{\vm}} \prod_{j=1}^{t_{n+1}}\nom{n_j+m_j}{n_j}^{(0)}.}
Here $\tilde{\vm}$ is defined in \mtil,
\eqn\newddefn{{\cal N}=\cases{0&if $b=1$\cr
c(j_s)& if $b=s$\cr
c(t_{n+1}+1)+c(t_n)-{1\over 2}\kie{n} & if $b=p'-1$}} 
and 
\eqn\bP{P=\cases{ 0 & if $b=1,s$\cr 1 & if $b=p'-1$.}}
${\bf n}$ and ${\bf m}$ are related by \mnsys~with
\eqn\buequ{\eqalign{
&{\bu}=-{\bf {\bar E}}^{(t)}_{1,n}+{\bf {\bar e}}_{j_s}-
{\bf \bar E}^{(t)}_{1+\mu_s,n}~~~{\rm for}~~b=1,p'-1\cr
&{\bu}=2{\bf \bar e}_{j_s}-2{\bf \bar E}^{(t)}_{1+\mu_s,n}
~~~{\rm for}~~b=s}}
where
\eqn\miwabar{{\bf \bar E}^{(t)}_{a,b}=\sum_{i=a}^{b}{\bf \bar e}_{t_i}.}

Even though $\vn, \vm \in \Z^{t_{1+n}}$, a bit of analysis shows that
effectively
\eqn\bmnrange{\eqalign{
&{\rm For}~~b=1,s:\cr
&n_i\geq \delta_{i,j_s}-\delta_{i,t_j},~~~t_j\geq j_s\cr
&m_i\geq 0\cr 
&{\rm For}~~b=p'-1:\cr
&n_i\geq 0,~m_i>0.
}}
Moreover, if $n_i$ takes on negative values then $m_i=0$.

Therefore, the threshold configurations for the three cases are 
\eqn\btconf{\eqalign{
&{\bf n}={\bf \bar e}_{j_{s}}-\theta(\m_s<n)
{\bf \bar E}^{(t)}_{1+\m_s,n},~m_{t_{n+1}}\equiv 0~({\rm mod}~2)~
{\rm for}~~b=1\cr
&{\bf n}=0,~m_{t_{n+1}}\equiv 1~({\rm mod}~2)~{\rm for}~~b=p'-1\cr
&{\bf n}={\bf \bar e}_{j_s}-\theta(\m_s<n){\bf \bar
E}^{(t)}_{1+\mu_s,n},~m_{t_{n+1}}\equiv
0~~({\rm mod}~2)~~{\rm for}~~b=s}}
with the corresponding thresholds from \partprob
\eqn\bthres{\eqalign{
L_{\rm tr}=&\sum_{i=1}^{\m_s}(y_i-y_{i-1})+l_{j_{s}}=s-1
~~~{\rm for}~~b=1\cr
L_{\rm tr}=&\sum_{i=1}^{n}(y_i-y_{i-1})+\sum_{i=\m_s+1}^n(y_i-y_{i-1})
-l_{j_{s}}+l_{1+t_{n+1}}\cr
=&p'-1-s~~~{\rm for}~~b=p'-1\cr
L_{\rm tr}=&0~~~{\rm for}~~b=s.
}}
with $l_j$ defined by ~\ldef~and $\m_s$ by \sdecom.
Finally the threshold for $B_{r(b),s}(L,b)$ is
\eqn\bthrsB{L_{\rm tr}=|s-b|.}

Hence for $b=1,~L=0,1,\ldots,s-2$ and $L+1+s\equiv 0 (\mod 2)$ 
we have from \bequ
\eqn\bequa{
0=\sum_{s'=1 \atop s'\equiv s (\mod 2)}^{p'-1}
k_{s,s'} {\tilde B}_{1,s'}(L,1).
}
Clearly, \bequa~and \bthrsB~imply that
\eqn\bequb{k_{s,s'}=0~~{\rm for}~~
1\leq s'\leq s-1,~s'\equiv s (\mod 2).}
Analogously, if we evaluate \bequ~at $b=p'-1,~
L=0,1,2,\ldots,p'-2-s$ and $L+p'-1+s\equiv 0 (\mod 2)$
we obtain
\eqn\bequc{
0=\sum_{s'=1 \atop s' \equiv s (\mod 2)}^{p'-1}
k_{s,s'} {\tilde B}_{p-1,s'}(L,p'-1).
}
\bthrsB~and \bequc~imply that
\eqn\bequd{
k_{s,s'}=0~~{\rm for}~~1+s\leq s',s'\equiv s (\mod 2).}
Combining \bequb~and \bequd~yields
\eqn\beque{
k_{s,s'}=\delta_{s',s}q^{a(j_{s})}.}
To determine $a(j_{s})$ we evaluate \bequ~at $b=s,~L=0$ to get 
\eqn\bequf{
F_{r(s),s}(L=0,b=s)=q^{a(j_{s})}
\tilde{B}_{r(s),s}(L=0,b=s)}
or since $F_{r(s),s}(L=0,b=s)=q^{c(j_{s})}$ and 
${\tilde B}_{r(s),s}(L=0,b=s)=1$ we have
\eqn\bequg{
a(j_{s})=c(j_{s}).}
This finally leads to
\eqn\bequh{
F_{r(b),s}(L,b)=q^{c(j_s)}\tilde{B}_{r(b),s}(L,b)}
where $L+b+s \equiv 0 (\mod 2).$

Analogously, we can show that
\eqn\bequj{
F_{r(b),s}(L,b)=0~~{\rm for}~~L+b+s\equiv 1 (\mod 2).}

This concludes the proof of the boundary conditions which leads to the
equality of the fermionic and the bosonic form of $M(p,p')$ polynomial
character identities.

\newsec{Rogers-Schur-Ramanujan type identities for $M(p,p')$ minimal models}

We now may give the Rogers-Schur-Ramanujan type identities for the the four
cases of ~\results. From ~\bequh~we see that all identities are of the
form
\eqn\finalid{F_{r(b),s}(L,b)=q^{{1\over
2}(\phi_{r(b),s}-\phi_{r(s),s})+c(j_s)}B_{r(b),s}(L,b).}
We thus state the results by giving the appropriate 
functions $F_{r(b),s}(L,b).$

{\bf 10.1. Rogers-Schur-Ramanujan identities for $b=l^{(\m)}_{1+j_{\m}}$ a
pure Takahashi lengths}

The Rogers-Schur-Ramanujan identities for $b=l^{(\m)}_{1+j_{\m}}$ a pure 
Takahashi
length follow immediately from theorem 1 of sec. 8 and the proof of the
boundary conditions of sec. 9. Since 
\eqn\tone{\eqalign{
&F_{1,s}(L,1)=f_s(L,-\E{1}{n})\cr
&F_{1,s}(L,2)=f_s(L,\e{1}-\E{1}{n})
}}
we infer from \ptheo~ that
\eqn\pRR{F_{r(b),s}(L,l^{(\m)}_{1+j_{\m}})=q^{c(j_{\m})}
f_s(L,\e{j_{\m}}-\E{\m+1}{n})}
where $c(j_{\m})$ was defined in \pcj~and $r(b)=\delta_{\m,0}
+\tilde{l}^{(\m)}_{1+j_{\m}}$. Thus we have explicitly proven the
result announced in~\rbm.  We would like to stress that (10.3) holds 
even in the cases where some or all $\nu _i =1$ or $\nu_n=0$ (or both)
with appropriate modification of the exponent $c(j{_\mu } )$.

{\bf 10.2. The vicinity of the Takahashi length}

Using the fermionic recursive properties of sec. 5 and 6 and the
Theorem 1 (of sec. 8), it is easy to 
extend the Rogers-Schur-Ramanujan identities of the previous subsection 
to case 2 of~\results.

{\bf Case 2a:} $l^{(\m)}_{1+j_{\m}}-\nu_0+1\leq b \leq  
l^{(\m)}_{1+j_{\m}}-1,~~\m\geq 1$

\eqn\va{\eqalign{
& F_{r(b),s}(L, l^{(\m)}_{1+j_{\m}}-j_0)=
q^{c(j_{\m})} \left( q^{{\nu_0-1-\kio{\m}\over 4}}
\tilde{f}_s(L,\e{j_0-1}-\E{1}{n}+\e{-1+j_{\m}})\right.\cr
&~~~~~+\theta(\m\geq 2)
\sum_{i=2}^{\m} q^{{\nu_0-1-\kie{i}\over 4}} \tilde{f}_s(L,
\e{j_0-1}-\E{1}{i}+\e{-1+t_i}+\e{j_{\m}}-\E{\m+1}{n})\cr
&~~~~~\left. +f_s(L,\e{\nu_0-j_0}-\e{t_1}+\e{j_{\m}}-\E{\m+1}{n})\right)~~
1\leq j_0 \leq \nu_0-1\cr}}

{\bf Case 2b:} $l^{(\m)}_{1+j_{\m}}+1\leq b\leq  
l^{(\m)}_{1+j_{\m}}+\nu_0+1,~~\m\geq 1$
\eqn\vaa{\eqalign{
&F_{r(b),s}(L,l^{(\m)}_{1+j_{\m}}+1+j_0)=
q^{c(j_{\m})} \left( q^{-{\nu_0+1\over 4}+{3\over 4}\kio{\m}}
f_s(L,\e{j_0}-\E{1}{n}+\e{1+j_{\m}})\right.\cr
&~~~~+\theta(\m\geq 2) \sum_{i=2}^{\m} q^{-{\nu_0-\kio{i}\over 4}}
f_s(L,\e{j_0}-\E{1}{i}+\e{-1+t_i}+\e{j_{\m}}-\E{\m+1}{n})\cr
&~~~~~\left. +\tilde{f}_s(L,\e{\nu_0-j_0-1}-\e{t_1}+\e{j_{\m}}
-\E{\m+1}{n}) \right),
~~0\leq j_0 \leq \nu_0\cr}}
where recalling ~\ftildeff~we note that case 2a agrees with ~\pRR~if $j_0=0$.
$r$ in this section is a function of $b$ as given in \results.

{\bf 10.3. Further cases}

The Propositions 1 and 2 of sec. 8 are very powerful and can be used
to to study all cases $1\leq b \leq p'-1$. As an illustration we present
the results for the fermionic forms of cases 3 and 4 of ~\results. The
details of the derivation are in appendix C.

To state the results we need several definitions (with $\alpha$ and
$\beta$ defined by the Takahashi decomposition of case 4 in ~\results)
 
{\bf Definitions}
\eqn\xxi{\eqalign{
&{\bf j}=|j_{\alpha},j_{\alpha+1},\ldots,j_{\beta}>\cr
&{\bf i}=|i_{\alpha},i_{\alpha+1},\ldots,i_{\beta-1},0>}}
with $i_k=0,1 $ for $ \alpha \leq k \leq \beta -1 $.
When $i_{\alpha}=0$ we define $a_i$ and $b_j$ from
\eqn\oxirepa{
{\bf i}=|\underbrace{0,\cdots,0}_{a_1},\underbrace{1,\cdots,1}_{b_2},
\underbrace{0,\cdots,0}_{a_2},\cdots ,\underbrace{1,\cdots
,1}_{b_l},
\underbrace{0,\cdots ,0}_{a_l}>}
where $1\leq a_i,~1\leq i\leq l$ and $1\leq b_j,~2\leq j\leq l$. If $l=1$
then ${\bf i}=|\underbrace{0, \cdots, 0}_{a_1}>$. Note that 
${\alpha}-1+a_1+\sum_{j=2}^{l}(a_j+b_j)=\beta$. 
Similarly when $i_{\alpha}=1$ we write
\eqn\oxirepb{{\bf i}=|\underbrace{1,\cdots, 1}_{b_1},\underbrace{0,\cdots
,0}_{a_1},\underbrace{1,\cdots , 1}_{b_2},\underbrace{0,\cdots ,
0}_{a_2},\cdots ,\underbrace{1,\cdots, 1}_{b_l},
\underbrace{0,\cdots ,0}_{a_l}>}
where $1\leq a_j,~b_j$ with $1\leq j\leq l$. From  $\bf i$ we further define
\eqn\xxR{R_{i_{\m}}(j_{\m}+1)=\cases{ j_{\m}+1 & if $i_{\m}=0$\cr
t_{\m+1}-(j_{\m}-t_{\m})-1 & for $i_{\m}=1.$}}

{\bf Results}

With these definitions we may have the following results for the
fermionic polynomials of cases 3 and 4 of ~\results

{\bf Case 3 of~\results}: 

Here $\alpha=0$ and we have
\eqn\xxRR{\eqalign{
F_{r(b),s}(L,b)&=q^{c_{(3)}({\bf j})}\left( 
\sum_{i_1,\ldots,i_{\beta-1}=0,1 \atop i_0=0}
q^{rf_{(3,1)}({\bf i})}
f_s(L,\e{j_0+i_1}-\e{t_1}+{\bf u}_{(3)}({\bf i},{\bf j}))
\right.\cr
&+\left. \sum_{i_1,\ldots,i_{\beta-1}=0,1 \atop i_0=1}
q^{rf_{(3,2)}({\bf i})}
\tilde{f}_s(L,\e{\nu_0-j_0-i_1-1}-\e{t_1}
+{\bf u}_{(3)}({\bf i},{\bf j}))\right)\cr}}
where
\eqn\xxc{c_{(3)}({\bf j})=\sum_{\m=1}^{\beta}
(c(j_{\m})-{\nu_0+1\over 4}+{3\over 4}\kio{\m})}
\eqn\xxu{{\bf u}_{(3)}({\bf i},{\bf j})=
\sum_{\m=1}^{\beta-1} \e{R_{i_{\m}}(j_{\m}+1)+|i_{\m+1}-i_{\m}|
-|i_{\m}-i_{\m-1}|}+\e{1+j_{\beta}-i_{-1+\beta}}-\E{2}{n},}
for $i_0=0$ we set
\eqn\xxrelfaca{
rf_{(3,1)}({\bf i})={1\over 4}\delta_{a_1,1}+{1\over 2}\sum_{j=2}^l 
(-1)^{a_1+\sum_{k=2}^{j-1}(a_k+b_k)}\kie{b_j}.}
and for $i_0=1$
\eqn\xxrelfacb{
rf_{(3,2)}({\bf i})={\nu_0+(-1)^{b_1}\over 4}-{\delta_{b_1,1}\over 4}+
{1\over 2}\sum_{j=2}^l (-1)^{\sum_{k=2}^{j-1}(a_k+b_k)}\kie{b_j}}
where we have used the convention that $\sum_{i=a}^{b}=0$ if $b<a$.

{\bf Case 4 of ~\results}: 

Here $\alpha \geq 1$ and we have
\eqn\oxRR{\eqalign{F_{r(b),s}(L,b)&=
 q^{c_{(4)}({\bf j})} \left( \sum_{i_{\alpha+1},\ldots,i_{\beta-1}=
0,1 \atop i_{\alpha}=0}
q^{rf_{(4,1)}({\bf i})} f_s(L,{\bf u}_{(4,1)}({\bf i},{\bf j}))\right.\cr
&~~~+\sum_{i_{\alpha+1},\ldots,i_{\beta-1}=0,1 \atop i_\alpha=1}
\left( q^{rf_{(4,2)}({\bf i})+\delta_{1,\alpha}({L\over 2}-{\nu_0\over 4})} 
f_s(L,{\bf u}_{(4,2)}({\bf i},{\bf j}))\right. \cr
&~~~~~~~~~~~~~~\left. 
\left. +q^{rf_{(4,2)}({\bf i})-{(-1)^{\alpha}+\delta_{1,\alpha}\over 4}} 
f_s(L,{\bf u}_{(4,3)}({\bf i},{\bf j})) \right) \right) \cr}}
where
\eqn\oxc{
c_{(4)}({\bf j})=c(j_\alpha)+\sum_{\m=\alpha+1}^{\beta}(c(j_{\m})
-{\nu_0+1\over 4}+{3\over 4}\kio{\m})}
\eqn\oxua{\eqalign{
{\bf u}_{(4,1)}({\bf i},{\bf j})&
=\e{j_{\alpha}+i_{\alpha+1}}+{\bf u}_{(4)}({\bf i},
{\bf j})\cr
{\bf u}_{(4,2)}({\bf i},{\bf j})&=
-\e{t_{\alpha}}+\e{t_{\alpha+1}-(j_{\alpha}-t_\alpha)
-1-i_{\alpha+1}}+{\bf u}_{(4)}({\bf i},{\bf j})\cr
{\bf u}_{(4,3)}({\bf i},{\bf j})&=\e{t_{\alpha}-1}
-\e{t_{\alpha}}+\e{t_{\alpha+1}-
(j_\alpha-t_{\alpha})-i_{\alpha+1}}+{\bf u}_{(4)}({\bf i},{\bf j})\cr
{\bf u}_{(4)}({\bf i},{\bf j})&=\sum_{\m=\alpha+1}^{\beta-1}
\e{R_{i_{\m}}(j_{\m}+1)+|i_{\m+1}-i_{\m}|-|i_{\m}-i_{\m-1}|}
+\e{j_{\beta}+1-i_{\beta-1}}-\E{\alpha+1}{n}}}
for $i_\alpha=0$
\eqn\oxrfa{
rf_{(4,1)}({\bf i})=
{1\over 2}(-1)^\alpha \sum_{j=2}^l (-1)^{a_1+\sum_{k=2}^{j-1}
(a_k+b_k)}\kie{b_j}}
and for $i_{\alpha}=1$
\eqn\oxrfb{\eqalign{
rf_{(4,2)}({\bf i})&=(-1)^\alpha \left({(-1)^{b_1}\over 4}+{1\over 2}
\sum_{j=2}^{l} (-1)^{\sum_{k=1}^{j-1}(a_k+b_k)}\kie{b_j}\right)\cr}}

{\bf 10.4 Character identities}

It remains to take the limit $L\rightarrow \infty$ to produce
character identities from the polynomial identities. These character
identities are somewhat simpler because ${\tilde f}_s(L,{\bf {\tilde
u}})$ vanishes as $L\rightarrow \infty.$ We remark that in this limit
the explicit dependence on $b$ vanishes and only the dependence on $r$
remains. All polynomial identities which have different values of $b$
but the same value of $r$ give the same character identity. Thus we
obtain the following fermionic forms for the bosonic forms of the
characters $q^{{1\over 2}
(\phi_{r({\bf j}),s}-\phi_{r(s),s})+c(j_s)}B_{r({\bf j}),s}(q).$
In all formulas of this section $r$ is a function of 
${\bf j}$ as  given in \results.

{\bf Cases 1 and 2a:}
\eqn\chara{F_{r({\bf j}),s}(q)=q^{c(j_{\m})}f_s({\bf e}_{j_{\m}}-\E{\m+1}{n})}

{\bf Case 2b:}
\eqn\charb{\eqalign{F_{r({\bf j}),s}(q)&=
q^{c(j_{\m})} \left( q^{-{\nu_0+1\over 4}+{3\over 4}\kio{\m}}
f_s(-\E{2}{n}+\e{1+j_{\m}})\right.\cr
&~~~~\left. +\theta(\m\geq 2)\sum_{i=2}^{\m} q^{-{\nu_0-\kio{i}\over 4}}
f_s(-\E{2}{i}+\e{-1+t_i}+\e{j_{\m}}-\E{\m+1}{n})\right)\cr}}

{\bf Case 3:}
\eqn\charc{F_{r({\bf j}),s}(q)=
q^{c_{(3)}({\bf j})} \sum_{i_1,\ldots,i_{\beta-1}=0,1 \atop i_0=0}
q^{rf_{(3,1)}({\bf i})}
f_s({\bf u}_{(3)}({\bf i},{\bf j}))}

{\bf Case 4 with $\alpha=1$}
\eqn\chard{\eqalign{F_{r({\bf j}),s}(q)&=
 q^{c_{(4)}({\bf j})} \left( \sum_{i_{2},\ldots,i_{\beta-1}=
0,1 \atop i_{1}=0}
q^{rf_{(4,1)}({\bf i})} f_s({\bf u}_{(4,1)}({\bf i},{\bf j}))\right.\cr
&~~~\left. +\sum_{i_{2},\ldots,i_{\beta-1}=0,1 \atop i_1=1}
q^{rf_{(4,2)}({\bf i})} 
f_s({\bf u}_{(4,3)}({\bf i},{\bf j})) \right) \cr}}

{\bf Case 4 with $\alpha\geq 2$}
\eqn\chare{\eqalign{F_{r({\bf j}),s}(q)&=
q^{c_{(4)}({\bf j})} \left( \sum_{i_{\alpha+1},\ldots,i_{\beta-1}=
0,1 \atop i_{\alpha}=0}
q^{rf_{(4,1)}({\bf i})} f_s({\bf u}_{(4,1)}({\bf i},{\bf j}))\right.\cr
&~~~+\sum_{i_{\alpha+1},\ldots,i_{\beta-1}=0,1 \atop i_\alpha=1}
\left( q^{rf_{(4,2)}({\bf i})} 
f_s({\bf u}_{(4,2)}({\bf i},{\bf j}))\right. \cr
&~~~~~~~~~~~~~~\left. 
\left. +q^{rf_{(4,2)}({\bf i})-{(-1)^{\alpha}\over 4}} 
f_s({\bf u}_{(4,3)}({\bf i},{\bf j})) \right) \right) \cr}}

\newsec{Reversed parity identities}

For all the identities  presented thus far we have started
with the term with $b=1$ where the fermionic polynomial ~\flq~has $m_{t_{1+n}}$
restricted to even values. Indeed $m_{t_{1+n}}$ has even parity for
all the results presented in this paper as long 
as $b\leq l^{(n)}_{1+t_{1+n}}=p'-2y_n$.
These fermionic sums are what were called even sums 
~$F_{r(b),s}^{(e)}(L,b)$
in ~\rbm. According to theorems 1 and 2 of sec.8 it is equally
possible to start with the fermionic polynomial which appears in the
first equation in~\Ffidb~with $m_{t_{1+n}}$
odd  and the
entire construction of this paper may be carried out exactly as
before. The only difference is that the parity of $m_{t_{1+n}}$ will
be reversed in all the formulas presented earlier, which amounts to the
replacement ${\bf u}\rightarrow {\bf u}+{\bf e}_{1+t_{n+1}}.$ 
The step which must be different is the
analysis of the boundary conditions and the computation of the normalization
constant. This is  done below and  calling 
the reversed parity fermionic polynomials
$F^{(o)}_{r(b),s}(L,b)$ ~\Ffidb~ 
we find from each of the identities proven
above the corresponding reversed parity identity
\eqn\pariden{F_{r(b),s}^{(o)}(L,b)=q^{a({\bar s},r)}B_{r(b),{\bar s}}(L,b)}
where ${\bar s}=p'-s$ with $s=l^{(\mu_s)}_{1+j_s}$ 
and $a({\bar s},r)$ is given by the two
equivalent expressions
\eqn\paridena{\eqalign{a({\bar s},r)&=c(j_s)+c(1+t_{n+1})+c(t_n)-{1\over
2}\kie{n}\cr
&+{1\over 2}(\phi_{r,{\bar s}}-\phi_{r(s),s}+\phi_{1,s}-\phi_{p-1,{\bar
s}})\cr
&=c(j_s)+{1\over 2}(\phi_{r,{\bar s}}-\phi_{r(s),s})+{1\over
4}(\phi_{p-1,s}-\phi_{p-1,{\bar s}}+\phi_{1,s}-\phi_{1,{\bar s}}).\cr}}
We note that when $b=l^{(\mu)}_{1+j_\mu}$ ~\pariden~gives the result
announced in ~\rbm.

{\bf Computation of the normalization constant}

{}From sec. 9 it suffices to determine the constant $a({\bar s})$ in the
identity
\eqn\parrevc{F^{(o)}_{r(b),s}(L,b)=q^{a(\bar s)}{\tilde B}_{r(b),{\bar
s}}(L,b)}
since from the definition of ${\tilde B}_{r(b),s}(L,b)$~\mobose~it
follows that
\eqn\parrevd{a({\bar s},r)=a({\bar s})+{1\over 2}(\phi_{r,{\bar
s}}-\phi_{r({\bar s}),{\bar s}})}

In ~\parrevc~we first set $b=p'-1,~r=p-1$ to obtain
\eqn\parrevda{F^{(o)}_{p-1,s}(L,p'-1)=q^{a({\bar s})+{1\over
2}(\phi_{{p-1},{\bar s}}-\phi_{r({\bar s}),{\bar s}})}B_{{p-1},{\bar
s}}(L,p'-1).}
On the other hand from ~\prfb~and the last equation in ~\Ffidb~
\eqn\parreve{\eqalign{F^{(o)}_{p-1,s}(L,p'-1)&=q^{c(1+t_{1+n})+c(t_n)-{1\over
2}\kie{n}}F_{1,s}(L,1)\cr
&=q^{c(1+t_{1+n})+c(t_n)-{1\over 2}\kie{n}+c(j_s)+{1\over
2}(\phi_{1,s}-\phi_{r(s),s})}B_{1,s}(L,1)\cr}}
where in the last line we used~\finalid. Comparing ~\parrevd~and
\parreve~and using the symmetry of $B_{r(b),s}(L,b)$~\blsym~with $r(b)=b=1$
\eqn\parrevf{B_{1,s}(L,1)=B_{p-1,{\bar s}}(L,p'-1)}
we find the first expression for $a({\bar s})$
\eqn\parrevg{\eqalign{a({\bar s})&=
c(1+t_{1+t_n})+c(t_n)-{1\over 2}\kie{n}+c(j_s)\cr
&+{1\over 2}(\phi_{1,s}-\phi_{r(s),s}+\phi_{r({\bar s}),{\bar
s}}-\phi_{p-1,{\bar s}})\cr}.}

A second expression for $a({\bar s})$ is obtained by setting
$b=1,~r=1$ in ~\parrevc~
\eqn\parrevh{F_{1,s}^{(o)}(L,1)=q^{a(\bar s)+
{1\over 2}(\phi_{1,{\bar s}}-\phi_{r({\bar s}),{\bar s}})}
B_{1,\bar{s}}(L,1).}
However we also find from the second line  of ~\Ffida, the first
line of ~\Ffidb~and ~\finalid~
\eqn\parrevi{\eqalign{F_{1,s}^{(o)}(L,1)&=q^{-c(1+t_{1+n})-c(t_n)+{1\over
2}\kie{n}} F_{p-1,s}(L,p'-1)\cr
&=q^{-c(1+t_{1+n})-c(t_n)+{1\over 2}\kie{n}
+c(j_s)+{1\over 2}(\phi_{p-1,s}-\phi_{r(s),s})}B_{p-1,s}(L,p'-1).\cr}}
Comparing ~\parrevh~and ~\parrevi~and using the symmetry~\blsym~with
$s=b=1$ and $s\rightarrow p'-s$
\eqn\parrevk{B_{1,{\bar s}}(L,1)=B_{p-1,s}(L,p'-1)} 
we obtain
\eqn\parrevj{\eqalign{a({\bar s})&=c(j_s)-c(1+t_{1+n})-c(t_n)+{1\over
2}\kie{n}\cr
&+{1\over 2}(\phi_{p-1,s}-\phi_{r(s),s}+\phi_{r({\bar s}),{\bar
s}}-\phi_{1,{\bar s}}).\cr}}

Equations ~\parrevg~and ~\parrevj~imply the following identity 
\eqn\parrevk{c(1+t_{1+n})+c(t_n)-{1\over 2}\kie{n}={1\over
4}(\phi_{p-1,s}+\phi_{p-1,{\bar s}}-\phi_{1,s}-\phi_{1,{\bar s}})}
and using this identity we obtain~\paridena~from ~\parrevg~and ~\parrevd.

\newsec{The dual case $p<p'<2p$}

In the XXZ chain
~\hamxxz~we see from ~\deldef~that the regime $p<p'<2p$ corresponds to
$\Delta>0$ and since the XXZ chain has the symmetry
$H_{XXZ}(\Delta)=-H_{XXZ}(-\Delta)$ the $(m,n)$ systems for $\pm
\Delta$ are the same. Thus we consider the relation between $M(p,p')$
and $M(p'-p,p')$ which is obtained by the transformation $q\rightarrow
q^{-1}$ in the finite $L$ polynomials $F_{r(b),s}(L,b;q)$ 
and $B_{r(b),s}(L,b;q).$ 
To implement this transformation it is mandatory that the
$L\rightarrow \infty$ limit be taken only in the final step.

Consider first $B_{r(b),s}^{(p,p')}(L,b;q)$ given by~\polyfb~where we have
made the dependence on $p$ and $p'$ explicit. Then if we note that from the
definitions ~\newqbind~
\eqn\qmosym{{n+m\atopwithdelims[]
m}^{(0,1)}_{q^{-1}}=q^{-mn}{n+m\atopwithdelims[] m}_q^{(0,1)}}
we find
\eqn\bmosym{B^{(p,p')}_{r(b),s}(L,b,q^{-1})=q^{-{L^2\over
2}}q^{(b-s)^2\over 4}B_{b-r(b),s}^{(p'-p,p')}(L,b;q).}

Similarly we let $q\rightarrow q^{-1}$ in $f_s(L,{\bf  u};q)$ and
${\tilde f}_s(L,{\bf  \tilde u};q)$ given
by~\flq~and~\ftildedf~to obtain
\eqn\fmosym{\eqalign{f_s(L,{\bf u};q^{-1})&
=q^{-{L^2\over 4}-{L\over 2}\sum_{j=1}^{\nu_0}u_j}f_{s}^{(d)}(L,{\bf u};q)\cr
{\tilde f}_s(L,{\bf \tilde u};q^{-1})&
=q^{-{L^2\over 4}}{\tilde f}_{s}^{(d)}(L,{\bf \tilde u};q)\cr}}
where
\eqn\fmosyma{\eqalign{f_s^{(d)}(L,{\bf u};q)=&\sum_{\vm \in 2\Z^{t_{1+n}}
+\vw(u_{1+t_{n+1}},{\bu})}
q^{-{1\over 2}{\bf m^TMm}+{\bf A'}^T{\bf m}+C'}\cr
&~~~~~~~\times\prod_{j=1}^{t_{n+1}}
{(({\bf I}_{t_{n+1}}+{\bf M}){\bf m}+{{\bf{\bar u}}\over 2}+{L\over 2}{\bf
\bar e}_1)_j\atopwithdelims[] m_j}_q^{(1)},}}
$\bf M$ was defined in ~\mnsys, \syn,and the $t_{n+1}$-dimensional
vector ${\bf A}'$ is 
\eqn\apdefn{{\bf A}'={\bf A}^{'(b)}+{\bf A}^{'(s)}}
with 
\eqn\apdefna{\eqalign{A_k^{'(b)}&=\cases{-{1\over 2}u_k&for $k$ in an
odd zone\cr
0&for $k$ in an even zone}\cr
A_k^{'(s)}&=\cases{-{1\over 2}{\bar{u}}_k(s)&for $k$ in an even zone\cr
0&for $k$ in an odd zone\cr}}}
and
\eqn\ctnewdf{C'=-{1\over 8}\left({\bf \bar{u}}^T(s)_{+}{\bf B}{\bf
\bar{u}}(s)_{+}+{\bf {\bar y}}^T{\bf B}{\bf {\bar y}}\right)}
$\bf{\bar u}$ is given by \ubardf~and ~\sstring,~
and
\eqn\miwatau{{\bf {\bar y}} =\sum_{i=1}^{\nu_0}{\bf \bar e}_i u_i,}
\eqn\ftmosyma{\eqalign{
 {\tilde f}_s^{(d)}(L,{\bf e}_{\nu_0-j_0-1}&
-\e{\nu_0}+\vu'_1;q)\cr
&=\cases{ 
0&$j_0=\nu_0$\cr
q^{{L+j_0\over 2}}
[q^{-{2L+1\over 4}}f^{(d)}_s(L+1,{\bf e}_{\nu_0-j_0}-
{\bf e}_{\nu_0}+\vu'_1;q)&~\cr
~~~~~~~~~~-f^{(d)}_s(L,{\bf e}_{1+\nu_0-j_0}-
{\bf e}_{\nu_0}+\vu'_1;q)]&
$1\leq j_0 \leq \nu_0-1,$\cr
q^{{L\over2}}f^{(d)}_s(L,{\bf e}_{\nu_0-1}
-{\bf e}_{\nu_0}+\vu'_1,q)&~\cr
~~~~~~~~~~-(q^{L}-1)q^{-{1\over 4}}
f^{(d)}_s(L-1,\vu'_1;q)&for $j_0=0$}
}}
and the vector $\vu'_1\in Z^{1+t_{1+n}}$ 
was defined as any $1+t_{1+n}$-dimensional 
vector with no
components in the zero zone.
We note that for the vectors $\bf u$ in ~\polyform~that
$\sum_{j=1}^{\nu_0}u_j=0$ unless 
\eqn\novan{{\bf u}=-{\bf e}_{\nu_0}+{\bf u}'_1}
in which case the sum is $-1.$

 Thus since the $L^2$ dependent factors in the transforms ~\bmosym~and~\fmosym~
 are the same we
obtain from all the polynomial identities of $M(p,p')$ of the form~\finalid~
polynomial identities for  $M(p'-p,p')$ of the form
\eqn\traida{F_{r(b),s}^{(d)}(L,b;q)=q^{-{1\over
2}(\phi_{r(b),s}-\phi_{r(s),s})-{1\over 2}c(j_s)+{(b-s)^2\over
4}}B_{b-r(b),s}^{(p'-p,p')}(L,b;q)}
where $F_{r(b),s}^{(d)}(L,b;q)$ for $M(p'-p,p')$ is obtained from 
the corresponding $F_{r(b),s}(L,b;q)$ of sec. 10 as 
\eqn\traidb{F_{r(b),s}^{(d)}(L,b;q)=\sum_{{\bf u}\in U(b)}q^{-c_{\bf u}}
f_s^{(d)}(L,{\bf u};q)+\sum_{\tilde{\vu}\in \tilde{U}(b)}
q^{-{\tilde c}_{\bf{\tilde  u}}}
{\tilde f}_s^{(d)}(L,{\bf{\tilde  u}};q)}
where the exponents $c_{\bf u},~{\tilde c}_{\bf \tilde u}$ and the
sets $U(b),~{\tilde U}(b)$ as defined in ~\polyform~are explicitly 
given in sec. 10.

We note in particular that the term $-{L\over 2}\sum_{j=1}^{\nu_0}u_j$
in $f_s(L,{\bf u};q^{-1})$ in~\fmosym~is needed in case 4 of ~\results
~with $\alpha=1$ in order to cancel the explicit factor of
$q^{-{L\over 2}}$ in ~\oxRR.

In the limit $L\rightarrow \infty$ we find from~\fmosyma~
\eqn\duallima{\lim_{L\rightarrow \infty}f_s^{(d)}(L,{\bf
u};q)=f_s^{(d)}({\bf u};q)}
with
\eqn\duallimb{\eqalign{f_s^{(d)}({\bf u};q)&=\sum_{\vm \in 2\Z^{t_{1+n}}
+\vw(u_{1+t_{1+n}},{\bu})}
q^{-{1\over 2}{\bf m^TMm}+{\bf A'}^T{\bf m}+C'}\cr
&~~~~~~~~\times{1\over (q)_{m_1}}\prod_{j=2}^{t_{n+1}}
{(({\bf I}_{t_{n+1}}+{\bf M}){\bf m}
+{{\bf{\bar u}}\over 2})_j\atopwithdelims[] m_j}_q^{(1)}}}
and from ~\ftmosyma~
\eqn\duallimc{\lim_{L\rightarrow \infty}{\tilde f}_s^{(d)}(L,{\bf
e}_{\nu_0-j_0-1}-{\bf e}_{\nu_0}+\vu'_1;q)=\cases{q^{{j_0\over
2}-{1\over 4}}f_s^{(d)}({\bf
e}_{\nu_0-j_0}-{\bf e}_{\nu_0}+\vu'_1;q)& $j_0\neq
\nu_0$\cr
0&$j_0=\nu_0.$}} 
Thus we find from ~\traida, \traidb~the corresponding character
identities
\eqn\duallimd{F_{r(b),s}^{(d)}(b;q)=q^{-{1\over
2}(\phi_{r(b),s}-\phi_{r(s),s})-{1\over 2}c(j_s)+{(b-s)^2\over
4}}B_{b-r(b),s}^{(p'-p,p')}(b;q)}
and
\eqn\duallime{F_{r(b),s}^{(d)}(b;q)=\sum_{\vu\in U(b)}q^{-c_{\bf u}}
f_s^{(d)}({\bf u};q)+\sum_{\tilde{\vu}\in \tilde{U}(b)}
q^{-{\tilde c}_{\bf{\tilde u}}+{j_0\over 2}-{1\over 4}}
f_s^{(d)}(\tilde{\vu}';q)}
where $\tilde{\vu}'$ is obtained from \mirror~by replacing $j_0$ by $j_0-1$.
We note that in contrast to the case $p'>2p$ the terms involving
${\tilde {\bf u}}$ contribute in the $L\rightarrow \infty$ limit.

\newsec{Discussion}

The techniques developed in this paper are of great generality and
may be extended to derive many further results which extend those
presented  in the previous sections. We will thus conclude this
paper by discussing several extensions and applications which will be
elaborated elsewhere.

{\bf 13.1 Extension of the fundamental fermionic polynomials \flq}

There are two features in the definition of the polynomials
$f_s(L,{\bf u};q)$ which may strike one as perhaps arbitrary and
unmotivated; namely the fact that the vector ${\bf A}$ of
~\lincomplete-~\linterms~ is of a different form depending on whether
the components of 
$\vu(b)$ and $\bu{(s)}$ are in odd or even zones and the choice of the
branch $(0)$ or $(1)$ of the $q$--binomial coefficients and it may be
asked whether other choices are possible.

For the choice of branch of the $q$--binomials it is clear from appendix A
that the fermionic recursive properties of sec. 5 will be equally valid
if the interchange of branches $(0)\leftrightarrow (1)$ is made
everywhere. The only place the branch is explicitly used is in sec. 9
where the polynomials $F_{r(b),s}(L,b)$ are identified with a linear
combination of the bosonic polynomials ${\tilde B}_{r(b),s}(L,b).$ A
change of branch will lead to different linear combinations.

In a similar fashion different linear terms are possible which will
also lead to different linear combinations. An analogous phenomena was
studied in some detail for the model $SM(2,4\nu)$ in ~\rbmb.

{\bf 13.2 Extension to arbitrary s}

Throughout this paper we have for simplicity confined our attention
to values of $s$ given by a pure Takahashi length. However, an
examination of the proof of appendix A shows that in fact the 
fundamental fermionic recursion
relations are valid for all vectors ${\bf \bar{u}}(s).$  With this
generalization all values of $s$ can be treated in a manner parallel
to the way general values of $b$ were treated above and the final
result is a double sum over both sets of vectors ${\bf u}(b)$ and
${\bf \bar{u}}{(s)}$. This generalization will be presented in full elsewhere.

{\bf 13.3 Extension of m,n system}

The $(\vm,\vn)$ system~\mnsys~of this paper has the feature that the term
$L/2$ appears only in the first component of the equation. However it
is also useful to consider $L/2$ in other (possibly several)
positions. Each choice will lead to further identities. Certain
features of this extension have been seen in ~\rbmo, \ranne, \role, \rbg. 
As long as $L/2$ remains in the zero zone we are in the regime of weak
anisotropy in the sense of ~\rbgs~and we will obtain identities for
coset models    $(A_1^{(1)})_N\times(A_1^{(1)})_M/(A_1^{(1)})_{N+M}$
with  $N$ integer and $M$ fractional. However if $L/2$ is in a higher
zone we are in the region of strong anisotropy and 
identities for cosets with both fractional levels  can be
obtained. 

{\bf 13.4 Bailey pairs} 

We conclude by demonstrating that all polynomial identities of this
paper may be restated in terms of new Bailey pairs and thus, by means of
Bailey's lemma~\randpac,~\rbail,~\raab, we may use the results of this paper to produce
bose/fermi character identities for models other than $M(p,p')$.

{\bf Definition of Bailey pair}

A pair of sequences $(\alpha_j, \beta_j)$ is said to form a Bailey
pair relative to $a$ if

\eqn\baia{
\beta_n=\sum_{j=0}^{n}{\alpha_j\over (q)_{n-j}(aq)_{n+j}}.}
All of our polynomial identities may be cast into this form. This is
easily seen, following \rfoda,~by writing the bosonic polynomial~\polyfb~
(where we suppress the dependence of $r$ on $b$) 
with $L=2n+b-s$ (when $L+b-s$ even) as
\eqn\bailb{\eqalign{B_{r,s}(2n+b-s,b;q)=(q^{b-s+1})_{2n}
\sum_{j=-\infty}^{\infty}&\left(
{q^{j(jpp'+rp'-sp)}\over (q)_{n-jp'}(q^{b-s+1})_{n+jp'}}\right.\cr
&\left. -{q^{(jp+r)(jp'+s)}\over (q)_{n-(jp'+s)}
(q^{b-s+1})_{n+(jp'+s)}}\right).\cr}}
Thus comparing the identity~\rfbpolyid~with ~\baia~
we obtain the following Bailey pair relative to $a=q^{b-s}$
\eqn\bailc{\eqalign{\beta_n&=q^{-{1\over
2}(\phi_{r,s}-\phi_{r(s),s})-c(j_s)}F_{r,s}(2n+b-s,b,q)/(aq)_{2n}\cr
\alpha_{n}
&=\cases{q^{j(jpp'+rp'-sp)} & for $n=jp',~~(j\geq 0)$\cr
q^{j(jpp'-rp'+sp)} & for $n=jp'+s-b,~~(j\geq 1)$\cr
-q^{(jp+r)(jp'+s)} & for $n=jp'+s,~~(j\geq 0)$\cr
-q^{(jp-r)(jp'-s)} & for $n=jp'-b,~~(j\geq 1)$\cr
}
}}
If two of the restrictions in \bailc~are the same the formula should be read
as the sum of both.

The utility of the Bailey pair follows from the lemma of
Bailey~\rbail

{\bf Lemma}

Given a Bailey pair and a pair $
(\gamma_n,\delta_n)$ satisfying
\eqn\baile{\gamma_n=\sum_{j=n}^{\infty}{\delta_j\over
(q)_{j-n}(aq)_{j+n}}}
then
\eqn\bailf{\sum_{n=0}^{\infty}\alpha_n \gamma_n=\sum_{n=0}^{\infty} \beta_n
\delta_n.}

We note in particular two  sets of $(\gamma_n, \delta_n)$
pairs. One is the original pair of Bailey~\rbail~
\eqn\bailg{\eqalign{\gamma_n&={(\rho_1)_n(\rho_2)_n(aq/\rho_1\rho_2)^n\over
(aq/\rho_1)_n(aq/\rho_2)_n(q)_{N-n}(aq)_{N+n}}\cr
\delta_n&={(\rho_1)_n(\rho_2)_n(aq
\rho_1\rho_2)_{N-n}(aq/\rho_1\rho_2)^n\over (q)_{N-n}(aq/\rho_1)_N(aq/\rho_2)_N}}}
which has been used to produce characters of the $N=1$ and $N=2$
supersymmetric models~\rbmsb~ and a new pair of~\rsw~
which gives characters of the coset models
$(A_1^{(1)})_N\times(A_1^{(1)})_M/(A_1^{(1)})_{N+M}$ where $M$ may be
fractional.

The presentation of the detailed consequences of this observation are
too lengthy for inclusion in this paper and will be published
separately~\rbmsw.

\bigskip
\bigskip
{\bf Acknowledgments}

The authors are grateful to G.E. Andrews for his interest and encouragement, 
to K. Voss and S.O. Warnaar for discussions and 
careful reading of the manuscript, and
to T. Miwa for many helpful comments. 
One of us (AB) is pleased to acknowledge the hospitality of
the ITP of SUNY Stony Brook where part of this work was done.
This work is supported in part by the NSF under DMR9404747.
\vfill
\eject

\appendix{A} {Proof of fundamental fermionic recursive properties}

In this appendix we prove the fundamental fermionic recursive properties
\freczero~for $0\leq j_0<\nu_0$ and \frecrelc~in detail. 
The proof of the remaining cases can be carried out in a similar fashion and
will be left to the reader. 

We begin by introducing some shorthand notations.
We use the vectors $\vu_0(j_{\m})$ and $\vu_{\pm 1}(j_{\m})$ 
as defined in \uvec.
Furthermore we need the operator which projects onto the components
of a vector in an odd zone
\eqn\apro{\pi_i=\cases{1 & for $t_{1+2l}+1\leq i \leq t_{2+2l}$\cr
0 & otherwise}.}
and
\eqn\aproa{\tilde{\pi}=1-\pi}
which projects onto components of even zones. 
We also need the vector
\eqn\aO{\O_{a,b}=\cases{0 & if $b<a$\cr
\sum_{i=a}^b \pi_i{\bf \bar e}_i.}}

We also make the following definition which generalizes the fermionic
polynomial \flq
\eqn\aferm{q^{\delta(\vn,\vm)}\bra{\bf \bar A}{\bf \bar B}
=\sum_{\vn,m_{t_{1+n}}\equiv P 
(\mod 2)} q^{\delta(\vn,\vm)+\Phi_s(\vn,\vm,\vu,L)}
\prod_{j=1}^{t_{n+1}}\nom{n_j+m_j+{\bar A}_j}{n_j+{\bar B}_j}_q^{(0)}}
where $P=0,1$ with $P\equiv {u}_{1+t_{1+n}}$ (mod 2), and $n_j,m_j$ are related 
by \mnsys~for given $L$ and ${\bu}$
\eqn\aubar{{\bu}=\sum_{i=1}^{t_{n+1}}({\bf {u}}_0(j_\mu))_i {\bf \bar
e}_i+\bu(s),~~1+t_{\m}\leq j_{\m}\leq t_{1+\m}}
 and
\eqn\aphi{\Phi_s(\vn,\vm,\vu,L)=Q(\vn,\vm)+Lf(\vn,\vm),} 
where the quadratic term $Q$ is defined in \phib~and the linear term $Lf$ 
is defined through $\vA$ in \lincomplete-\linterms. We note in particular
\eqn\apart{\bra{0}{0}=f_s(L,\vu_0(j_{\m})).}

Finally we define
\eqn\amnu{\left\{\vn,\vm,{\bu}\right\}_L = {\rm~set~of~solutions~
to}~(\vm,\vn)-{\rm system~\mnsys~with}~L,{\bu}.}

It will be necessary to make variable changes in $\vn,\vm$ in \aferm. Hence
we will need some identities relating different objects. For example
it is easy to check that the set of solutions $\{\vn,\vm,{\bu}\}_L$
is equal to the set of solutions $\{\vn,\vm,{\bu}+{\bf \bar e}_1\}_{L-1}$. One
may also show that $\Phi_s(\vn,\vm,\vu,L)=\Phi_s(\vn,\vm,\vu+\e{1},L-1)$.
Thus
\eqn\aa{f_s(L,\vu)=f_s(L-1,\vu+\e{1})}
which proves \freczero~for $j_0=0$.
One may also show that
\eqn\aex{\eqalign{\{\vn,\vm,{\bu}\}_L-\{{\bf \bar e}_1,0,0\}
&= \{\vn,\vm,{\bu}\}_{L-2}\cr
m_1+\Phi_s(\vn,\vm,\vu,L)&=(L-1)+\Phi_s(\vn-{\bf \bar e}_1,\vm, {\bf
u}, L-2)}}
where $\vu$ is any vector such that $\sum_{i=1}^{\nu_0}u_i=0$. 
{}From this follows
\eqn\ab{q^{m_1}\bra{-{\bf \bar e}_1}{-{\bf \bar e}_1}=q^{L-1}f_s(L-2,\vu_0(j_{\m})).}

Similarly one can prove the identities
\eqn\ac{\eqalign{
&1:~\bra{-2\bEl{1}{\nu_0}-\be{\nu_0+1}}{-\be{\nu_0+1}}
=f_s(L-2,\vu_0(j_{\m}))\cr
&2:~\bra{-\be{\nu_0}-\bEl{1}{\nu_0+1}}{-\be{\nu_0+1}}=f_s(L-1,\e{\nu_0-1}
-\e{\nu_0}+\vu_0(j_{\m}))\cr
&3:~\bra{-\bEl{1}{j_0}}{0}=f_s(L-1,\vu_1(j_0))\cr
&4:~
q^{m_{j_0}-m_{j_0-1}+1}
\bra{-\bEl{1}{j_0-1}+\be{j_0-1}-\be{j_0}}{\be{j_0-1}-\be{j_0}}
=f_s(L-1,\vu_{-1}(j_0)),~~j_0>1\cr
&5:~q^{m_{j_0}-m_{j_0-1}+1}\bra{-2\bEl{1}{j_0-1}+\be{j_0-1}-\be{j_0}}
{\be{j_0-1}-\be{j_0}}
=f_s(L-2,\vu_{0}(j_0)),~~j_0>1\cr
&6:~q^{\O_{1,j_{\m}}^T\vn}\bra{-\bEl{1}{j_{\m}}}{0}
=q^{{L-1\over 2}-{\nu_0+1\over 4}
+{3\over 4}\kio{\m}}f_s(L-1,\vu_1(j_{\m})-\E{1}{\m}),~~0<\m\leq n\cr
&7:~q^{\O_{1,j_{\m}-1}^T\vn+(m_{j_{\m}}-m_{j_{\m}-1}+1)
\kie{\m}}\bra{-\bEl{1}{j_{\m}-1}
+\be{j_{\m}-1}-\be{j_{\m}}}{\be{j_{\m}-1}-\be{j_{\m}}}\cr
&~~~~~=
q^{{L-1\over 2}-{\nu_0-\kie{\m}\over 4}}f_s(L-1,\vu_{-1}(j_{\m})-\E{1}{\m}),~~
0<\m\leq n\cr
&8:~q^{\O_{1,t_k}^T\vn+(n_{t_k}+m_{1+t_k})\kie{k}}
\bra{-\be{t_k}-\bEl{1}{1+t_k}}{-\be{1+t_k}}\cr
&~~~~~=q^{{L-1\over 2}-{\nu_0-\kio{k}\over 4}}
f_s(L-1,\vu_0(j_{\m})+\e{-1+t_k}-\E{1}{k}),~~2\leq k\leq \m\leq n.
}}
The derivation of \ac~requires the verification 
the following intermediate equations
\eqn\achelp{\eqalign{
1':&~\{\vn,\vm,{\bu}\}_L-\{\be{\nu_0+1},2\bEl{1}{\nu_0},0\}
=\{\vn,\vm,{\bu}\}_{L-2}\cr
&~\Phi_s(\vn,\vm,\vu,L)=\Phi_s(\vn-\be{\nu_0+1},\vm-2\bEl{1}{\nu_0},\vu,L-2)\cr
2':&~\{\vn,\vm,{\bu}\}_L-\{\be{\nu_0+1},\be{\nu_0}+\bEl{1}{\nu_0},0\}
=\{\vn,\vm,\be{\nu_0-1}-\be{\nu_0}+{\bu}\}_{L-1}\cr
&~\Phi_s(\vn,\vm,\vu,L)=\Phi_s(\vn-\be{1+\nu_0},\vm-\be{\nu_0}-\bEl{1}{\nu_0},
\vu+\e{\nu_0-1}-\e{\nu_0},L-1)\cr
3':&~\{\vn,\vm,{\bu}\}_L-\{0,\bEl{1}{j_0},0\}=
\{\vn,\vm,{\bu}+\be{j_0+1}-\be{j_0}\}_{L-1}\cr
&~\Phi_s(\vn,\vm,\vu,L)=\Phi_s(\vn,\vm-\bEl{1}{j_0},\vu+\e{j_0+1}-\e{j_0},L-1)
\cr
4':&~\{\vn,\vm,{\bu}\}_L+\{\be{j_0-1}-\be{j_0},-\bEl{1}{j_0-1},0\}
=\{\vn,\vm,{\bu}+\be{j_0-1}-\be{j_0}\}_{L-1}\cr
&~(m_{j_0}-m_{j_0-1}+1)+\Phi_s(\vn,\vm,\vu,L)\cr
&~~~~~=\Phi_s(\vn+\be{j_0-1}-\be{j_0},
\vm-\bEl{1}{j_0-1},\vu+\e{j_0-1}-\e{j_0},L-1)\cr
5':&~\{\vn,\vm,{\bu}\}_L+\{\be{j_0-1}-\be{j_0},-2\bEl{1}{j_0-1},0\}
=\{\vn,\vm,{\bu}\}_{L-2},~~j_0>1\cr
&~\Phi_s(\vn,\vm,\vu_0(j_0),L)+1+m_{j_0}-m_{j_0-1}\cr
&~~~~~=\Phi_s(\vn+\be{j_0-1}-\be{j_0},\vm-2\bEl{1}{j_0-1},\vu_0(j_0),L-2)\cr
6':&~\{\vn,\vm,{\bu}\}_L-\{0,\bEl{1}{j_{\m}},0\}
=\{\vn,\vm,{\bu}+\be{j_{\m}+1}-\be{j_{\m}}-\bE{1}{\m}\}_{L-1}\cr
&~\Phi_s(\vn,\vm,\vu_0(j_{\m}),L)+\O_{1,j_{\m}}^T\vn\cr
&~~~~~={L-1\over 2}-{\nu_0+1\over 4}+{3\over 4}\kio{\m}
+\Phi_s(\vn,\vm-\bEl{1}{j_{\m}},\vu_1(j_{\m})-\E{1}{\m},L-1),~\m\not=0\cr
7':&~\{\vn,\vm,{\bu}\}_L+\{\be{j_{\m}-1}-\be{j_{\m}},-\bEl{1}{j_{\m}-1},0\}
=\{\vn,\vm,{\bu}+\be{j_{\m}-1}-\be{j_{\m}}-\bE{1}{\m}\}_{L-1}\cr
&~\Phi_s(\vn,\vm,\vu_0(j_{\m}),L)+\O_{1,j_\m-1}^T\vn+\kie{\m}(1+m_{j_{\m}}
-m_{j_{\m}-1})={L-1\over 2}\cr
&~~~~~-{\nu_0-\kie{\m}\over 4}+\Phi_s(\vn+\be{j_{\m}-1}-\be{j_{\m}},
\vm-\bEl{1}{j_{\m}-1},\vu_{-1}(j_{\m})-\E{1}{\m},L-1)\cr
8':&~\{\vn,\vm,{\bu}\}_L-\{\be{1+t_k},\be{t_k}+\bEl{1}{t_k},0\}
=\{\vn,\vm,{\bu}+\be{-1+t_k}-\bE{1}{k}\}_{L-1}\cr
&~\Phi_s(\vn,\vm,\vu_0(j_{\m}),L)+\O_{1,t_k}^T\vn+(n_{t_k}+m_{1+t_k})
\kie{k}={L-1\over 2}\cr
&~~~~~-{\nu_0-\kio{k}\over 4}+
\Phi_s(\vn-\be{1+t_k},\vm-\be{t_k}-\bEl{1}{t_k},\vu_0(j_{\m})+\e{-1+t_k}
-\E{1}{k},L-1)
}}
where ${\bf {\bar E}}_{a,b},
~ {\bf E}^{(t)}_{a,b},~{\bf\bar E}^{(t)}_{a,b}$ are given 
by~\eevecdef ~\etdefn~and ~\miwabar.

After these preliminaries we are in the position to prove
the fermionic recurrences \freczero~for $0\leq j_0<\nu_0$ and \frecrelc. 
Our method is the technique of telescopic expansion introduced in
\rb~and further developed in \rbmo~and \ranne. With this technique we may
expand $f_s(L,\vu_0(j_{\m}))$
\eqn\atele{
f_s(L,\vu_0(j_{\m}))= \bra{0}{0}
=q^{\O_{1,j_{\m}}^T\vn}\bra{-\bEl{1}{j_{\m}}}{0}
+\sum_{l=1}^{j_{\m}}q^{\O_{1,l}^T\vn
+\tilde{\pi}_l m_l}\bra{-\bEl{1}{l}}{-\be{l}}.
}
where \qbinrra~ has been used in odd zones and \qbinrrb~has been used in even
zones.

We shall further need the telescopic expansions
\eqn\ateleb{\eqalign{
&q^{\kio{k}n_{t_{1+k}}} \bra{-\be{t_{1+k}}-
\bEl{1}{1+t_{1+k}}}{-\be{1+t_{1+k}}}=\cr
&~~~\sum_{l=d_k}^{-1+t_{1+k}} q^{\O_{l+1,t_{1+k}}^T\vn+(m_l-1)\kie{k}}
\bra{-\bEl{1}{t_{1+k}}-\bEl{l}{1+t_{1+k}}}{-\be{l}-\be{1+t_{1+k}}}\cr
&~~~+q^{\O_{d_k,t_{1+k}}^T\vn}\bra{-\bEl{1}{t_{1+k}}-\bEl{d_k}{1+t_{1+k}}}
{-\be{1+t_{1+k}}}
}}
and
\eqn\atelec{\eqalign{
&q^{\kie{\m}}\bra{-\bEl{1}{j_{\m}-1}+\be{j_{\m}-1}-\be{j_{\m}}}
{\be{j_{\m}-1}-\be{j_{\m}}}=\cr
&~~~\sum_{l=d_k}^{j_{\m}-1} q^{\O_{l+1,j_{\m}-1}^T(\vn+\be{j_{\m}-1})
+\tilde{\pi}_l m_l}
\bra{-\bEl{1}{j_{\m}}-\bEl{l}{j_{\m}-2}}{-\be{l}+\be{j_{\m}-1}-\be{j_{\m}}}\cr
&~~~+q^{\kie{\m}+\O_{1+t_{\m},j_{\m}-1}^T(\vn+\be{j_{\m}-1})}
\bra{-\bEl{1}{j_{\m}-2}-\bEl{1+t_{\m}}{j_{\m}}}
{\be{j_{\m}-1}-\be{j_{\m}}}
}}
with $d_k=1+t_k+\delta_{k,0}$, where again \qbinrra~has been used in odd zones 
and \qbinrrb~has been used in even zones.

Using case 6 in \ac~we may rewrite \atele~as
\eqn\atelea{
f_s(L,\vu_0(j_{\m}))=q^{{L-1\over 2}-{\nu_0+1\over 4}+{3\over 4}\kio{\m}}
f_s(L-1,\vu_1(j_{\m})-\E{1}{\m})+\sum_{k=0}^{\m} I_k}
where
\eqn\adefI{I_k=\sum_{l=d_k}^{g_k} q^{\O_{1,l}^T\vn+\tilde{\pi}_l m_l}
\bra{-\bEl{1}{l}}{-\be{l}},~~~0\leq k\leq \m}
and $g_k=\cases{t_{1+k} & for $k<\m$\cr
j_{\m} & for $k=\m$}$
and may now process the terms $I_k~(0\leq k\leq \m-1)$ as follows. In the 
$l$th term ($l\not= d_k$) in \adefI~we perform the change of variables
\eqn\avca{\eqalign{
\vn &\to \vn+\be{l}-\be{l-1}-\be{1+t_{1+k}}\cr
\vm &\to \vm-2\bEl{l}{t_{1+k}}
}}
to obtain
\eqn\aI{\eqalign{
I_k&=q^{\O_{1,t_k}^T\vn+m_{d_k}\kie{k}}\bra{-\bEl{1}{d_k}}{-\be{d_k}}\cr
&+\sum_{l=d_k}^{-1+t_{1+k}}q^{(\O_{1,t_{1+k}}^T+\O_{1+l,t_{1+k}}^T)\vn
+m_{1+t_{1+k}}\kio{k}+(m_l-1)\kie{k}}\bra{-\bEl{l}{t_{1+k}}-\bEl{1}{1+t_{1+k}}}
{-\be{l}-\be{1+t_{1+k}}}.
}}

For $k=0,~\m>0$ we use \ab~and \ac~1,2 and obtain
\eqn\aIo{\eqalign{
I_0&=q^{m_1}\bra{-\be{1}}{-\be{1}}+\bra{-\be{\nu_0}-
\bEl{1}{1+\nu_0}}{-\be{\nu_0+1}}
-\bra{-2\bEl{1}{\nu_0}-\be{1+\nu_0}}{-\be{1+\nu_0}}\cr
&=f_s(L-1,\e{\nu_0-1}-\e{\nu_0}+\vu_0(j_{\m}))
+(q^{L-1}-1)f_s(L-2,\vu_0(j_{\m})).
}}

Applying the telescopic expansion \ateleb~to the last sum in \aI, we obtain
for $k\not=0$
\eqn\aIa{\eqalign{
I_k&=q^{\O_{1,t_k}^T\vn+m_{d_k}\kie{k}}\bra{-\bEl{1}{d_k}}{-\be{d_k}}\cr
&+q^{\O_{1,t_{1+k}}^T\vn+m_{1+t_{1+k}}\kio{k}} \left(
q^{n_{t_{1+k}}\kio{k}}\bra{-\be{t_{1+k}}-\bEl{1}{1+t_{1+k}}}{-\be{1+t_{1+k}}}
\right.\cr
&~~~~~~~~~~~~~\left. -q^{\O_{d_k,t_{1+k}}^T\vn}
\bra{-\bEl{1}{t_{1+k}}-\bEl{d_k}{1+t_{1+k}}}{-\be{1+t_{1+k}}}\right).
}}

For $k\not=0$, we perform one more change of variables 
\eqn\avcb{\eqalign{
&\vn\to \vn+\be{1+t_{1+k}}-\be{t_k}-\be{1+t_k}\cr
&\vm\to \vm+2\bEl{1+t_k}{t_{1+k}}
}}
in the last term of \aIa~to arrive at
\eqn\aIb{\eqalign{
I_k&=q^{\O_{1,t_{1+k}}^T\vn+(m_{1+t_{1+k}}+n_{t_{1+k}})\kio{k}}
\bra{-\be{t_{1+k}}-\bEl{1}{1+t_{1+k}}}{-\be{1+t_{1+k}}}\cr
&+q^{\O_{1,t_k}^T\vn+m_{1+t_k}\kie{k}}\bra{-\bEl{1}{1+t_k}}{-\be{1+t_k}}\cr
&-q^{\O_{1,t_k}^T\vn+m_{1+t_k}\kie{k}+(m_{t_k}-1)\kio{k}}
\bra{-\be{t_k}-\bEl{1}{1+t_k}}{-\be{t_k}-\be{1+t_k}}\cr
&=\sum_{a=k}^{k+1} q^{\O_{1,t_a}^T\vn+(n_{t_a}+m_{1+t_a})\kie{a}}
\bra{-\be{t_a}-\bEl{1}{1+t_a}}{-\be{t_a}-\be{1+t_a}}
}}
where to get the last line \qbinrra~or \qbinrrb~has been used for $k$ even or
odd, respectively. Thus
recalling \ac~2,8 we have for $1\leq k\leq \m-1$
\eqn\aIc{
I_k=\sum_{a=k}^{k+1}q^{(1-\delta_{a,1})({L-1\over 2}-{\nu_0-\kio{a}
\over 4})}f_s(L-1,\e{t_a-1}-\E{1}{a}+\vu_0(j_{\m})).}

It remains to process $I_{\m}$. To this end we perform the change of variables
\eqn\avcc{\eqalign{
\vn\to \vn -\be{l-1}+\be{l}+\be{j_{\m}-1}-\be{j_{\m}}\cr
\vm\to \vm-2\bEl{l}{j_{\m}-1}
}}
in the $l$th term ($l\not=d_{\m}$) appearing in the sum \adefI~with $k=\m$. 
This yields
\eqn\aId{\eqalign{
I_{\m}&=q^{\O_{1,j_{\m}-1}^T\vn+m_{d_k}\kie{\m}}\bra{-\bEl{1}{d_{\m}}}
{-\be{d_{\m}}}\cr
&+\sum_{l=d_{\m}}^{j_{\m}-1} q^{\O_{1,j_{\m}-1}^T\vn
+\O_{l+1,j_{\m}-1}^T(\vn+\e{j_{\m}-1})
+(m_l+m_{j_{\m}}-m_{j_{\m}-1})\kie{\m}}\bra{-\bEl{1}{j_{\m}}
-\bEl{l}{j_{\m}-2}}{-\be{l}+\be{j_{\m}-1}-\be{j_{\m}}}
}}
where we used the empty sum convention. Using the telescopic expansion
\atelec~we find
\eqn\aIe{\eqalign{
I_{\m}&=q^{\O_{1,j_{\m}-1}^T\vn+m_{d_k}\kie{\m}}\bra{-\bEl{1}{d_{\m}}}
{-\be{d_{\m}}}\cr
&+q^{\O_{1,j_{\m}-1}^T\vn+(m_{j_{\m}}-m_{j_{\m}-1}+1)\kie{\m}} \left(
\bra{-\bEl{1}{j_{\m}-1}+\be{j_{\m}-1}-\be{j_{\m}}}{\be{j_{\m}-1}-\be{j_{\m}}}
\right.\cr
&~~~~~~~~~~~\left.-q^{\O_{1+t_{\m},j_{\m}-1}^T(\vn+\be{j_{\m}-1})}
\bra{-\bEl{1}{j_{\m}-2}-\bEl{1+t_{\m}}{j_{\m}}}{\be{j_{\m}-1}-
\be{j_{\m}}} \right)
}}
For $\m=0$ \aIe~yields
\eqn\aIf{\eqalign{
I_0&=q^{m_1}\bra{-\be{1}}{-\be{1}}+\theta(j_0>1)q^{m_{j_0}-m_{j_0-1}+1}\cr
&~~~~\times\left( \bra{-\bEl{1}{j_0-1}+\be{j_0-1}-\be{j_0}}
{\be{j_0-1}-\be{j_0}}
-\bra{-2\bEl{1}{j_0-1}+\be{j_0-1}-\be{j_0}}{\be{j_0-1}-\be{j_0}}\right)\cr
&=f_s(L-1,\vu_{-1}(j_0))+(q^{L-1}-1)f_s(L-2,\vu_{0}(j_0))
}}
where we used \aa, \ab, \ac~4,5.

For $\m>0$ we need to change variables again
\eqn\acvd{\eqalign{
&\vn\to \vn -\be{t_{\m}}-\be{1+t_{\m}}-\be{j_{\m}-1}+\be{j_{\m}}\cr
&\vm\to\vm+2\bEl{1+t_{\m}}{j_{\m}-1}
}}
in the last term of the right hand side of \aIe~to get
\eqn\aIg{\eqalign{
I_{\m}&=q^{\O_{1,j_{\m}-1}^T\vn+m_{1+t_{\m}}\kie{\m}}
\bra{-\bEl{1}{1+t_{\m}}}{-\be{1+t_{\m}}}\cr
&+q^{\O_{1,j_{\m}-1}^T\vn+(m_{j_{\m}}+m_{j_{\m}-1}+1)\kie{\m}}
\bra{-\bEl{1}{j_{\m}-1}+\be{j_{\m}-1}-\be{j_{\m}}}{\be{j_{\m}-1}-
\be{j_{\m}}}\cr
&-q^{\O_{1,j_{\m}-1}^T\vn+m_{1+t_{\m}}\kie{\m}+(m_{t_{\m}}-1)\kio{\m}}
\bra{-\be{t_{\m}}-\bEl{1}{1+t_{\m}}}{-\be{t_{\m}}-\be{1+t_{\m}}}\cr
&=q^{\O_{1,j_{\m}-1}^T\vn+(m_{j_{\m}}-m_{j_{\m}-1}+1)\kie{\m}}
\bra{-\bEl{1}{j_{\m}-1}+\be{j_{\m}-1}-\be{j_{\m}}}{\be{j_{\m}-1}-
\be{j_{\m}}}\cr
&+q^{\O_{1,j_{\m}-1}^T\vn+(m_{1+t_{\m}}+n_{t_{\m}})\kie{\m}}
\bra{-\be{t_{\m}}-\bEl{1}{1+t_{\m}}}{-\be{1+t_{\m}}}
}}
where we used the elementary binomial recurrences \qbinrra~and \qbinrrb.
Recalling equations 7,8 of ~\ac~we finally obtain for $\m>0$
\eqn\afin{\eqalign{
I_{\m}&=q^{{L-1\over 2}-{\nu_0-\kio{\m}\over 4}}f_s(L-1,\vu_0(j_{\m})+
\e{t_{\m}-1}
-\E{1}{\m})\cr
&+q^{{L-1\over 2}-{\nu_0-\kie{\m}\over 4}}f_s(L-1,\vu_{-1}(j_{\m})-\E{1}{\m}).
}}

The desired results \freczero~(for $2\leq j_0\leq \nu_0-1$) 
and \frecrelc~follow
by combining formulas \atelea, \aIo, \aIc, \aIf~and \afin.

\appendix{B}{Proof of the initial conditions for Propositions 1 and 2}

In this appendix we prove the initial conditions for propositions 1 and 2
as discussed in sec. 8. These proofs explain the origin of
the factors $q^{{a\over 2}(-1)^{\m}}$ in propositions 1 and 2 and also
demonstrate that at the dissynchronization point the form
~\polyform~does not hold.

{\it Proof of proposition 1 with $a=-1$ and $\m=2$}

To prove proposition 1 with $a=-1,~\mu=2$ we will use the following
decomposition of the flow $\ev{y_2-2}{-1}$
\eqn\demappa{\psfig{file=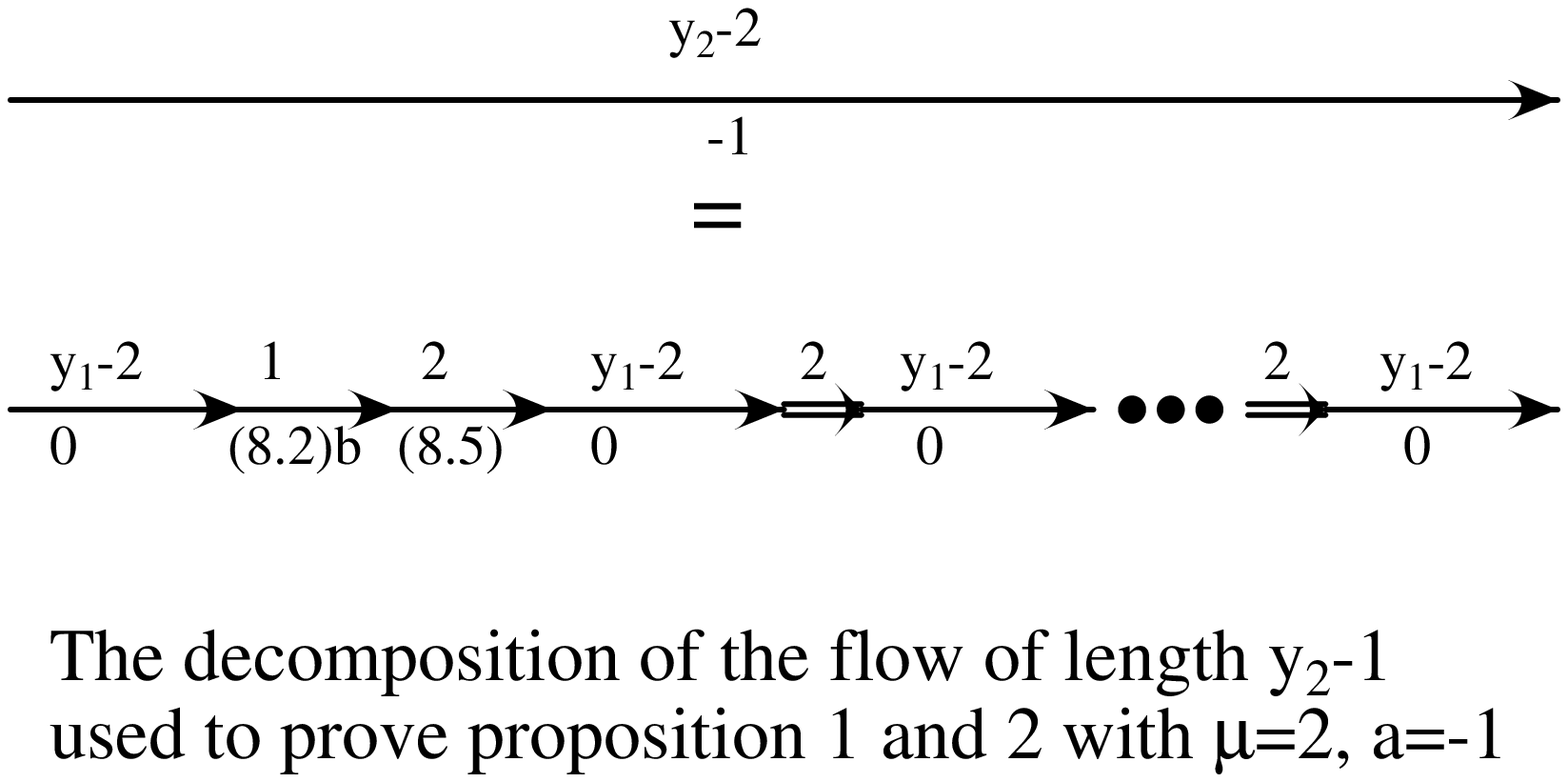,width=3.5in}}
whose proof is elementary and will be omitted.

We recall that  $\u{2}$ is any $1+t_{n+1}$--dimensional vector 
with $(\u{2})_j=0$ for 
$j\leq t_2.$ Then first using proposition 1 with $\mu=1,~a=0$ and then 
using ~\freczero~ with $j_0=\nu_0$
we obtain 
\eqn\pia{\eqalign{
& \left\{ f_s(L,-\e{t_1}-\e{t_2}+\u{2}), f_s(L,\e{1}-\e{t_1}-\e{t_2}+\u{2})
\right\} \cr
&\ev{y_1-2}{0} 
\left\{ f_s(L,\e{-1+t_1}-\e{t_1}-\e{t_2}+\u{2}),
f_s(L,-\e{t_2}+\u{2}) \right\}\cr
&\ev{1}{\pprbb b}
\left\{ f_s(L,-\e{t_2}+\u{2}),q^{-{L\over 2}}f_s(L,\e{1+t_1}-\e{t_2}+\u{2})
+(q^{L\over 2}-{1\over q^{L\over 2}})f_s(L-1,-\e{t_2}+\u{2}) \right\}
}}
It may be seen that the second polynomial in the last pair 
of ~\pia~for which $b$ is at the dissynchronization value $b_{2,-1}$ 
is not of the form~\polyform.

We may further evolve using flow $\ev{2}{\newgra}$ to find  
\eqn\pic{\eqalign{
& \left\{ f_s(L,-\e{t_2}+\u{2}),q^{-{L\over 2}}f_s(L,\e{1+t_1}-\e{t_2}+\u{2})
+(q^{L\over 2}-{1\over q^{L\over 2}})f_s(L-1,-\e{t_2}+\u{2}) \right\} \cr
\ev{2}{\newgra } & q^{-{1\over 2}}\left( q^{-{\nu_0-2\over 4}}\left\{
 f_s(L,-\e{t_1}+\e{2+t_1}-\e{t_2}+\u{2}),f_s(L,\e{1}-\e{t_1}+\e{2+t_1}-\e{t_2}+\u{2})\right\}\right.\cr
&\left.+\left\{P_0(L,1,\e{1+t_1}-\e{t_2}+\u{2}), 
P_1(L,1,\e{1+t_1}-\e{t_2}+\u{2}), \right\}\right)
}}
where we used ~\frecrelb~ with $j_1=1+t_1$ and the first equation in
~\freczero~ 
as well as the definitions ~\ftildedf,~\ppms,~\ppmsh. 
Now applying propositions 1 and 2 with $a=0$ and $\m=1$ 
for $\ev{y_1-2}{0}$ and lemma 2.1 for $\evs{2}$ according to ~\demappa~
we finally obtain upon recalling the definition ~\ppmsl
\eqn\pid{
 \left\{ f_s(L,-{\bf E}^{(t)}_{1,2}  +\u{2}), 
f_s(L,-{\bf E}^{(t)}_{1,2}-\e{t_2}+\u{2})
\right\}\ev{y_2-2}{-1} q^{-{1\over 2}} 
q^{c(t_2)}\left\{ P_1(L,2,{\bf u}'_2),f_s(L,\u{2}) \right\}
}
Therefore proposition 1 for $a=-1$ and $\m=2$ is proven.
We note that the factor $q^{-{1\over 2}}$ is the origin of the factor 
$q^{{a\over 2}(-1)^{\m}}$ in proposition 1 with $a=-1$.

{\it Proof of proposition 2 with $a=-1$ and $\m=2$}

For the proof of proposition 2 with $a=-1,~\mu=2$ we shall again use
the decomposition~\demappa. However, in the present case it is
convenient to  rewrite the piece $\ev{1}{\pprbb b}\ev{2}{\newgra}$ in
the bottom graph of ~\demappa~using the identity
\eqn\demappb{\ev{1}{\pprbb b}\ev{2}{\newgra}=\evs{2}\ev{1}{\pprbb a}}
To simplify our analysis, let us further assume that $\nu_1\geq
3,~\nu_0\geq 2.$ We start by applying propositions 
1 and 2 with $a=0$ and $\m=1$
\eqn\pie{\eqalign{
& \left\{ P_0(L,2,\u{2}),P_1(L,2,\u{2}) \right\} \cr
= &  q^{-{\nu_0\over 4}}\left\{ f_s(L,\e{-1+t_2}-{\bf
E}^{(t)}_{1,2}+\u{2}),f_s(L,{\bf e}_1+\e{-1+t_2}-{\bf
E}^{(t)}_{1,2}+\u{2})\right\}\cr
&+\left\{P_0(L,1,\u{2}),P_1(L,1,\u{2}) \right\}\cr
&\ev{y_1-2}{0}  
q^{-{\nu_0\over 4}}\left\{P_{-1}(L,1,{\bf e}_{-1+t_2}-{\bf
e}_{t_2}+\u{2})+q^{{L\over 2}-{\nu_0-2\over 4}}f_s(L,
-{\bf e}_{t_1}+\u{2})\right.,\cr
&~~~~\left. f_s(L,{\bf e}_{-1+t_2}-{\bf e}_{t_2}+\u{2})\right\}\cr}}

Evolving this once more according to $\evs{2}$ 
and using lemma 2.2 with $\mu=1,~j_1=1+t_1$ we find from ~\pie~
\eqn\pif{\eqalign{&\left\{P_0(L,2,\u{2}),P_1(L,2,\u{2})\right\}\aw{y_1}\cr
& q^{-{\nu_0\over 2}}\left\{f_s(L,\e{-2+t_2}-{\bf E}^{(t)}_{1,2}+\u{2}),
 f_s(L,\e{1}+\e{-2+t_2}-{\bf E}^{(t)}_{1,2}+\u{2})\right\}\cr
&+q^{-{\nu_0\over 4}}
\left\{P_0(L,1,\e{-1+t_2}-\e{t_2}+\u{2}),
P_1(L,1,\e{-1+t_2}-\e{t_2}+\u{2}) \right\}
}}
Evolving now according to  $\ev{1}{\pprbb a}$  and using the fermionic
recurrence ~\freczero~ for $f_s$ and ~\reccft~ for ${\tilde f}_s$ as well as
definitions ~\ppms,~\ppmsh~ we derive for $\nu_0\geq 2,~\nu_1 \geq 3$
\eqn\pig{\eqalign{
&\left\{P_0(L,2,\u{2}),P_1(L,2,\u{2})\right\}\ev{y_1+1}{(B.1)}\cr
&q^{-{\nu_0\over 2}} \left\{f_s(L,{\bf e}_1+{\bf e}_{-2+t_2}-{\bf
E}^{(t)}_{1,2}+\u{2}),f_s(L,{\bf e}_2+{\bf e}_{-2+t_2}-{\bf
E}^{(t)}_{1,2}+\u{2})\right\}\cr
& +q^{-{\nu_0\over 4}}
\left\{{\tilde f}_s(L,{\bf e}_{-2+t_1}+{\bf e}_{-1+t_2}-{\bf
E}^{(t)}_{1,2}+\u{2}),{\tilde f}_s(L,{\bf e}_{-3+t_1}+{\bf e}_{-1+t_2}-{\bf
E}^{(t)}_{1,2}+\u{2})\right\}}}
where the arrow in ~\pig~represents the flow ~$\ev{y_2-2}{-1}$
restricted to the first $1+y_1$ steps.

We now note that the polynomials in the rhs of ~\pig~and in the rhs of
~\initiald~differ only by the overall factor of $q^{-1/2}.$ Furthermore,
according to ~\ppdis~the flow $\ev{y_2-2}{-1}$ and the flow
$\ev{y_2-2}{0}$ have the same steps starting with the step $y_1+2.$
These observations along with proposition 2 for $a=0,~\mu=2$ clearly
imply that
\eqn\pij{\left\{ P_0(L,2,\u{2}),P_1(L,2,\u{2}) \right\} 
\ev{y_2-2}{-1}  q^{-{1\over 2}} \left\{
q^{c(t_2)+{L-\nu_0\over 2}} f_s(L,-{\bf E}^{(t)}_{1,2}+\u{2}),0 \right\}}
which concludes the proof of proposition 2 with $a=-1,~\m=2,~\nu_0\geq
2,~\nu_1\geq 3$. The proof of the remaining cases with $\nu_1=2$ and $\nu_1\geq 3,
\nu_0=1$ can be carried  out in a similar fashion and is left to the
reader. Finally, we note that the factor $q^{-1/2}$ in ~\pij~is the
origin of the factor $q^{a{(-)^{\mu}\over 2}}$ in proposition 2 with
$a=-1,~\mu.$

{\it Proof of proposition 1 with $a=-1$ and $\m=3$}

To prove proposition 1 with $a=-1,~\mu=3$ we shall use a
decomposition of $\ev{{y_3}-2}{-1}$ obtained from~\pppropa~with
$\mu=2,~a=-1$ by replacing the second arrow $\evs{2}$ from the left by
the composite arrow $\ev{1}{\pprbb a}\ev{1}{\pprbb b}.$ This
replacement is necessary because ~\pppropa~is not valid as it stands
for $\mu=2.$ In what follows the notation $\ev{b-2}{-1},~y_2< b \leq
y_3$ is understood as the flow $\ev{y_3-2}{-1}$ restricted to the
first $b-2$ steps.

For convenience we consider $\nu_0\geq 2$ and we start 
by applying proposition 1 with $\mu=2,~a=0$ and subsequently
using ~\anotherthing~  with $\mu=2,~\u{2}=-{\bf e}_{t_3}+\u{3}$ and lemma 2.1 with
$\mu=1,~j_1=t_2$ to obtain
\eqn\pik{\eqalign{
& \left\{ f_s(L,-{\bf E}^{(t)}_{1,3}+\u{3}), 
f_s(L,\e{1}-{\bf E}^{(t)}_{1,3}+\u{3})\right\}\ev{y_2}{-1} \cr
  & q^{c(t_2)-{\nu_0-1\over 4}} 
\left\{f_s(L,-{\bf e}_{t_1}+\e{1+t_2}-{\bf e}_{t_3}+\u{3}),
f_s(L,\e{1}-{\bf e}_{t_1}+\e{1+t_2}-\e{t_3}+\u{3})\right\}\cr
&+q^{c(t_2)}\left\{ P_{0}(L,1,-\e{t_3}+\u{3}),
P_{1}(L,1,-\e{t_3}+\u{3})\right\}.\cr}} 
Application of propositions 1 and 2 
with $a=0$ and $\m=1$ to the rhs of~\pik~yields
\eqn\pil{\eqalign{&\left\{f_s(L,-{\bf E}^{(t)}_{1,3}+\u{3}),f_s(L,{\bf
e}_{1}-{\bf E}^{(t)}_{1,3}+\u{3})\right\}\ev{y_2+y_1-2}{-1}\cr
&q^{c(t_2)-{\nu_0-1\over 4}}  \left\{ 
f_s(L,\e{-1+t_1}-\e{t_1}+\e{1+t_2}-\e{t_3}+\u{3})
+q^{{L\over 2}-{\nu_0-1\over 4}} f_s(L,-\e{t_1}-\e{t_3}+\u{3}),\right.\cr
& ~~~~~~~~\left. f_s(L,\e{1+t_2}-\e{t_3}+\u{3}) \right\}
}}
Next evolving the rhs of ~\pil~according to $\ev{1}{\pprbb a}$ and
making use of ~\frecrelb~ with $j_2=1+t_2$ we derive
\eqn\pim{\eqalign{
& \left\{ f_s(L,-{\bf E}^{(t)}_{1,3}+\u{3}), 
f_s(L,\e{1}-{\bf E}^{(t)}_{1,3}+\u{3})\right\} \cr
&\ev{y_2+y_1-1}{-1}  q^{c(t_2)-{\nu_0-1\over 4}}
\left\{ f_s(L,\e{1+t_2}-\e{t_3}+\u{3}),
q^{{L\over 2}-{\nu_0+1\over 4}} f_s(L,\e{2+t_2}-{\bf E}^{(t)}_{1,3}
+\u{3})\right.\cr
&~~\left. +q^{{L\over 2}-{\nu_0\over 4}} f_s(L,\e{-1+t_2}
+\e{1+t_2}-{\bf E}^{(t)}_{1,3}+\u{3})+
 f_s(L,{\bf e}_{-1+t_1}+{\bf e}_{t_2}+{\bf e}_{1+t_2}
-{\bf E}^{(t)}_{1,3}+\u{3}).\right\}}}
To proceed further we observe that the first two steps of
$\ev{1}{\pprbb b}\ev{y_2-2}{1}$ are the same as $\evs{2}.$ The
straightforward use of properties ~\freczero,~\reccft~and 
definition ~\ftildedf~ then gives
\eqn\appbnew{\left\{f_s(L,-{\bf E}^{(t)}_{1,3}+\u{3}), f_s(L,{\bf
e}_{1}-{\bf E}^{(t)}_{1,3}+\u{3})\right\}\ev{y_2+y_1+1}{-1}
q^{1\over 2}\left\{z_1(L,\u{3}),z_2(L,\u{3})\right\}}
where
\eqn\appbnewb{\eqalign{&z_{a}(L,\u{3})=
q^{c(t_2)-{\nu_0-1\over 4}}\left(
q^{-{\nu_0+1\over 4}} f_s(L,\e{a}+\e{2+t_2}-{\bf E}^{(t)}_{1,3}+\u{3})+
\right. \cr 
&\left.q^{-{\nu_0\over 4}}  
f_s(L,\e{a}+\e{-1+t_2}+\e{1+t_2}-{\bf E}_{1,3}^{(t)}
+\u{3})+\tilde{f}_s(L,\e{t_1-1-a}+\e{t_2}+\e{1+t_2}-{\bf
E}_{1,3}^{(t)}+\u{3})\right) \cr}}
with $a=1,2.$ On the other hand, as can be seen from case 3 of~\tzcc~
with  $j_2=t_2+1,~j_0=1,2$ and $-{\bf E}_{4,n}^{(t)}$ replaced by
$\u{3}$  
\eqn\pin{
 \left\{ f_s(L,-{\bf E}_{1,3}^{(t)}+\u{3}), 
f_s(L,\e{1}-{\bf E}_{1,3}^{(t)}+\u{3})\right\} 
\aw{y_2+y_{1}+1}\left\{z_1(L,\u{3}),z_2(L,\u{3})\right\}.}
Comparing ~\appbnew~and ~\pin~we see that the expressions in the rhs
of these expressions differ  only by a factor of $q^{1/2}.$ According
to ~\ppdis~with $\mu=3,~a=-1$ the flow $\ev{y_3-2}{-1}$ and the 
flow $\ev{y_3-2}{0}$ have identical steps starting with the step
$y_2+y_1+2.$ This fact together with ~\appbnew,~\pin~ and proposition
1 with $\mu=3,a=0$ imply that 
\eqn\appfinal{\left\{f_s(L,-{\bf E}_{1,3}^{(t)}+\u{3}),f_s(L,{\bf
e}_{1}-{\bf E}_{1,3}^{(t)}+\u{3})\right\}\ev{y_3-2}{-1}
q^{c(t_3)+{1\over 2}}\left\{P_{-1}(L,3,\u{3}),f_s(L,\u{3})\right\}}
which concludes the proof of proposition 1 with $a=-1,\m=3,~\nu_0\geq 2$.
The proof of the special case  $\nu_0=1$ follows along 
similar lines and will be
omitted.

\appendix{C}{Proof of 10.10 and 10.15}

We sketch the proof of ~\xxRR~ and ~\oxRR~ in several
stages. First we use the results of sec. 8 to 
prove (C.12) and (C.13). From theorem 1 of sec. 8 and (C.13) we then
prove ~\oxRR~for
the case $\beta=\alpha+1$ and after that we use 
both (C.12) and (C.13) to prove
the general case of ~\oxRR~with $1\leq\alpha\leq \beta-2$ by 
induction. We  conclude by proving~\xxRR.

It is useful to introduce the notation  ${\tilde j}_\mu=j_{\mu}-t_{\mu},~
1+t_\mu\leq j_{\mu}\leq t_{\mu+1}+\delta_{\mu,n}.$ 
In addition, we define the flow $\ev{b'-2}{a}$ with
$a=0,-1$ and $y_{\mu}<b'\leq y_{\mu+1}$ as part of the flow
$\ev{y_{\mu+1}-2}{a}$ restricted to the first $b'-2$ steps. 
We also note the equality of flows
\eqn\appa{\ev{{\tilde j}_{\mu}y_{\mu}}{a}=\ev{
{\tilde j}_{\mu}y_{\mu}-2}{a}\evs{2},~~2\leq{\tilde j}_{\mu}\leq
\nu_{\mu}+\delta_{\mu,n},~2-|a|\leq \mu \leq n}
and
\eqn\appb{\ev{y_{\mu}}{a}=\ev{y_{\mu}-2}{0}\evs{2},~~
2-|a|\leq \mu \leq n}
with  $a=0,-1.$
These follow from ~\rpropiii~for $a=0$ if we note   
 the
Takahashi decomposition 
\eqn\newapp{{\tilde
j}_{\mu}y_{\mu}=
\cases{l^{(\mu)}_{j_{\mu}}+l^{(\mu-1)}_{t_{\mu}}&for~$2\leq{\tilde j}_{\mu}\leq
\nu_{\mu}+\delta_{\mu,n},~1\leq \mu \leq n$\cr
l^{(\mu-1)}_{1+t_{\mu}}&for~${\tilde j}_{\mu}=1,~1\leq \mu \leq n.$\cr}}
The identical argument works for $a=-1$ if we replace $j_{\mu}$ by $1+j_{\mu}.$

Next we use ~\appa~and \appb~ along with
\pc -\pda~and lemma 2.2 of sec. 8 to obtain
\eqn\appc{\eqalign{\{P_0(L,\mu+1,{\bf u}'_{\mu+1})&,P_1(L,\mu+1,{\bf
u}'_{\mu+1})\}\ev{{\tilde j}_{\mu}y_{\mu}}{a}\cr
q^{c(j_\mu)-c(t_{\mu-1})-{\nu_0-\kie{\mu}\over
4}+{1\over 2}a(-1)^{\mu+1}}&
\{f_s(L,-\E{1}{\mu+1}+{\bf e}_{t_{1+\mu}-{\tilde j}_\mu -1}+{\bf
u}'_{\mu+1}),\cr
&f_s(L,{\bf e}_1-\E{1}{\mu+1}+{\bf e}_{t_{1+\mu}-
{\tilde j}_\mu -1}+{\bf u}'_{\mu+1})\}\cr
+q^{c(j_\mu)-c(t_{\mu-1})+{1\over 2}a(-1)^{\mu+1}}\{&P_0(L,\mu,{\bf e}_{t_{1+\mu}-{\tilde
j}_{\mu}}-{\bf e}_{t_{1+\mu}}+{\bf u}'_{\mu+1}),\cr
&P_1(L,\mu,{\bf e}_{t_{1+\mu}-{\tilde j}_{\mu}}-
{\bf e}_{t_{1+\mu}}+{\bf u}'_{\mu+1})\}\cr}}
with $1+|a|\leq {\tilde j}_{\mu}\leq \nu_{\mu}-2,~2\leq \mu\leq n-1,
~a=0,-1.$ Furthermore, from ~\ppdis~ and ~\pppropb~we derive
\eqn\appd{\ev{l^{(\mu)}_{1+j_{\mu}}-2}{a}=\ev{{\tilde
j}_{\mu}y_{\mu}}{a}\ev{y_{\mu-1}-2}{0}~~~~{\rm for}~2\leq \mu \leq
n,~1\leq {\tilde j}_{\mu}\leq \nu_{\mu}-1+2\delta_{\mu,n},~a=0,-1.}
Then combining~\appd~ with ~\appc, Prop. 1 from sec. 8 and ~
\pc~with $a=0$ and  $\mu$ replaced by $\mu-1$ gives
\eqn\appe{\eqalign{\{P_0(L,\mu+1,{\bf u}'_{\mu+1}),&P_1(L,\mu+1,{\bf
u}'_{\mu+1})\}\ev{l^{(\mu)}_{1+j_{\mu}}-2}{a}\cr
q^{c(j_{\mu})-{\nu_0-\kie{\mu}\over 4}+{1\over 2}a(-1)^{\mu+1}}
\{&Y(L,\mu,j_{\mu},{\bf
u}'_{\mu +1}),
f_s(L,-\E{\mu}{\mu+1}+{\bf e}_{t_{1+\mu}-{\tilde j}_{\mu}-1}+{\bf
u}'_{\mu+1})\cr
&~~~~~~+q^{-{(-1)^{\mu}\over 4}}
f_s(L,{\bf e}_{-1+t_{\mu}}-\E{\mu}{\mu+1}+
{\bf e}_{t_{1+\mu}-{\tilde j}_{\mu}}+{\bf u}'_{\mu+1})\}\cr}}
with $1+|a|\leq {\tilde j}_{\mu} \leq \nu_{\mu}-2,~~2\leq \mu\leq
n-1,~a=0,-1$ and
\eqn\appf{\eqalign{Y(L,\mu,j_{\mu},{\bf u}'_{\mu+1})=&
P_{-1}(L,\mu-1,-\E{\mu}{\mu+1}+{\bf e}_{t_{1+\mu}-{\tilde
j}_{\mu}-1}+{\bf u}'_{\mu+1})\cr
&+q^{(-1)^{\mu}\over 4}P_{-1}(L,\mu-1,{\bf e}_{-1+t_{\mu}}-\E{\mu}{\mu+1}+
{\bf e}_{t_{1+\mu}-{\tilde j}_{\mu}}+{\bf u}'_{\mu+1})\cr
&+q^{{L\over 2}-{\nu_0-\kie{\mu}\over 4}}f_s(L,-\E{1}{\mu-1}+
{\bf e}_{t_{1+\mu}-{\tilde j}_{\mu}}+{\bf u}'_{\mu+1}).\cr}}
The analogues of ~\appe~and ~\appf~for $\mu=1$ are 
\eqn\appg{\eqalign{&\{P_0(L,2,{\bf u}'_2),P_1(L,2,{\bf
u}'_2)\}\ev{l^{(1)}_{1+j_1}-2}{a}\cr
&q^{c(j_1)-{\nu_0\over 4}+{a\over 2}}\{Y(L,1,j_1,{\bf u}'_2),q^{{L\over
2}-{\nu_0\over 4}}f_s(L,{\bf e}_{t_2-{\tilde j}_1-1}-\E{1}{2}+{\bf
u}'_2)+\cr
&~~~~~~~~~~~~~~~f_s(L,{\bf e}_{-1+t_1}+{\bf e}_{t_2-{\tilde j}_1}-
\E{1}{2}+{\bf u}'_2)\}\cr}}
and
\eqn\apph{Y(L,1,j_1,{\bf u}'_{2})=f_s(L,{\bf e}_{t_2-{\tilde j}_1}
-{\bf e}_{t_2}+{\bf u}'_2)}
with $1+|a|\leq {\tilde j}_1\leq \nu_1-2,~a=0,-1.$
To derive~\appg-\apph~ for $a=0$  
we compare case 3 of ~\tzcc~with $j_0=0,1$ with case 2 of ~\tzcc~with
$j_{0}=\nu_0-1,$ and $j_1~{\rm replaced~by}~j_1-1$ and 
case 4 of ~\tzcc~to find the desired result
after  replacing ${\bf e}_{j_2}-\E{3}{n}$ by ${\bf u}'_2.$
By comparing (B.7)-(B.8) and ~\appg~with $j_1=2,a=0$ and recalling
that according to ~\ppdis~ the flows $\ev{y_2-2}{0}$ and
$\ev{y_2-2}{-1}$ have the same steps, starting with the step
$1+y_0+y_1$ we easily verify ~\appg~for $a=-1,~2\leq {\tilde
j}_{1}\leq \nu_1-2.$

Next we apply the flow $\ev{l^{(\mu-1)}_{2+j_{\mu-1}}-2}{-1}$ 
to the right hand side (rhs)
of ~\appc~with $a=0$. Taking into account ~\ppe~with 
$a=-1,$ 
$\mu$ replaced by $\mu-1$
and $j_{\mu}$ replaced by $j_{\mu-1}+1$ and ~\appe~with $a=-1,$ $\mu$ replaced
by $\mu-1,~\mu\geq 3$ or ~\appg~with $a=-1$
with $j_1$ replaced by $1+j_1$ for $\mu=2$ we derive
\eqn\appi{\eqalign{&{\rm rhs~of~\appc}\ev{l^{(\mu-1)}_{2+j_{\mu-1}}-2}{-1}
q^{c(j_\mu)+c(1+j_{\mu-1})-{\nu_0-\kio{\mu}\over 4}+{1\over 2}(-1)^{\mu+1}}\times\cr
&~~~\Bigl(\{Y(L,\mu-1,1+j_{\mu-1},{\bf e}_{t_{1+\mu}-{\tilde j}_{\mu}}-{\bf
e}_{t_{1+\mu}}+{\bf u}'_{\mu+1}),\cr
&~~~~q^{\delta_{\mu,2}({L\over 2}-{\nu_0\over 4})}f_s(L,-\E{\mu-1}{\mu+1}
+{\bf e}_{t_{\mu}-{\tilde j}_{\mu-1}-2}+{\bf e}_{t_{1+\mu}-{\tilde
j}_{\mu}}+{\bf u}'_{\mu+1})\cr
&~~~~+q^{{(-1)^{\mu}-\delta_{\mu,2}\over 4}}f_s(L,{\bf e}_{-1+t_{\mu-1}}
-\E{\mu-1}{\mu+1}+{\bf e}_{t_{\mu}-{\tilde j}_{\mu-1}-1}+
{\bf e}_{t_{1+\mu}-{\tilde j}_{\mu}}+{\bf u}'_{\mu+1})\}\cr
&~~~~+q^{(-1)^\mu\over 4}\{X(L,\mu-1,1+j_{\mu-1},{\bf
e}_{t_{1+\mu}-{\tilde j}_{\mu}-1}-{\bf e}_{t_{1+\mu}}
+{\bf u}'_{\mu+1}),\cr
&~~~~~~~~~~~~~~~f_s(L,-\E{\mu}{\mu+1}+{\bf e}_{1+j_{\mu-1}}+{\bf
e}_{t_{1+\mu}-{\tilde j}_{\mu}-1}+{\bf u}'_{\mu+1})\}\Bigr)\cr}}
with $1+t_{\mu-1}\leq j_{\mu-1}\leq t_{\mu}-3,~1+t_{\mu}\leq j_{\mu}
\leq t_{\mu+1}-2,~2\leq \mu \leq n-1$ and
\eqn\appj{\eqalign{X(L,\mu,j_{\mu},{\bf u}'_{\mu+1})=\theta(j_{\mu}>y_1)
&\Bigl(P_{-1}(L,\mu-\delta_{1+t_{\mu},j_{\mu}},{\bf e}_{j_{\mu}}-{\bf
e}_{t_{1+\mu}}+{\bf u}'_{\mu+1})+\cr
&q^{{L\over 2}-{\nu_0-\kie{\mu}\over 4}}f_s(L,{\bf
e}_{-1+j_{\mu}}-\E{1}{1+\mu}+{\bf u}'_{1+\mu})\Bigr)\cr
&+\delta_{\mu,1}\delta_{j_1,y_1}f_s(L,-{\bf e}_{t_2}+{\bf
u}'_2)\cr}}
for $1\leq \mu\leq n.$

To proceed further we note the following equality of flows:
\eqn\appk{\aw{{\tilde j}_{\mu}y_{\mu}}~~\ev{l^{(\mu-1)}_{2+j_{\mu-1}}-2}{-1}=
\aw{l^{(\mu)}_{1+j_{\mu}}+l^{(\mu-1)}_{1+j_{\mu-1}}-2}~~=~~
\aw{l^{(\mu)}_{1+j_{\mu}}-2}\evs{2}\aw{l^{(\mu-1)}_{1+j_{\mu-1}}-2}}
with $\mu\geq 2,~1\leq {\tilde j}_{\mu}\leq
\nu_{\mu}-2+2\delta_{\mu,n},
~1\leq {\tilde
j}_{\mu-1}\leq \nu_{\mu-1}-1.$ The first equation in ~\appk~
follows from ~\pppropb~with $a=0$. The second equation in ~\appk~
is just properties (3.12)-(3.13) of the map $b\rightarrow r$ with
$b=l^{(\mu)}_{1+j_{\mu}}.$ 
Eqn. ~\appk~together with ~\appe~with $a=0$ and \appi~imply the flow
\eqn\appl{\eqalign{&\{Y(L,\mu,j_{\mu},{\bf u}'_{\mu+1}),
f_s(L,-\E{\mu}{\mu+1}+{\bf e}_{t_{1+\mu}-
{\tilde j}_{\mu}-1}+{\bf u}'_{\mu+1})\cr
&~~~~~~~~~~~~~~~~~~~~+q^{-{(-1)^{\mu}\over 4}}f_s(L,{\bf
e}_{-1+t_{\mu}}-\E{\mu}{\mu+1}+{\bf e}_{t_{1+\mu}-{\tilde
j}_{\mu}}+{\bf u}'_{\mu+1})\}\cr
&\evs{2}\aw{l^{(\mu-1)}_{1+j_{\mu-1}}-2}q^{c(j_{\mu-1})-{\nu_0+1\over
4}+{3\over 4}\kio{\mu}}\times\cr
&\Big(\{Y(L,\mu-1,1+j_{\mu-1},{\bf e}_{t_{1+\mu}-{\tilde j}_{\mu}}-
{\bf e}_{t_{1+\mu}}+{\bf u}'_{\mu+1}),\cr
&q^{\delta_{\mu,2}({L\over 2}-{\nu_0\over 4})}f_s(L,-\E{\mu-1}{\mu+1}
+{\bf e}_{t_{\mu}-{\tilde j}_{\mu-1}-2}+{\bf
e}_{t_{1+\mu}-{\tilde j}_{\mu}}+{\bf u}'_{1+\mu})\}\cr
&+q^{{(-1)^{\mu}\over 4}-{\delta_{\mu,2}\over 4}}
f_s(L,{\bf e}_{-1+t_{\mu-1}}-\E{\mu-1}{\mu+1}+{\bf e}_{t_{\mu}-{\tilde
j}_{\mu-1}-1}+{\bf e}_{t_{1+\mu}-{\tilde j}_{\mu}}+{\bf
u}'_{\mu+1})\}\cr
&+q^{{(-1)^{\mu}\over 4}}\{X(L,\mu-1,1+j_{\mu-1},{\bf
e}_{t_{1+\mu}-{\tilde j}_{\mu}-1}-{\bf e}_{t_{1+\mu}}+{\bf
u}'_{\mu+1}),\cr
&~~~~~~~~~~f_s(L,-\E{\mu}{\mu+1}+{\bf e}_{1+{\tilde j}_{\mu-1}}+{\bf
e}_{t_{1+\mu}-{\tilde j}_{\mu}-1}+{\bf u}'_{\mu+1})\}\Bigr)\cr}}
with $1+t_{\mu-1}\leq j_{\mu-1}\leq t_{\mu}-3,~1+t_{\mu}\leq
j_{\mu}\leq t_{\mu+1}-2,~~2\leq \mu \leq n-1.$

Next we use the recursive properties of sec. 5 along with
~\ppe~and ~\appe~with $a=0,~\mu$ replaced by $\mu-1$ and ~\appg~with
$a=0$ for $\mu=2$ to obtain
\eqn\appm{\eqalign{&\{X(L,\mu,j_{\mu},{\bf u}'_{\mu+1}),f_s(L,{\bf
e}_{j_\mu}-{\bf e}_{t_{1+\mu}}+{\bf u}'_{1+\mu})\}\evs{2}
\aw{l^{(\mu-1)}_{1+j_{\mu-1}}-2}q^{c(j_{\mu-1})-{\nu_0+1\over
4}+{3\over 4}\kio{\mu}}\times\cr
&\Bigl(\{X(L,\mu-1,j_{\mu-1},{\bf e}_{1+j_{\mu}}-{\bf
e}_{t_{1+\mu}}+{\bf u}'_{\mu+1}),f_s(L,-\E{\mu}{\mu+1}+
{\bf e}_{j_{\mu-1}}+{\bf e}_{1+j_{\mu}}+{\bf u}'_{1+\mu})\}\cr
&+q^{(-1)^{\mu}\over 4}\{Y(L,\mu-1,j_{\mu-1},{\bf e}_{j_\mu}-{\bf
e}_{t_{1+\mu}}+{\bf u}'_{\mu+1}),\cr
&~~~~~~~~~~q^{\delta_{\mu,2}({L\over 2}-{\nu_0\over 4})}f_s(L,-\E{\mu-1}{\mu+1}
+{\bf e}_{t_\mu-{\tilde j}_{\mu-1}-1}+{\bf e}_{j_{\mu}}+{\bf
u}'_{\mu+1})\cr
&~~~~~~~~~~+q^{{(-1)^{\mu}\over 4}-{\delta_{\mu,2}\over 4}}f_s(L,{\bf
e}_{-1+t_{\mu-1}}-\E{\mu-1}{\mu+1}+{\bf e}_{t_\mu-{\tilde j}_{\mu-1}}
+{\bf e}_{j_\mu}+{\bf u}'_{\mu+1})\}\Bigr)\cr}}
with $1+t_{\mu-1}\leq j_{\mu-1}\leq t_{\mu}-2,~1+t_{\mu}\leq
j_{\mu}\leq t_{1+\mu}-1+\delta_{\mu,n},~2\leq \mu\leq n.$ Theorem 1
(of sec. 8) with $P=0$ and ~\appm~with ${\bf u}'_{\mu+1}$ specified to
be $-\E{\mu+2}{n}$ along with
\eqn\appn{\aw{l^{(\mu)}_{1+j_{\mu}}+l^{(\mu-1)}_{1+j_{\mu-1}}-2}=
\aw{l^{(\mu)}_{1+j_{\mu}}-2}\evs{2}\aw{l^{(\mu-1)}_{1+j_{\mu-1}}-2}}
imply
\eqn\appo{\eqalign{&\{f_s(L,-\E{1}{n}),f_s(L,{\bf e}_1-\E{1}{n})\}
\aw{l^{(\mu)}_{1+j_{\mu}}+l^{(\mu-1)}_{1+j_{\mu-1}}-2}
q^{c(j_\mu)+c(j_{\mu-1})-{\nu_0+1\over 4}+{3\over
4}\kio{\mu}}\times\cr
&\Bigl(\{X(L,\mu-1,j_{\mu-1},{\bf
e}_{1+j_{\mu}}-\E{1+\mu}{n}),f_s(L,{\bf e}_{j_{\mu-1}}+{\bf e}_{1+j_\mu}
-\E{\mu}{n})\}\cr
&+q^{(-1)^{\mu}\over 4}\{Y(L,\mu-1,j_{\mu-1},{\bf
e}_{j_{\mu}}-\E{1+\mu}{n}),\cr
&~~~~~~~~~q^{\delta_{\mu,2}({L\over 2}-{\nu_0\over 4})}
f_s(L,{\bf e}_{t_{\mu}-{\tilde j}_{\mu-1}-1}+{\bf e}_{j_{\mu}}-\E{\mu-1}{n})\cr
&~~~~~~~~~+q^{{(-1)^{\mu}\over 4}-{\delta_{\mu,2}\over 4}}
f_s(L,{\bf e}_{-1+t_{\mu-1}}+{\bf e}_{t_\mu-{\tilde
j}_{\mu-1}}+{\bf e}_{j_{\mu}}-\E{\mu-1}{n})\}\Bigr)\cr}}
with $1+t_{\mu-1}\leq j_{\mu-1} \leq t_{\mu}-2,~~1+t_{\mu} \leq
j_{\mu}\leq t_{\mu+1}+\delta_{\mu,n}-1,~~2 \leq \mu \leq n.$ The
second terms of the right hand side of ~\appo~ agrees with 
the right hand side of 
~\oxRR~when $\beta=\alpha+1=\mu$ which thus completes the proof of
this special case.

To complete the proof of ~\oxRR~in the general case $1\leq \alpha \leq
\beta -2$ we use 
the equality of flows
\eqn\appp{\aw{b(\alpha+1,\beta)-2}\evs{2}
\aw{l^{(\alpha)}_{1+j_{\alpha}}-2}=\aw{b(\alpha,\beta)-2}}
where
\eqn\appq{\eqalign{b(\alpha,\beta)&
=\sum_{\mu=\alpha}^{\beta}l^{(\mu)}_{1+j_{\mu}};~~~1\leq
\alpha,~1+\alpha\leq \beta \leq n\cr
&1+t_{\mu}\leq j_{\mu} \leq t_{\mu+1}-3,~~\alpha \leq \mu \leq \beta-2\cr
&1+t_{\beta-1} \leq j_{\beta-1} \leq t_{\beta}-2\cr
&1+t_{\beta} \leq j_{\beta} \leq t_{\beta+1}-1+\delta_{\beta , n}.\cr}}
Eqn.\appp~is just properties (3.12)-(3.13) of the $b\rightarrow r$ map
with $b=b(\alpha+1,\beta)$. 
Then by starting with ~\appo~and by use of~\appl~and \appm~
we may prove by induction $(\alpha+1\rightarrow \alpha)$ that
\eqn\appr{\eqalign{&\{f_s(L,-\E{1}{n}), f_s(L,{\bf e}_1-\E{1}{n})\}
\aw{b(\alpha,\beta)-2}\cr
&\sum_{i_{\alpha+1},\cdots,i_{\beta-1}=0,1\atop
i_{\alpha}=0}q^{c_{(4)}({\bf j})+rf_{(4,1)}({\bf
i})}\{X(L,\alpha,j_{\alpha}+i_{\alpha+1},{\bf u}_{(4)}({\bf i}, {\bf j})
+{\bf e}_{t_{1+\alpha}}),f_s(L,{\bf u}_{(4,1)}({\bf i},{\bf
j}))\}\cr
&+\sum_{i_{\alpha+1},\cdots,i_{\beta-1}=0,1\atop i_{\alpha}=1}
q^{c_{(4)}({\bf j})+rf_{(4,2)}({\bf i})}
\{Y(L,\alpha,j_{\alpha}+i_{\alpha+1},{\bf u}_{(4)}({\bf i},{\bf
j})+{\bf e}_{t_1+\alpha}),\cr
&~~~q^{\delta_{1,\alpha}({L\over 2}-{\nu_0\over 4})}
f_s(L,{\bf u}_{(4,2)}({\bf i},{\bf j}))
+q^{-{(-1)^{\alpha}+\delta_{1,\alpha}\over 4}}
f_s(L,{\bf u}_{(4,3)}({\bf i},{\bf j}))\}\cr}}
where $c_{(4)}({\bf j}),~rf_{(4,1)}({\bf i}),~rf_{(4,2)}({\bf i}),~{\bf
u}_{(4)}({\bf i},{\bf j}),
{\bf u}_{(4,1)}({\bf i},{\bf j}),
{\bf u}_{(4,2)}({\bf i},{\bf j}).~{\rm and}~
{\bf u}_{(4,3)}({\bf i},{\bf j})$ are defined in sec. 10.3. The second
terms on the righthandside of ~\appr~ agree with the right hand side of
~\oxRR~which thus completes the proof.

It remains to prove~\xxRR. To do this we first compare cases 1 and 2
of ~\tzcb~ with $\E{3}{n}$ replaced by ${\bf u}'_{2}$ to infer that
\eqn\apps{\eqalign{&\{X(L,1,j_1,{\bf u}'_2),f_s(L,{\bf e}_{j_1}-{\bf e}_{t_2}+{\bf
u}'_{2})\}\evs{2}\aw{l^{(0)}_{1+j_0}-2}\cr
&\{q^{-{\nu_0-2\over 4}}f_s(L,{\bf e}_{j_0-1}+{\bf
e}_{j_1+1}-\E{1}{2}+{\bf u}'_2)+{\tilde f}_s(L,{\bf
e}_{\nu_0-j_0}+{\bf e}_{j_1}-\E{1}{2}+{\bf u}'_{2}),\cr
&q^{-{\nu_0-2\over 4}}f_s(L,{\bf e}_{j_0}+{\bf e}_{j_1+1}-\E{1}{2}+{\bf
u}'_{2})+{\tilde f}_s(L,{\bf e}_{\nu_0-j_0-1}+{\bf
e}_{j_1}-\E{1}{2}+{\bf u}'_{2})\}\cr}}
with $1\leq j_0\leq \nu_0-1,~1+t_{1}\leq j_1 \leq t_2-1.$ Furthermore we have
\eqn\appt{\eqalign{&\{X(L,1,j_1,{\bf u}'_2),f_s(L,{\bf e}_{j_1}
-{\bf e}_{t_2}+{\bf u}'_{2})\}\ev{1}{\pprbb~b}\cr
&\{f_s(L,{\bf e}_{j_1}-{\bf e}_{t_2}+{\bf u}'_{2}),q^{-{\nu_0-2\over 4}}
f_s(L,{\bf e}_{j_1+1}-\E{1}{2}+{\bf u}'_2)+{\tilde f}_s(L,{\bf
e}_{\nu_0-1}+{\bf e}_{j_1}-\E{1}{2}+{\bf u}'_2)\}\cr}}
with $1+t_1\leq j_1 \leq t_2-1.$ Analogously by comparing 
cases 2 and 4 of ~\tzcc~with ${\bf
e}_{j_2}-\E{3}{n}$ replaced by ${\bf u}'_2$ we find
\eqn\appu{\eqalign{&\{Y(L,1,j_1,{\bf u}'_2),q^{{L\over 2}-{\nu_0\over
4}}
f_s(L,{\bf e}_{t_2-{\tilde j}_1-1}-\E{1}{2}+{\bf u}'_2)\cr
&~~~~~~~~~~~~~~~~~~~~~+f_s(L,{\bf e}_{-1+t_1}+{\bf e}_{t_2-{\tilde
j}_1}-\E{1}{2}+{\bf u}'_2)\}\evs{2}\aw{l^{(0)}_{1+j_0}-2}\cr
&\{q^{-{\nu_0-2\over 4}}f_s(L,{\bf e}_{j_0}+{\bf e}_{t_2-{\tilde
j}_1-1}-\E{1}{2}+{\bf u}'_2)+q^{1\over 2}{\tilde f}_s(L,{\bf
e}_{\nu_0-j_0-1}+{\bf e}_{ t_2-{\tilde j}_1}+{\bf u}'_2),\cr
&q^{-{\nu_0-2\over 4}}f_s(L,{\bf e}_{j_0+1}+{\bf e}_{t_2-{\tilde
j}_1-1}-\E{1}{2}+{\bf u}'_2)+q^{1\over 2}{\tilde f}_s(L,{\bf
e}_{\nu_0-{j}_0-2}+{\bf e}_{t_2-{\tilde j}_1}-\E{1}{2}+{\bf
u}'_{2})\}\cr}}
with $1\leq j_0\leq \nu_0-1, 1+t_1\leq j_1 \leq t_2-2$ and
\eqn\appv{\eqalign{&\{Y(L,1,j_1,{\bf u}'_2),q^{{L\over 2}-{\nu_0\over 4}}
f_s(L,{\bf e}_{t_2-{\tilde j}_1-1}-\E{1}{2}+{\bf u}'_2)\cr
&~~~~~~~~~~~~~~~~~~~~~~~~~~~+
f_s(L,{\bf e}_{-1+t_1}+{\bf e}_{t_2-{\tilde j}_1}-\E{1}{2}+{\bf
u}'_2)\}\ev{1}{\pprbb~b}\cr
&\{q^{{L\over 2}-{\nu_0\over 4}}f_s(L,{\bf e}_{t_2-{\tilde j}_1-1}
-\E{1}{2}+{\bf u}'_2)+f_s(L,{\bf e}_{-1+t_1}+{\bf e}_{t_2-{\tilde
j}_1}-\E{1}{2}+{\bf u}'_2),\cr
&q^{-{\nu_0-2\over 4}}f_s(L,{\bf e}_1+{\bf e}_{t_2-{\tilde
j}_1-1}-\E{1}{2}+{\bf u}'_2)+q^{1\over 2}{\tilde f}_s(L,{\bf
e}_{\nu_0-2}+{\bf e}_{t_2-{\tilde j}_1}-\E{1}{2}+{\bf u}'_2)\}\cr}}
with $1+t_1\leq j_1 \leq t_2-2.$
The result ~\xxRR~ follows from combining ~\appr~with $\alpha=1$ with
~\apps -\appv.

\vfill
\eject
\listrefs

\vfill\eject

\bye
\end